\documentclass{article} 
\usepackage{pdfpages}
\usepackage{geometry} 
\usepackage{amsmath}  
\usepackage{graphicx} 
\usepackage{float}
\usepackage{color}
\usepackage{amsmath}
\usepackage{amssymb}
\usepackage{amsthm}

\newtheorem{theorem}{Theorem}

\theoremstyle{definition}
\newtheorem{definition}{Definition}[section]

\theoremstyle{lemma}

\theoremstyle{remark}

\theoremstyle{boldremark}
\newtheorem*{boldremark}{Remark} 

\newtheorem{corollary}[theorem]{Corollary}

\newcommand{\maketitletwo}[2][]{\begin{center}
        \Large{\textbf{A General Solution to Bellman’s Lost-in-a-forest Problem}
            
            } 
        \vspace{5pt}
        
        \normalsize{Zhipeng Deng

        }        
        \vspace{15pt}
        
\end{center}}
\begin{document}
    \maketitletwo[5]

    \begin{abstract}
We present a general solution and formulation framework to Bellman’s lost-in-a-forest problem. The forest boundary is known and may take any shape. The starting point and the orientation are unspecified. We convert the problem into translation and rotation of the forest boundary. This transformation allows us to formulate this problem as a constrained minimization problem. Upon discretization, the problem becomes a variation of the traveling salesman problem or the Hamiltonian path problem. We leverage discrete optimization and derive several nontrivial results consistent with those from previous papers. This method is general, and we also extend the approach to related problems, including Moser's worm problem and the shortest opaque set problem.
\end{abstract}

\\

\noindent \textbf{Classification}

Optimization and Control (math.OC)

Metric Geometry (math.MG)

Discrete Mathematics (cs.DM)

Computational Geometry (cs.CG)

49K30 (Optimal Solutions in Calculus of Variations)

49Q10 (Optimization of Shapes Other Than Minimal Surfaces)

52A40 (Inequalities and Extremum Problems)

\tableofcontents

\section{Introduction}

Bellman's lost-in-a-forest problem is an unsolved minimization problem in geometry, introduced in 1955 by Richard E. Bellman \cite{Gross1955} \cite{forest wiki} \cite{Croft2012}. The problem is commonly stated as follows: 
\\

``A hiker is lost in a forest whose shape and dimensions are precisely known. What is the best path for the hiker to follow to escape the forest?"
\\

It is typically assumed that the hiker does not know the starting point or the facing orientation. The optimal path is defined as the one that minimizes the worst-case distance the hiker must travel to reach the forest boundary. There are also other ways to express this problem. In earlier studies, the problem was treated as a minimax problem \cite{Gluss1961}. The reason why this problem is very difficult is that for any shape, we need to consider all the possible starting points and orientation, and the space for variation is very large.

While practical non-contrived real-world applications are not immediately evident, this problem belongs to a broader class of geometric optimization problems, including search strategies with practical significance \cite{Baezayates1993}. A bigger motivation for study has been the connection to Moser's worm problem. It was featured in a list of 12 problems described by Scott W. Williams as ``million buck problems" because he believed that resolving them would yield techniques worth at least a million dollars to mathematics \cite{Williams2000}.

A proven solution exists for only a limited number of shapes. For example, Isbell provided solutions for a straight line/half-plane \cite{Gross1955} \cite{Gluss1961} \cite{Isbell1957} \cite{Joris1980} \cite{Finch2019}, while Zalgaller solved the unit stripe case \cite{Zalgaller2005}. Solutions for a square and a circle have also been established \cite{Finch2004}. It has been conjectured that for convex regular polygons, the diagonal is the optimal path \cite{Ward2008}. However, for an equilateral triangle, Besicovitch's zigzag path was found to be optimal \cite{Besicovitch1965}. Bellman himself also proposed zigzag strategies involving straight-line segments. John W. Ward \cite{Ward2008} and Steven R. Finch \cite{Finch2004} have published detailed papers describing the problem and progress in addressing it. However, most approaches rely on geometry-based methods. A general solution in the form of a geometric algorithm, which could take any forest shape as input and return the optimal escape path, remains highly challenging. Additionally, some studies have applied numerical methods to specific cases, such as isosceles triangles \cite{Gibbs2016}, and approximation techniques \cite{Kübel2021}.

After literature review, we have not found any prior work employing continuous or discrete optimization or variational methods to solve this problem.

\section{Proof of concept}

Figures 1 and 2 demonstrate the proof of concept presented in this paper. The key to our success lies in transforming the escape path from the forest boundary—across various orientations and locations—into a path that does not move by the translation and rotation of the forest boundary. Previous researchers did not recognize the inspiration for this idea. In this way, we circumvent the need to rely exclusively on boundary shapes and geometric relationships, allowing a broader optimization-based perspective. With this perspective, we establish that for any translation or rotation of the forest boundary, the escape path must pass through an escape point on the boundary. Consequently, we reformulate the original problem into finding the shortest path that connects all possible escape points.

\begin{boldremark}
For each escape point and orientation, we only need to consider the first escape point. This means once escaping the forest boundary, subsequent paths are no longer considered and necessary.
\end{boldremark}

\begin{figure}[H]
  \centering
  \includegraphics[width=15cm]{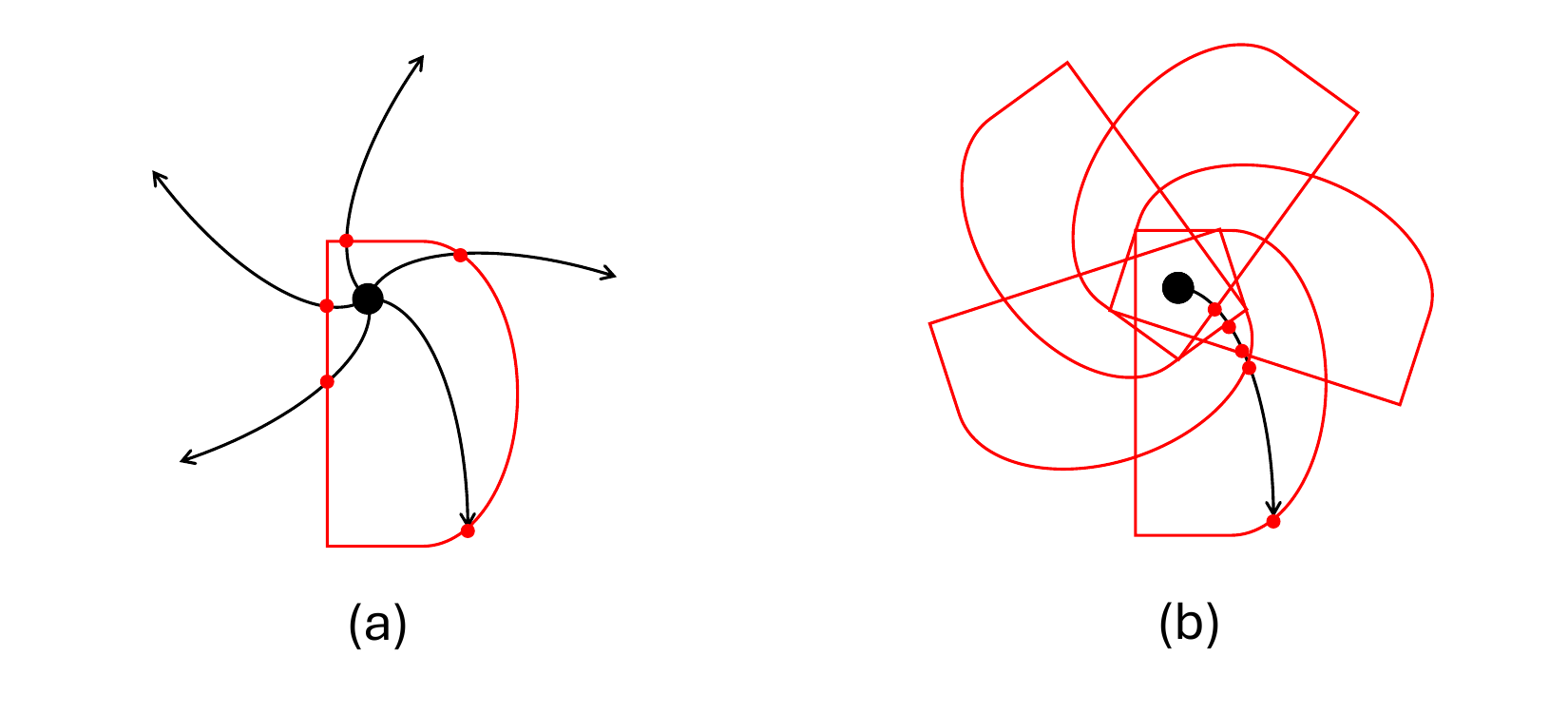}
  \caption{Starting from a known point with unknown orientation. We convert the escape path out of the forest boundary in different orientations to the rotation of the forest boundary. Therefore, the rotated escape points are all on the escape path. (Black curve represents the escape path, black dot is the starting point, red dots are the escape points, and red curve is the forest boundary)}
\end{figure}

\begin{figure}[H]
  \centering
  \includegraphics[width=15cm]{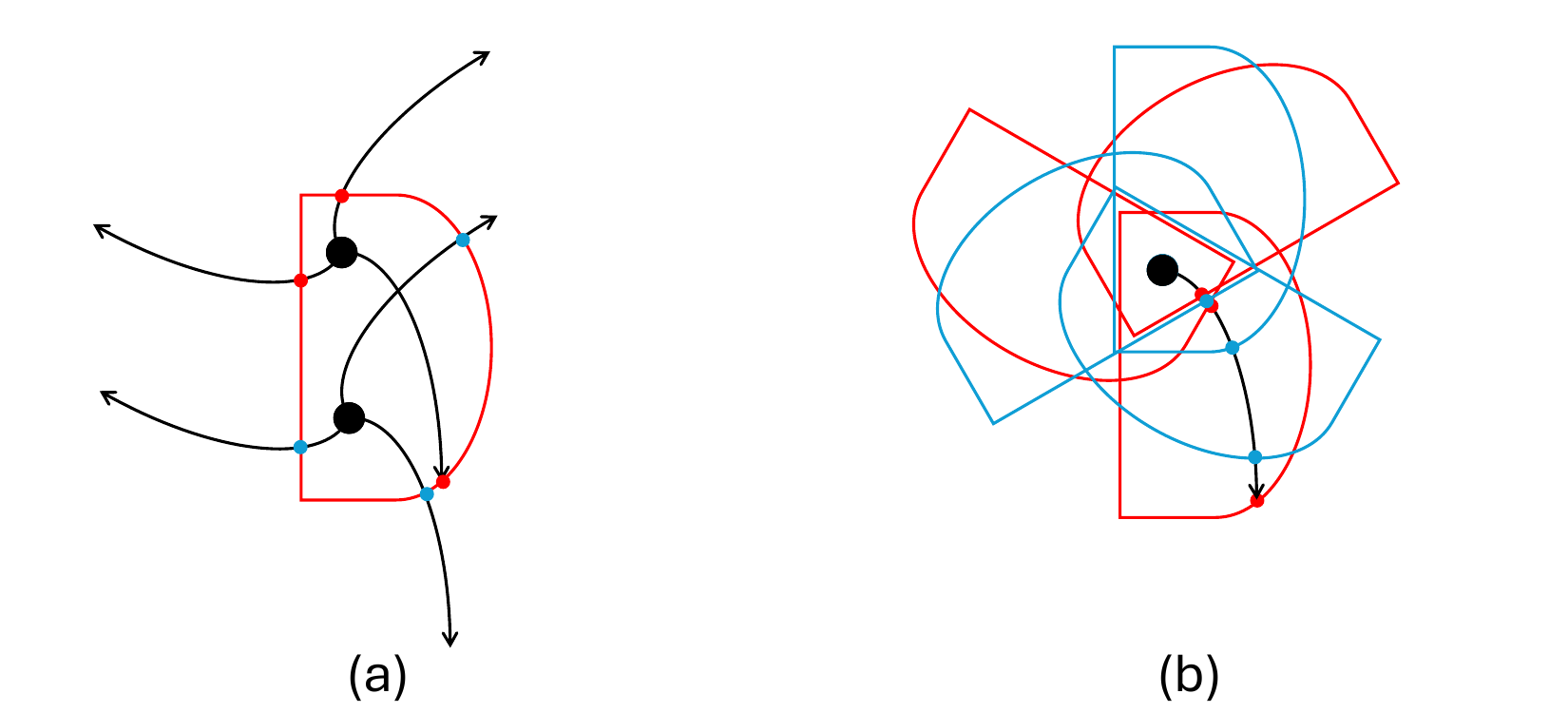}
  \caption{Starting from two possible known points with unknown orientation. We convert the path out of the forest boundary in different orientations and locations to the translation and rotation of the forest boundary. The rotated and translated escape points are all on the path. (Black curve represents the escape path, black dot is the starting point, red/blue dots are the escape points, and red/blue curve is the forest boundary)}
\end{figure}

\begin{definition}[Translation and rotation of forest boundary]
In this context, rotation involves rotating all forest boundaries such that the initial heading orientation of the escape path remains consistent. Similarly, translation involves shifting all forest boundaries so that the initial starting point of the escape path aligns consistently across cases.
\end{definition}

\begin{boldremark}
Consider a special Euclidean group SE(2) with rigid motions with combinations of translations and rotations applied to the boundary, the existence and uniqueness of the corresponding escape points are preserved due to the invariance of the Euclidean metric. 
\end{boldremark}

\section{Notations and definitions}
\subsection{Notations}
\begin{flushleft}
$F$: A given function representing the ``forest boundary" \\
$f$: A continuous trajectory/path function to escape the boundary \\
$(x,y)$: point in two-dimensional planar coordinate system \\
$R, R_{worm},R_O$: A region \\
$(s_x, s_y), (s_{kx}, s_{ky})$: Starting point \\
$\alpha$: Initial heading orientation \\
$\theta, \varphi $: angle \\
$L, L_N, L_{MN}$: Curve length \\
$i, j, k, h, u, v$: Positive integer index \\
$N, M, m$: Positive integer \\
$F_i, F_{ki}$: Rotated, rotated and translated forest boundary \\
$\textbf{a}$: Order of the escape points that escape path going through \\
$S_N, S_{MN}, S_H$: Symmetric group of all permutations \\
$s_{i,j}, b_{i,j}$: Binary decision variables \\
$u_i$: Auxiliary variables \\
$\phi_1, \phi_2, \phi_3, \delta_1, \delta_2, \delta_3, r$: Given parameters for lines and circle sector\\
$A$: Area of a region\\
$A_m$: Minimum area of a shape that can cover every unit-length curve\\
$P_i$: Partition of ${1,2,...,MN}$\\
\end{flushleft}

\subsection{Definitions for Bellman’s Lost-in-a-forest Problem}
In the Introduction Section, we provided a natural language description of the problem. Here, we define its mathematical formulation precisely below. 

We always assume the two-dimensional planar xy coordinate system.

\begin{definition}[Forest boundary]
Consider a two-dimensional plane $\mathbb{R}^2$ and a curve defined implicitly by an equation:

\begin{equation}
F(x, y) = 0
\end{equation}

where $F: \mathbb{R}^2 \to \mathbb{R}$ is a given function. The curve $F(x, y) = 0$ represents the ``forest boundary" aiming to find or escape.
\end{definition}

\begin{boldremark}
The forest boundary equation $F(x, y)$ can be not closed, degenerated, discontinuous, non-differentiable, non-convex, as we will only use in the minimization constraints. If the boundary is degenerated and discontinuous, the interpretation of ``escape" may not make sense. But search makes more sense. In this paper, we still use escape. 
\end{boldremark}

\begin{definition}[Escape path]
We aim to determine a continuous trajectory or path function without abrupt jumps or breaks as

\begin{equation}
f\left [ x(t),y(t) \right ] 
\end{equation}

such that for any possible unknown initial position and orientation, the path will intersect the boundary and there exists an escape point.

\begin{boldremark}
We need to assume the escape path exists, and the path must be integrable as length is finite and well-defined for integrability and regularity conditions. The escape path is an open simple curve (non-self-intersecting) and piecewise $C^1$.
\end{boldremark}
\end{definition}

\begin{definition}[Starting region]
We define a two-dimensional region $R$ which is the union of the starting points.
\end{definition}

\begin{definition}[Starting point]
More formally, let the starting region consist of an arbitrary unknown starting point:

\begin{equation}
(s_x, s_y) \in R \in\mathbb{R}^2
\end{equation}
\end{definition}

\begin{boldremark}
The start points and region $(s_x, s_y)\in R$ are not necessarily located inside the boundary defined by $F(x, y) < 0$. The original natural language description seems to be required, but it is not necessary.
\end{boldremark}

\begin{definition}[Orientation]
We define an arbitrary unknown initial heading orientation at an angle to the x-axis:
\begin{equation}
\alpha \in [0, 2\pi]
\end{equation}
\end{definition}

\begin{definition}[Length of escape path]
The curve length $L$ of the function $f$ is defined as:

\begin{equation}
L = \int \sqrt{x'^2 + y'^2} \, dt
\end{equation}

\end{definition}

\begin{definition}[Bellman’s lost-in-a-forest problem]
The requirement of the problem is to escape or reach the boundary is for $f$ to satisfy the condition that, regardless of the initial state $(s_x, s_y, \alpha)$, the path $f$ will eventually intersect the boundary $F(x, y) = 0$. The goal is to determine the minimal length of $f$.

The formal solution to the problem is as follows:
\begin{equation}
\begin{split}
\text{minimize } &L\\
\text{subject to:} &F[x_\alpha(t), y_\alpha(t)] = 0, \forall (s_x, s_y) \in R, \alpha \in [0, 2\pi], 
\end{split}
\end{equation}

The escape points are $[x_\alpha(t), y_\alpha(t)]$, and all of these points lie on the path.

In the continuous case, the objective is to identify the functions that minimize the curve length functional $L$.
\end{definition}

\subsection{Definitions for Weak Form I}
The original problem is very challenging because of the arbitrary unknown starting point and orientation. Hence, we now consider two relaxed variants, refer to Weak Forms I and II of Bellman’s lost-in-a-forest problem.

\begin{definition}[Weak Form I of Bellman’s lost-in-a-forest problem]
As illustrated in Figure 1, Weak Form I of Bellman’s lost-in-a-forest problem assumes a known starting point $(s_x, s_y)$, while the orientation $\alpha \in [0, 2\pi]$ is unknown. 
\end{definition}

\begin{definition}[A natural number $N$ for different orientation]
We define $N\in \mathbb{N}$ to be a natural number. Let the orientation uniformly distribute across $N$ possible choices from $0$ to $2\pi$.
\end{definition}

\begin{definition}[Escape orientation in Weak Form I]
With \((s_x, s_y) \in R\) being a given known starting point. Let
\begin{equation}
\alpha \in \left\{\frac{2\pi i}{N} \, \middle| \, i = 0, 1, 2, \ldots, N-1\right\}
\end{equation}
represent unknown initial heading orientation with \(N\) possible choices. Note that angle increments $2\pi/N$ approaching 0 for infinite $N$.
\end{definition}

\begin{definition}[Rotated forest boundary]
With the rotation matrix, rotated forest boundary $F_i$ is defined as:
\begin{equation}
\begin{split}
&F_i(x, y) = \\
&F \left[s_x + (x - s_x) \cos \frac{2\pi i}{N} - (y - s_y) \sin \frac{2\pi i}{N}, s_y + (x - s_x) \sin \frac{2\pi i}{N} + (y - s_y) \cos \frac{2\pi i}{N}\right]=0,\\
&\forall i \in \{0, 1, 2, \ldots, N-1\}
\end{split}
\end{equation}
\end{definition}

\begin{definition}[Escape points in Weak Form I]
The escape points $\left ( x_i,y_i \right ) \in \mathbb{R}^2$ on each rotated forest boundary are denoted as 
\begin{equation*}
(x_0,y_0),(x_1,y_1),...,(x_{N-1},y_{N-1})
\end{equation*}
\end{definition}

Since each escape point lies on the rotated forest boundary, the following condition holds:
\begin{equation}
F_i(x_i, y_i) = 0, \forall i \in \{0, 1, 2, \ldots, N-1\},
\end{equation}

\begin{definition}[Order of escape points on escape path]
We only know escape points are on the escape path. However, the order is unknown. 

Let:
\begin{equation}
\mathbf{a} = \{a_0, a_1, \ldots, a_{N-1}\}, \mathbf{a} \in S_N,
\end{equation}
where $S_N$ is the symmetric group of all permutations of $\{0, 1, 2, \ldots, N-1\}$.

$\mathbf{a}$ is the order in which the escape path visits the escape points.
\end{definition}

\begin{definition}[Length of escape path for Weak Form I]
Since the shortest curve between two points is a line segment. We define the length of escape path for Weak Form I  $L_N$ computed as the sum of the Euclidean distances from the origin to the first escape point and between the subsequent $N$ escape points in order.
\end{definition}
$L_N$ can be written as
\begin{equation}
L_N=\sqrt{(x_{a_0} - 0)^2 + (y_{a_0} - 0)^2} + \sum_{i=1}^{N-1} \sqrt{(x_{a_i} - x_{a_{i-1}})^2 + (y_{a_i} - y_{a_{i-1}})^2},\\
\end{equation}

\begin{definition}[Solution to Weak Form I of Bellman’s lost-in-a-forest problem]
The objective is to find the minimal length $L_N$ of the path $f$ for $N$ orientations to escape/find the boundary:
\begin{equation}
\begin{split}
\text{minimize } &L_N,\\
\text{subject to:} &F_i(x_i, y_i) = 0, \forall \alpha \in \left\{ \frac{2\pi i}{N} \, \bigg| \, i = 0, 1, 2, \ldots, N-1 \right\}.
\end{split}
\end{equation}

In the discrete case, the goal is to find the location and order of points $(x_i, y_i)$ that minimize the curve length $L_N$. 
\end{definition}

\subsection{Definitions for Weak Form II}
\begin{definition}[Weak Form II of Bellman’s lost-in-a-forest problem]
In Weak Form II, it is the case starting from several known points, but the orientation remains unknown. 

Same as in Weak Form I, the orientation is still from $0$ to $2\pi$ with $N$ uniformly possible choices. Let $\alpha \in \left\{ \frac{2\pi i}{N} \, \bigg| \, i = 0, 1, 2, \ldots, N-1 \right\}$ represent unknown initial heading orientation with $N$ possible choices.

\begin{definition}[A natural number $M$ for different starting points]
We define $M \in \mathbb{N}$ to be a natural number. With start points $(s_{kx}, s_{ky}) \in R \in \mathbb{R}^2$, $k = 1, 2, \ldots, M$, be given $M$ known starting points. 
\end{definition}

This results in a total of $MN$ escape points $(x_{ki}, y_{ki})$ for the combination of $M$ starting points and $N$ orientations, where $k = 1, 2, \ldots, M$, $i = 0, 1, 2, \ldots, N-1$.

\begin{boldremark}
The start points $(s_{kx}, s_{ky}) \in \mathbb{R}^2$, $k = 1, 2, \ldots, M$, are not necessarily located inside the boundary defined by $F(x, y) < 0$.
\end{boldremark}

\begin{definition}[Rotated and translated forest boundary]
With the rotation matrix, rotated and translated forest boundary $F_{ki}$ is defined as:
\begin{equation}
\begin{split}
&F_{ki}(x, y) =\\
&F \left[s_{kx} + (x - s_{kx}) \cos \frac{2\pi i}{N} - (y - s_{ky}) \sin \frac{2\pi i}{N}, s_{ky} + (x - s_{kx}) \sin \frac{2\pi i}{N} + (y - s_{ky}) \cos \frac{2\pi i}{N}\right]=0
\end{split}
\end{equation}
\end{definition}

\begin{definition}[Escape points in Weak Form II]
The escape points on each rotated forest boundary are denoted as 
\begin{equation*}
(x_0,y_0),(x_1,y_1),...,(x_{MN-1},y_{MN-1})
\end{equation*}
\end{definition}

Since each escape point lies on the rotated and translated forest boundary, the following condition applies:
\begin{equation}
F_{ki}(x_{ki}, y_{ki}) = 0, \forall i \in \{0, 1, 2, \ldots, N-1\}, \, k \in \{1, 2, \ldots, M\},
\end{equation}
where $F_{ki}$ is the rotated and translated forest boundary.

Again, we only know escape points are on the escape path. However, the order is unknown. 
\begin{definition}[Order of escape points on escape path in Weak Form II]
Recall Definition 3.3, let
\[
\mathbf{a} = \{a_1, a_2, \ldots, a_{MN}\} \in S_{MN},
\]
where $S_{MN}$ is the symmetric group of all permutations of $\{1, 2, \ldots, MN\}$.

$\mathbf{a}$ is the order of the escape path going through the combinations of $M$ different starting points and $N$ orientations.
\end{definition}

\begin{definition}[Length of escape path for Weak Form I]
Similarly, we define the length of escape path for Weak Form II  $L_{MN}$ computed as the sum of the Euclidean distances from the origin to the first escape point and between the subsequent $MN-1$ escape points in order.

$L_{MN}$ can be written as
\begin{equation}
L_{MN}=\sqrt{(x_{a_1} - 0)^2 + (y_{a_1} - 0)^2} + \sum_{i=2}^{MN} \sqrt{(x_{a_i} - x_{a_{i-1}})^2 + (y_{a_i} - y_{a_{i-1}})^2}
\end{equation}
\end{definition}

\begin{definition}[Solution to Weak Form II of Bellman’s lost-in-a-forest problem]
The objective is to find the minimal length $L_{MN}$ of $f$ for $M$ starting points and $N$ orientations to escape/find the boundary:
\begin{equation}
\begin{split}
\text{minimize } &L_{MN},\\
\text{subject to:} &F_{ki}(x_{ki}, y_{ki}) = 0, \forall \alpha \in \left\{ \frac{2\pi i}{N} \, \bigg| \, i = 0, 1, 2, \ldots, N-1 \right\}, \forall (s_{kx}, s_{ky}) \in \mathbb{R}^2, \forall k \in \{1, 2, \ldots, M\}.
\end{split}
\end{equation}

\end{definition}

In the following sections, we will present solutions for Weak Forms I and II, followed by the solution to the original Bellman’s lost-in-a-forest problem.
\end{definition}

\section{Weak Form I - Starting from a known point with unknown orientations}
\subsection{Discrete Formulation}

Without loss of generality, we assume the known starting point is at the origin $(0, 0)$.

Recall Definition 3.7 and Figure 1, the solution involves determining the shortest path passing through the escape points starting from the origin $(0, 0)$. But the location and order of the points is unknown and needs to be solved.

\begin{definition}[Using all permutations for Weak Form I]

The shortest path length is calculated as the Euclidean distance from the origin $(0, 0)$ and these points in order.
 
The solution to Weak Form I can be written as:
\begin{equation}
\begin{split}
\text{minimize } &\sqrt{(x_{a_0} - 0)^2 + (y_{a_0} - 0)^2} + \sum_{i=1}^{N-1} \sqrt{(x_{a_i} - x_{a_{i-1}})^2 + (y_{a_i} - y_{a_{i-1}})^2},\\
\text{subject to:} &F_i(x_{a_i}, y_{a_i}) = 0, \forall i \in \{0, 1, 2, \ldots, N-1\}.\\
&\{a_0, a_1, \ldots, a_{N-1}\} \in S_N
\end{split}
\end{equation}
\end{definition}

The solution is to obtain the corresponding location and order. The global optimum is the solution of Weak Form I.

Optimized variables are the location of points $(x_i,y_i)$ and the order of escape points $a_i$. There are $2N$ continuous variables and $N$ integer variables in above optimization.

\begin{definition}[Using subscript for Weak Form I]
As the order of escape points $a_i$ is a permutation of $\{0, 1, 2, \ldots, N-1\}$, we can use subscript to rewrite the constraints. 

$a_i$ is different from each other and is exactly one of $0, 1, ..., N-1$. It can be represented by an $N\times N$ matrix $s_{ij}$ with only one 1 in each row and column and all other elements are 0.

Eq (17) can be transformed into a binary optimization problem with subscripts:
\begin{equation}
\begin{split}
\text{minimize } &\sqrt{(x_{a_0} - 0)^2 + (y_{a_0} - 0)^2} + \sum_{i=1}^{N-1} \sqrt{(x_{a_i} - x_{a_{i-1}})^2 + (y_{a_i} - y_{a_{i-1}})^2},\\
\text{subject to:} &F_i(x_{a_i}, y_{a_i}) = 0,  \forall i \in \{0, 1, 2, \ldots, N-1\}, \\
&a_i = \sum_{j=0}^{N-1} j s_{ij}, \forall i \in \{0, 1, 2, \ldots, N-1\}, \\
&\sum_{j=0}^{N-1} s_{ij} = 1, \forall i \in \{0, 1, 2, \ldots, N-1\}, \\
&\sum_{i=0}^{N-1} s_{ij} = 1, \forall j \in \{0, 1, 2, \ldots, N-1\}, \\
&s_{ij} =
\begin{cases}
1, & \text{if } a_i = j, \, \forall i, j \in \{0, 1, 2, \ldots, N-1\}, \\
0, & \text{otherwise}.
\end{cases}
\end{split}
\end{equation}
\end{definition}

The solution is to obtain the value of binary variables for the corresponding order, and the value of continuous variables for location. Optimized variables are the location of points $(x_i, y_i)$ and the  binary variables $s_{ij}$. There are $2N$ continuous variables and $N^2$ binary variables in above optimization.

As $a_i$ is on the subscript in Eq (18), it is still difficult to solve directly. Solving an order of point is similar to the traveling salesman problem (TSP) and Hamiltonian path problem. Thus we can use the related formulation.

\begin{definition}[Using Miller–Tucker–Zemlin formulation for Weak Form I]
In the traditional TSP, the points are predefined and fixed. In this problem, however, the points are constrained to lie on the translated and rotated boundary curve. Leveraging the Hamiltonian path problem and the open-loop TSP, we can apply the Miller–Tucker–Zemlin formulation \cite{Miller1960}.
We define the decision variables:
\begin{equation*}
b_{i,j} = 
\begin{cases} 
1, & \text{if order from } i \text{ to } j, \\
0, & \text{otherwise}.
\end{cases}
\end{equation*}

and auxiliary variables $u_i$.
\end{definition}

Miller–Tucker–Zemlin formulation can be written as
\begin{itemize}
    \item Departure from each point once:
    \begin{equation*}
    \sum_{j=0}^{N-1} b_{i,j} = 1, \forall i \in \{0, 1, 2, \ldots, N-1\}.
    \end{equation*}
        \item Arrival to each point once:
    \begin{equation*}
    \sum_{i=0}^{N-1} b_{i,j} = 1, \forall j \in \{0, 1, 2, \ldots, N-1\}.
    \end{equation*}
        \item No self-loop:
    \begin{equation*}
    b_{i,i} = 0, \forall i \in \{0, 1, 2, \ldots, N-1\}.
    \end{equation*}
        \item Subtour Elimination Constraints:
    \begin{equation*}
    u_i - u_j + 1 \leq (N-2)(1 - b_{i,j}), \forall i, j \in \{1, 2, \ldots, N-1\}.
    \end{equation*}
        \item Auxiliary Variables Bounds:
    \begin{equation*}
    0 \leq u_i \leq N-1, \forall i \in \{0, 1, 2, \ldots, N-1\}.
    \end{equation*}
        \item Binary decision variables:
    \begin{equation*}
    b_{i,j} \in \{0, 1\}.
    \end{equation*}
\end{itemize}

The optimization formulation for Weak Form I can be written as:
\begin{equation}
\begin{split}
\text{minimize} &\sqrt{(x_0 - 0)^2 + (y_0 - 0)^2} + \sum_{i=0}^{N-1} \sum_{j=0}^{N-1} b_{i,j} \sqrt{(x_i - x_j)^2 + (y_i - y_j)^2},\\
\text{subject to:} &\sum_{j=0}^{N-1} b_{i,j} = 1, \forall i \in \{0, 1, 2, \ldots, N-1\}, \\
&\sum_{i=0}^{N-1} b_{i,j} = 1, \forall j \in \{0, 1, 2, \ldots, N-1\}, \\
&b_{i,i} = 0, \forall i \in \{0, 1, 2, \ldots, {N-1}\}, \\
&u_i - u_j + 1 \leq (N-2)(1 - b_{i,j}), \forall i, j \in \{1, 2, \ldots, N-1\}, \\
&F_i(x_i, y_i) = 0, \forall i \in \{0, 1, 2, \ldots, N-1\} \\
&0 \leq u_i \leq N-1, \\
&b_{i,j} \in \{0, 1\}
\end{split}
\end{equation}

The global optimum is the solution of Weak Form I. There are $3N$ continuous variables and $N^2$ binary variables in above optimization.

\begin{boldremark}
In the Miller–Tucker–Zemlin formulation above, the constraints represent passing through each point once, no self-loop, no subtour, and escape points on the boundary.
\end{boldremark}

\begin{corollary}
Write the square root as a constraint, Eq (19) can be reformulated as:
\begin{equation}
\begin{split}
\text{minimize} &c_{00} + \sum_{i=0}^{N-1} \sum_{j=0}^{N-1} b_{i,j} c_{i,j},\\
\text{subject to:} &c_{00}^2 = (x_0 - 0)^2 + (y_0 - 0)^2, 0 \leq c_{00}\\
&c_{i,j}^2 = (x_i - x_j)^2 + (y_i - y_j)^2, 0 \leq c_{i,j}, \forall i, j \in \{0, 1, 2, \ldots, {N-1}\}, \\
&\sum_{j=0}^{N-1} b_{i,j} = 1, \forall i \in \{0, 1, 2, \ldots, {N-1}\}, \\
&\sum_{i=0}^{N-1} b_{i,j} = 1, \forall j \in \{0, 1, 2, \ldots, {N-1}\}, \\
&b_{i,i} = 0, \forall i \in \{0, 1, 2, \ldots, {N-1}\}, \\
&u_i - u_j + 1 \leq (N-2)(1 - b_{i,j}), \forall i, j \in \{1, 2, \ldots, {N-1}\}, \\
&F_i(x_i, y_i) = 0, \forall i \in \{0, 1, 2, \ldots, N-1\}  \\
&0 \leq u_i \leq {N-1}, \\
&b_{i,j} \in \{0, 1\}.
\end{split}
\end{equation}

The above optimization can be solved by existing mixed integer programming (MIP) solvers.
\end{corollary}

\begin{boldremark}
We have listed the equations for solving Weak Form I. The difficulty is that TSP or MIP is an NP-hard problem \cite{Karp1972} \cite{Gavish 1978}. In particular it is necessary but extremely difficult to find the values of all binary variables used to determine the order of points.
\end{boldremark}

\subsection{Examples and Results}

Eqs (17-20) can be applied to solve Weak Form I. However, the computational cost is significantly high.

In this subsection, we introduce some assumptions for the order of points and provide some examples and results. Them are nontrivial and consistent with previous papers.

\subsubsection{Search for one line/half-plane with unit distance}

For one line/half-plane with unit distance in earliest study, as established by Isbell \cite{Isbell1957} the problem can be described as 
\\

``It is equivalent to require that the path $C$ meet every tangent to the unit circle about $p$ or that the convex hull of $C$ contain the circle". 
\\

And another paper \cite{Melzak2007} mentioned the problem as 
\\

``What is the shortest arc in the plane, starting from the origin and has a point in common with every straight line in the plane at a unit distance from origin". 
\\

It is a special case of Weak Form I. 

\begin{definition}[Rotated line with unit distance]
Recall Definition 3.11, rotated line with unit distance can be 
\begin{equation}
x\cos \frac{2\pi i}{N} + y\sin \frac{2\pi i}{N} - 1 = 0, \forall i \in \{0, 1, 2, \ldots, N-1\}.
\end{equation}
\end{definition}

The discrete formulation based on Eq (17) can be written as:
\begin{equation}
\begin{split}
\text{minimize } &\sqrt{(x_{a_0} - 0)^2 + (y_{a_0} - 0)^2} + \sum_{i=1}^{N-1} \sqrt{(x_{a_i} - x_{a_{i-1}})^2 + (y_{a_i} - y_{a_{i-1}})^2},\\
\text{subject to:} &x_{a_i} \cos \frac{2\pi i}{N} + y_{a_i} \sin \frac{2\pi i}{N} - 1 = 0, \forall i \in \{0, 1, 2, \ldots, N-1\}.\\
&\{a_0, a_1, \ldots, a_{N-1}\} \in S_N
\end{split}
\end{equation}

If assuming the order of points ranges from $0$ to $N-1$ as $\left \{ a_0,a_1,...,a_{N-1} \right \} =\left \{ 0,1,...,N-1 \right \}$.

Discrete formulation of the minimization problem can be written as:
\begin{equation}
\begin{split}
\text{minimize } &\sqrt{x_0^2 + y_0^2}+\sum_{i=1}^{N-1} \sqrt{(x_i - x_{i-1})^2 + (y_i - y_{i-1})^2} \\
\text{subject to:} &x_i \cos \frac{2\pi i}{N} + y_i \sin \frac{2\pi i}{N} - 1 = 0, \forall i \in \{0, 1, 2, \ldots, N-1\}.
\end{split}
\end{equation}

By solving above optimization, the results are shown in Figure 3:
\begin{figure}[H]
  \centering
  \includegraphics[width=5cm]{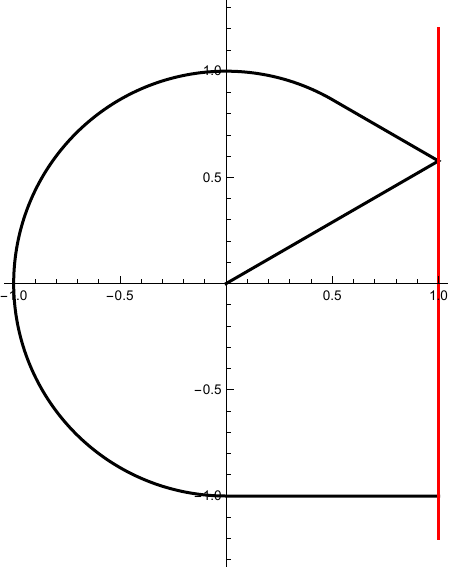}
  \caption{Results of search for one line/half plane with unit distance. (Black curve is escape path, and red curve is forest boundary)}
\end{figure}

\begin{boldremark}
The result consists of line segments and circular arcs. It was obtained and discussed in previous papers, but with geometric-based methods: Fig 1 in Ref \cite{Isbell1957},  Fig 2 in Ref \cite{Joris1980}, Fig 2 in page 151 in Ref \cite{Melzak2007}, Fig 6 in Ref \cite{Finch2004}, Fig 1(c) in Ref \cite{Zalgaller2005}. The exact solution is $\frac{7}{6} \pi +1+\sqrt{3}$.
\end{boldremark}

\begin{corollary}
With the assumption of order of points, when $N\to \infty $ the continuous formulation of the minimization problem can be written as:
\begin{equation}
\begin{split}
\text{minimize } &\sqrt{x(0)^2 + y(0)^2}+\int_0^{2\pi} \sqrt{x'(t)^2 + y'(t)^2} \, dt ,\\
\text{subject to:} & x(t) \cos t + y(t) \sin t - 1 = 0, \forall t \in [0, 2\pi].
\end{split}
\end{equation}

It is a constrained functional minimization problem, and it can be solved by calculus of variations. Mathematica code is in Appendix I.

Specifically, let functional
\begin{equation}
\mathcal{L}  = \sqrt{x'(t)^2 + y'(t)^2} + \lambda(t) \left[x(t)\cos t + y(t)\sin t - 1\right]
\end{equation}
where $\lambda(t)$ is the Lagrange multiplier.

Euler-Lagrange (EL) equations for $\mathcal{L}$ are:
\begin{equation}
\lambda(t)\cos t  + \frac{y'(t)\left[-y'(t)x''(t) + x'(t)y''(t)\right]}{\left[x'(t)^2 + y'(t)^2\right]^{3/2}} = 0
\end{equation}
\begin{equation}
\lambda(t)\sin t  + \frac{x'(t)\left[y'(t)x''(t) - x'(t)y''(t)\right]}{\left[x'(t)^2 + y'(t)^2\right]^{3/2}} = 0
\end{equation}

Substituting
\begin{equation}
x(t) = \frac{1 - y(t)\sin t}{\cos t}
\end{equation}

We obtain:
\begin{equation}
\left[\sin t - y(t)\right] \left\{\sec t \left[-3 + \cos 2t + 4y(t)\sin t\right] y'(t) + 4y'(t)^2 + 2\left[\sin t - y(t)\right] y''(t)\right\} = 0
\end{equation}

It is clear that the solution contains two parts:
\begin{equation}
\sin t - y(t) = 0
\end{equation}
or
\begin{equation}
\sec t \left[-3 + \cos 2t + 4y(t)\sin t\right] y'(t) + 4y'(t)^2 + 2\left[\sin t - y(t)\right] y''(t) = 0
\end{equation}
And the first part is the circle arc.
\end{corollary}

\begin{boldremark}
For other boundary functions not limited to this subsection, we can solve them in a similar way by modifying the constraints.
\end{boldremark}

\subsubsection{Search for a circle from exterior}
In this case, let outside distance to the circle is $r+s$ and the circle radius is $s$, and $r/s = 1$ the same setting in page 652 in \cite{Finch2004}. 

\begin{definition}[Rotated circle]
The rotated circle can be 
\begin{equation}
\left(x- \cos \frac{2\pi i}{N}\right)^2 + \left(y- \sin \frac{2\pi i}{N}\right)^2 - \left(\frac{1}{2}\right)^2 = 0\end{equation}
\end{definition}

Discrete formulation based on Eq (17) can be written as:
\begin{equation}
\begin{split}
\text{minimize } &\sqrt{(x_{a_0} - 0)^2 + (y_{a_0} - 0)^2} + \sum_{i=1}^{N-1} \sqrt{(x_{a_i} - x_{a_{i-1}})^2 + (y_{a_i} - y_{a_{i-1}})^2},\\
\text{subject to:} &\left(x_{a_i}- \cos \frac{2\pi i}{N}\right)^2 + \left(y_{a_i}- \sin \frac{2\pi i}{N}\right)^2 - \left(\frac{1}{2}\right)^2 = 0, \forall i \in \{0, 1, 2, \ldots, N-1\}.\\
&\{a_0, a_1, \ldots, a_{N-1}\} \in S_N
\end{split}
\end{equation}

If assuming the order of points ranges from $0$ to $N-1$ as $\left \{ a_0,a_1,...,a_{N-1} \right \} =\left \{ 0,1,...,N-1 \right \}$.

Discrete formulation of the minimization problem can be written as:
\begin{equation}
\begin{split}
\text{minimize } &\sqrt{x_0^2 + y_0^2}+\sum_{i=1}^{N-1} \sqrt{(x_i - x_{i-1})^2 + (y_i - y_{i-1})^2} ,\\
\text{subject to:} &\left( x_i - \cos \frac{2\pi i}{N} \right)^2 + \left( y_i - \sin \frac{2\pi i}{N} \right)^2 - \left( \frac{1}{2} \right)^2 = 0,  \forall i \in \{0, 1, 2, \ldots, N-1\}.
\end{split}
\end{equation}

By solving above optimization, the results are shown in Figure 4:
\begin{figure}[H]
  \centering
  \includegraphics[width=5cm]{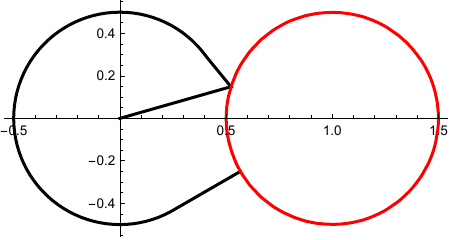}
  \caption{Results of search for a circle from exterior. (Black curve is escape path, and red curve is forest boundary)}
\end{figure}

\begin{boldremark}
The result consists of line segments and circular arcs. It was obtained and discussed in previous papers: Fig 1 in Ref \cite{Gluss1961}, Fig 7 in Ref \cite{Finch2004}, Fig 7 in Ref \cite{Baeza-Yates1988}, and Fig 3 in Ref \cite{Anderson2001}.
\end{boldremark}

With the assumption of order of points, when $N\to \infty $ continuous formulation of the minimization problem can be written as:
\begin{equation}
\begin{split}
\text{minimize } &\sqrt{x(0)^2 + y(0)^2}+\int_0^{2\pi} \sqrt{x'(t)^2 + y'(t)^2} \, dt ,\\
\text{subject to:} &\left[ x(t) - \cos t \right]^2 + \left[ y(t) - \sin t \right]^2 - \left( \frac{1}{2} \right)^2 = 0, \forall t \in [0, 2\pi].
\end{split}
\end{equation}

\subsubsection{Search for a circle from interior}
Let outside distance to circle be $r+s$ and circle radius be $s$, and $r/s = 0.2$ the same setting in page 652 in \cite{Finch2004}. 

Discrete formulation based on Eq (17) can be written as:
\begin{equation}
\begin{split}
\text{minimize } &\sqrt{(x_{a_0} - 0)^2 + (y_{a_0} - 0)^2} + \sum_{i=1}^{N-1} \sqrt{(x_{a_i} - x_{a_{i-1}})^2 + (y_{a_i} - y_{a_{i-1}})^2},\\
\text{subject to:} &\left(x_{a_i}- \cos \frac{2\pi i}{N}\right)^2 + \left(y_{a_i}- \sin \frac{2\pi i}{N}\right)^2 - 1.2^2 = 0, \forall i \in \{0, 1, 2, \ldots, N-1\}.\\
&\{a_0, a_1, \ldots, a_{N-1}\} \in S_N
\end{split}
\end{equation}

If assuming the order of points ranges from $0$ to $N-1$ as $\left \{ a_0,a_1,...,a_{N-1} \right \} =\left \{ 0,1,...,N-1 \right \}$.

Discrete formulation of the minimization problem can be written as:
\begin{equation}
\begin{split}
\text{minimize } &\sqrt{x_0^2 + y_0^2}+\sum_{i=1}^{N-1} \sqrt{(x_i - x_{i-1})^2 + (y_i - y_{i-1})^2} ,\\
\text{subject to:} &\left( x_i - \cos \frac{2\pi i}{N} \right)^2 + \left( y_i - \sin \frac{2\pi i}{N} \right)^2 - 1.2^2 = 0, \forall i \in \{0, 1, 2, \ldots, N-1\}.
\end{split}
\end{equation}

By solving above optimization, the results are shown in Figure 5:
\begin{figure}[H]
  \centering
  \includegraphics[width=5cm]{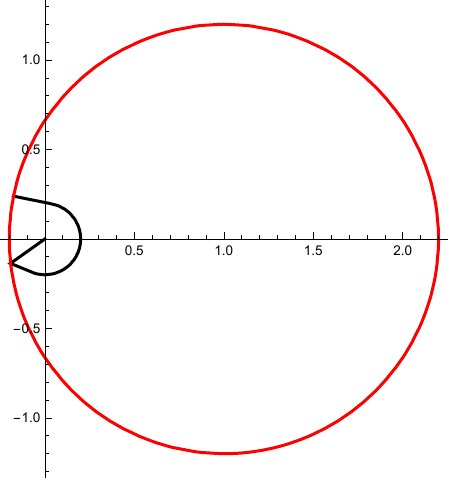}
  \caption{Results of search for a circle from interior. (Black curve is escape path, and red curve is forest boundary)}
\end{figure}

\begin{boldremark}
The result consists of line segments and circular arcs. It was obtained and discussed in previous papers: Fig 9 in Ref \cite{Finch2004}.
\end{boldremark}

With the assumption of order of points, when $N\to \infty $ continuous formulation of the minimization problem can be written as:
\begin{equation}
\begin{split}
\text{minimize } &\int_0^{2\pi} \sqrt{x'(t)^2 + y'(t)^2} \, dt + \sqrt{x(0)^2 + y(0)^2},\\
\text{subject to:} &\left[ x(t) - \cos t \right]^2 + \left[ y(t) - \sin t \right]^2 - 1.2^2 = 0, \forall t \in [0, 2\pi].
\end{split}
\end{equation}

\subsubsection{Search for a circle from interior when shortest path is not unique}

Let outside distance to circle be $r+s$ and circle radius be $s$, and $r/s = 0.333454$ the same setting in page 652 in \cite{Finch2004}.  

Discrete formulation based on Eq (17) can be written as:
\begin{equation}
\begin{split}
\text{minimize } &\sqrt{(x_{a_0} - 0)^2 + (y_{a_0} - 0)^2} + \sum_{i=1}^{N-1} \sqrt{(x_{a_i} - x_{a_{i-1}})^2 + (y_{a_i} - y_{a_{i-1}})^2},\\
\text{subject to:} &\left(x_{a_i}- \cos \frac{2\pi i}{N}\right)^2 + \left(y_{a_i}- \sin \frac{2\pi i}{N}\right)^2 - 1.500272^2 = 0, \forall i \in \{0, 1, 2, \ldots, N-1\}.\\
&\{a_0, a_1, \ldots, a_{N-1}\} \in S_N
\end{split}
\end{equation}

If assuming the order of points ranges from $0$ to $N-1$ as $\left \{ a_0,a_1,...,a_{N-1} \right \} =\left \{ 0,1,...,N-1 \right \}$.

Discrete formulation of the minimization problem can be written as:
\begin{equation}
\begin{split}
\text{minimize } &\sqrt{x_0^2 + y_0^2}+\sum_{i=1}^{N-1} \sqrt{(x_i - x_{i-1})^2 + (y_i - y_{i-1})^2} ,\\
\text{subject to:} &\left( x_i - \cos \frac{2\pi i}{N} \right)^2 + \left( y_i - \sin \frac{2\pi i}{N} \right)^2 - 1.500272^2 = 0, \forall i \in \{0, 1, 2, \ldots, N-1\}.
\end{split}
\end{equation}

By solving above optimization, the results are shown in Figure 6:
\begin{figure}[H]
  \centering
  \includegraphics[width=5cm]{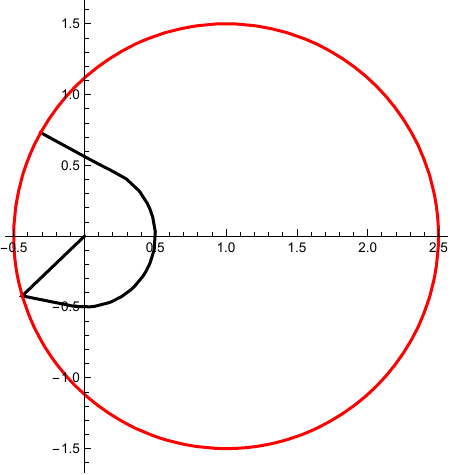}
  \caption{Results of search for a circle from interior when shortest path is not unique. (Black curve is escape path, and red curve is forest boundary)}
\end{figure}

\begin{boldremark}
In this setting, the solution is not unique. The length of optimal path in Figure 6 is 1. The diameter of the circle also represents a valid minimal solution \cite{Gross1955} \cite{Finch2004}. This is possible as there may be multiple global optima for Eqs (17-20).
\end{boldremark}

With the assumption of order of points, when $N\to \infty $ continuous formulation of the minimization problem can be written as:
\begin{equation}
\begin{split}
\text{minimize } &\sqrt{x(0)^2 + y(0)^2}+\int_0^{2\pi} \sqrt{x'(t)^2 + y'(t)^2} \, dt,\\
\text{subject to:} &\left[ x(t) - \cos t \right]^2 + \left[ y(t) - \sin t \right]^2 - 1.500272^2 = 0, \forall t \in [0, 2\pi].
\end{split}
\end{equation}

\subsubsection{Search for one point with given distance 1}
\begin{definition}[Rotated point with given distance 1]
The rotated point with given distance 1 can be 
\begin{equation}
\left(x- \cos \frac{2\pi i}{N}\right)^2 + \left(y- \sin \frac{2\pi i}{N}\right)^2 = 0 
\end{equation}
\end{definition}

Discrete formulation based on Eq (17) can be written as:
\begin{equation}
\begin{split}
\text{minimize } &\sqrt{(x_{a_0} - 0)^2 + (y_{a_0} - 0)^2} + \sum_{i=1}^{N-1} \sqrt{(x_{a_i} - x_{a_{i-1}})^2 + (y_{a_i} - y_{a_{i-1}})^2},\\
\text{subject to:} &\left(x_{a_i}- \cos \frac{2\pi i}{N}\right)^2 + \left(y_{a_i}- \sin \frac{2\pi i}{N}\right)^2 = 0, \forall i \in \{0, 1, 2, \ldots, N-1\}.\\
&\{a_0, a_1, \ldots, a_{N-1}\} \in S_N
\end{split}
\end{equation}

If assuming the order of points ranges from $0$ to $N-1$ as $\left \{ a_0,a_1,...,a_{N-1} \right \} =\left \{ 0,1,...,N-1 \right \}$.

Discrete formulation of the minimization problem can be written as:
\begin{equation}
\begin{split}
\text{minimize } &\sqrt{x_0^2 + y_0^2}+\sum_{i=1}^{N-1} \sqrt{(x_i - x_{i-1})^2 + (y_i - y_{i-1})^2} ,\\
\text{subject to:} &\left( x_i - \cos \frac{2\pi i}{N} \right)^2 + \left( y_i - \sin \frac{2\pi i}{N} \right)^2  = 0, \forall i \in \{0, 1, 2, \ldots, N-1\}.
\end{split}
\end{equation}

By solving above optimization, the results are shown in Figure 7:
\begin{figure}[H]
  \centering
  \includegraphics[width=5cm]{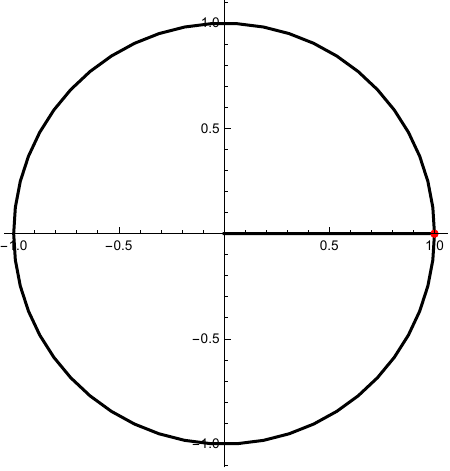}
  \caption{Results of search for one point with given distance 1. (Black curve is escape path, and red dot is forest boundary)}
\end{figure}

\begin{boldremark}
It was obtained and discussed in previous papers: Fig 8 in Ref \cite{Finch2004}, Fig 2(a) in Ref \cite{Melzak1973}. The exact solution is $1+2\pi$  as “move the known distance in any direction and then describe a circular path about the starting point” \cite{Melzak1973}. 
\end{boldremark}

With the assumption of order of points, when $N\to \infty $ continuous formulation of the minimization problem can be written as:
\begin{equation}
\begin{split}
\text{minimize } &\sqrt{x(0)^2 + y(0)^2}+\int_0^{2\pi} \sqrt{x'(t)^2 + y'(t)^2} \, dt ,\\
\text{subject to:} &\left[ x(t) - \cos t \right]^2 + \left[ y(t) - \sin t \right]^2 = 0, \forall t \in [0, 2\pi].
\end{split}
\end{equation}

\subsubsection{Search for a shape from interior with distance significantly short at certain angle, when shortest path is not unique}

David and Elmar \cite{Kübel2021} analyzed the dilemma of searching the area close-by versus going straight into one direction. ``While a breadth-first-search approach (BFS) will always result in an escape strategy when expanding the search-radii gradually, this is not true for a pure depth-first-search approach (DFS)". This is equivalent to considering the situations in 4.2.3 and 4.2.5 together. We construct a similar case to analyze. The forest boundary consists of a circle and a line segment. 

Assume that starting from the center of the circle. Length of line segment inside the circle is $1-\frac{1}{1+2\pi}$ , as Figure 8 shows. Then the path length of the straight radius is 1. According to 4.3.5, the path length of searching the line segment is also 1. 

\begin{definition}[Rotated circle and a line segment]
The rotated circle and a line segment can be 
\begin{equation}
\begin{split}
&\left[
\left( x_i \cos \frac{2\pi i}{N} - y_i \sin \frac{2\pi i}{N} \right)^2 
+ \left( x_i \sin \frac{2\pi i}{N} + y_i \cos \frac{2\pi i}{N} \right)^2 
- 1\right]\\
&\cdot 
\left[
\left( x_i \sin \frac{2\pi i}{N} + y_i \cos \frac{2\pi i}{N} \right)^2 
+ \left( \left( x_i \cos \frac{2\pi i}{N} - y_i \sin \frac{2\pi i}{N} - 1 \right) 
- \left( \frac{1}{1 + 2\pi} \right)^2 \right)^2 
\right] = 0 
\end{split}
\end{equation}
\end{definition}

We use product to unify circle and line constraints above.

The solution for this case with discrete formulation is
\begin{equation}
\begin{split}
\text{minimize}&\sqrt{x_0^2 + y_0^2} + \sum_{i=1}^{N-1} \sqrt{(x_i - x_{i-1})^2 + (y_i - y_{i-1})^2} \\
\text{subject to:}&\left[
\left( x_i \cos \frac{2\pi i}{N} - y_i \sin \frac{2\pi i}{N} \right)^2 
+ \left( x_i \sin \frac{2\pi i}{N} + y_i \cos \frac{2\pi i}{N} \right)^2 
- 1\right]\\
&\cdot 
\left[
\left( x_i \sin \frac{2\pi i}{N} + y_i \cos \frac{2\pi i}{N} \right)^2 
+ \left( \left( x_i \cos \frac{2\pi i}{N} - y_i \sin \frac{2\pi i}{N} - 1 \right) 
- \left( \frac{1}{1 + 2\pi} \right)^2 \right)^2 
\right] = 0,\\
& \forall i \in \{0, 1, 2, \ldots, N-1\}.
\end{split}
\end{equation}

By solving above optimization, the results are shown in Figure 8:
\begin{figure}[H]
  \centering
  \includegraphics[width=5cm]{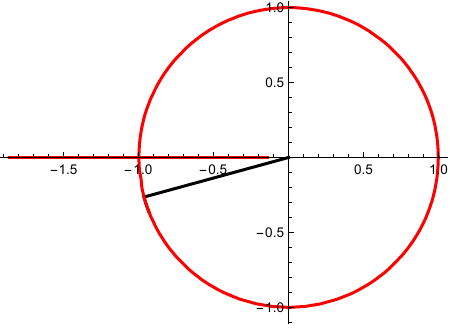} \\
  (a)\\
  \includegraphics[width=5cm]{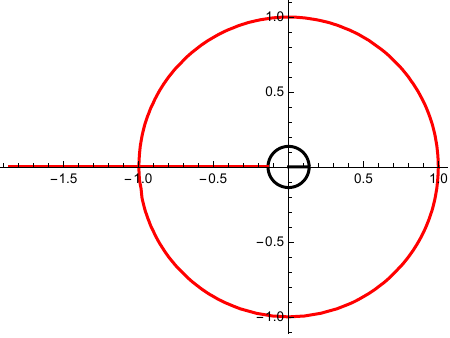}\\
  (b)
  \caption{Search for a shape from interior with distance significantly short at certain angle, when shortest path is not unique. (Black curve is escape path, and red curve is forest boundary consisting of a circle and a line segment)}
\end{figure}

\begin{boldremark}
In this setting, the solution is not unique. Both searching the circle and the line yield paths of the same length of 1, but with different shapes.
\end{boldremark}

\subsubsection{Search for two perpendicular lines with $\frac{1}{2}$ distance}
In this case, the discrete formulation based on Eq (17) can be written as:
\begin{equation}
\begin{split}
\text{minimize } &\sqrt{(x_{a_0} - 0)^2 + (y_{a_0} - 0)^2} + \sum_{i=1}^{N-1} \sqrt{(x_{a_i} - x_{a_{i-1}})^2 + (y_{a_i} - y_{a_{i-1}})^2},\\
\text{subject to:} &\left[-x_{a_i} \sin \frac{2\pi i}{N} + y_{a_i} \cos \frac{2\pi i}{N} + \frac{1}{2}\right]\cdot\left[-x_{a_i} \sin \left( \frac{2\pi i}{N} + \frac{\pi}{2} \right) + y_{a_i} \cos \left( \frac{2\pi i}{N} + \frac{\pi}{2} \right) + \frac{1}{2}\right] = 0, \\
&\forall i \in \{0, 1, 2, \ldots, N-1\}.\\
&\{a_0, a_1, \ldots, a_{N-1}\} \in S_N
\end{split}
\end{equation}
We use product to unify two line constraints.

If assuming order of points is $0, \frac{3N}{4}, 1, \frac{3N}{4} + 1, \ldots, \frac{N}{4}-1, N-1, \frac{N}{4}, \frac{N}{4}+1, \ldots, \frac{3N}{4}-1$, and $N$ can be divisible by 4.

Then note that the angle between the two lines is $\pi/2$, and we only need to use the rotation range of $3\pi/2$ for one straight line in the constraint. Because exceeding $3\pi/2$, the other line intersects the path.

Discrete formulation of the minimization problem can be simplified as:
\begin{equation}
\begin{split}
\text{minimize } &\sqrt{x_0^2 + y_0^2}+\sum_{i=1}^{N-1} \sqrt{(x_i - x_{i-1})^2 + (y_i - y_{i-1})^2} ,\\
\text{subject to:} &- x_i \sin \frac{3\pi i}{2N}  +  y_i \cos \frac{3\pi i}{2N}   +\frac{1}{2} = 0, \forall i \in \{0, 1, 2, \ldots, N-1\}.
\end{split}
\end{equation}

By solving above optimization, the results are shown in Figure 9:
\begin{figure}[H]
  \centering
  \includegraphics[width=5cm]{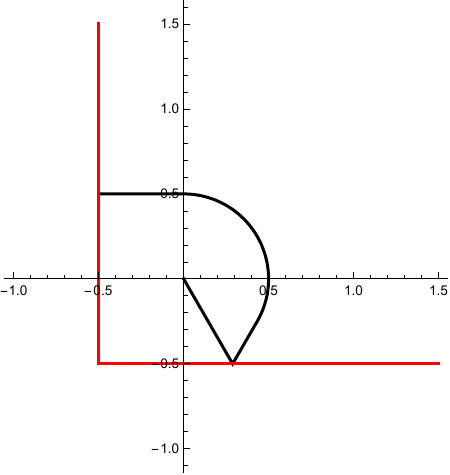}
  \caption{Results of search for two perpendicular lines with $\frac{1}{2}$ distance. (Black curve is escape path, and red curve is forest boundary)}
\end{figure}

\begin{boldremark}
The result consists of line segments and circular arcs.
\end{boldremark}

With the above assumptions, we only need to consider $ \left[ 0, \frac{3\pi}{2} \right]$. The continuous formulation of the minimization problem can be written as:
\begin{equation}
\begin{split}
\text{minimize } &\sqrt{x(0)^2 + y(0)^2}+\int_0^{3\pi/2} \sqrt{x'(t)^2 + y'(t)^2} \, dt ,\\
\text{subject to:} & -x(t) \sin t  +  y(t) \cos t  +\frac{1}{2} = 0, \forall t \in [0, 3\pi/2].
\end{split}
\end{equation}

\subsubsection{Search for two parallel lines with unit distance (unit strip) from the middle}

Discrete formulation based on Eq (17) can be written as:
\begin{equation}
\begin{split}
\text{minimize } &\sqrt{(x_{a_0} - 0)^2 + (y_{a_0} - 0)^2} + \sum_{i=1}^{N-1} \sqrt{(x_{a_i} - x_{a_{i-1}})^2 + (y_{a_i} - y_{a_{i-1}})^2},\\
\text{subject to:} &\left[
-x_{a_i} \sin \frac{2\pi i}{N} + y_{a_i} \cos \frac{2\pi i}{N} + \frac{1}{2}
\right]\cdot
\left[
-x_{a_i} \sin \left( \frac{2\pi i}{N} + \pi \right) + y_{a_i} \cos \left( \frac{2\pi i}{N} + \pi \right)+\frac{1}{2}
\right] = 0, \\
&\forall i \in \{0, 1, 2, \ldots, N-1\}.
\end{split}
\end{equation}

If assuming order of points is $0, N/2, 1, N/2+1, \ldots, N/2-1, N-1$, and $N$ can be divisible by 2.

Then note that the angle between the two parallel lines can be regarded as $\pi$, and we only need to use the rotation range of $\pi$ for one straight line in the constraint. Because exceeding $\pi$, the other line intersects the path.

Discrete formulation of the minimization problem can be simplified as:
\begin{equation}
\begin{split}
\text{minimize } &\sqrt{x_0^2 + y_0^2}+\sum_{i=1}^{N-1} \sqrt{(x_i - x_{i-1})^2 + (y_i - y_{i-1})^2} ,\\
\text{subject to:} & x_i \cos \frac{\pi i}{N}  +  y_i \sin \frac{\pi i}{N}  +\frac{1}{2} = 0, \forall i \in \{0, 1, 2, \ldots, N-1\}.
\end{split}
\end{equation}

By solving above optimization, the results are shown in Figure 10:
\begin{figure}[H]
  \centering
  \includegraphics[width=5cm]{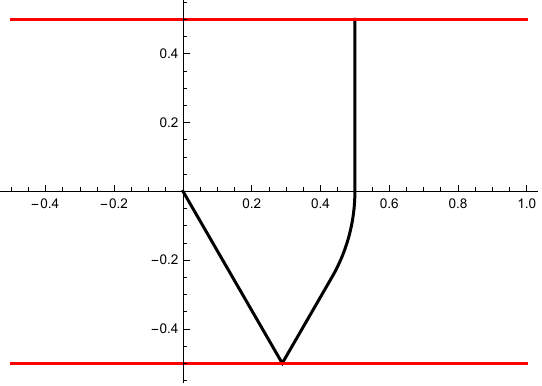}
  \caption{Results of search for two parallel lines with unit distance (unit strip) from the middle. (Black curve is escape path, and red curve is forest boundary)}
\end{figure}

\begin{boldremark}
The result consists of line segments and circular arcs. It was mentioned in Croft’s book page 41 \cite{Croft2012}.
\end{boldremark}

With the above assumptions, continuous formulation of the minimization problem can be written as:
\begin{equation}
\begin{split}
\text{minimize } &\sqrt{x(0)^2 + y(0)^2}+\int_0^{\pi} \sqrt{x'(t)^2 + y'(t)^2} \, dt ,\\
\text{subject to:} & x(t) \cos t +  y(t) \sin t  +\frac{1}{2} = 0, \forall t \in [0, \pi].
\end{split}
\end{equation}

\subsubsection{Search for two lines from the angle bisector in the middle}
In this subsection, we consider the case of searching for two lines from the angle bisector in the middle. In 4.2.7 we considered two perpendicular lines with an angle of $\pi/2$, and in 4.2.8 we considered two parallel lines with an angle of 0. This section is the general case.

Assume that two lines are at certain angle $\theta \in \left [ 0,2\pi  \right ] $. Assume the distance from the starting point to the two straight lines is $1/2$. Then we only need to use the rotation range of $2\pi-\theta$ for one straight line in the constraint. Because exceeding $2\pi-\theta$, the other line intersects the path.

Discrete formulation of the minimization problem can be written as
\begin{equation}
\begin{split}
\text{minimize } &\sqrt{x_0^2 + y_0^2}+\sum_{i=1}^{N-1} \sqrt{(x_i - x_{i-1})^2 + (y_i - y_{i-1})^2} ,\\
\text{subject to:} & x_i \cos \frac{(2\pi-\theta) i}{N}  +  y_i \sin \frac{(2\pi-\theta) i}{N}  +\frac{1}{2} = 0, \forall i \in \{0, 1, 2, \ldots, N-1\}.
\end{split}
\end{equation}

By solving above optimization, the results are shown in Figure 11:
\begin{figure}[H]
  \centering
  (a)\includegraphics[width=5cm]{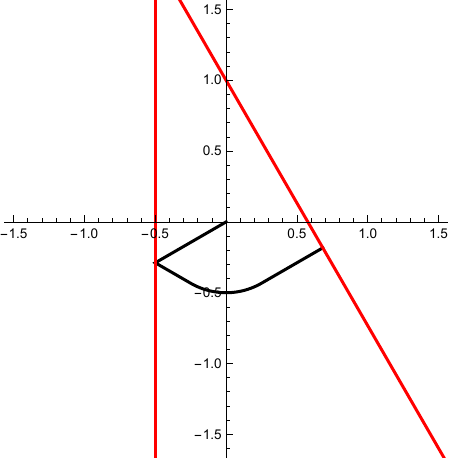}
  (b)\includegraphics[width=5cm]{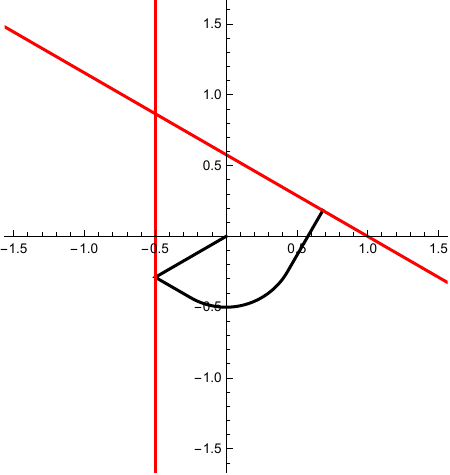}\\
  (c)\includegraphics[width=5cm]{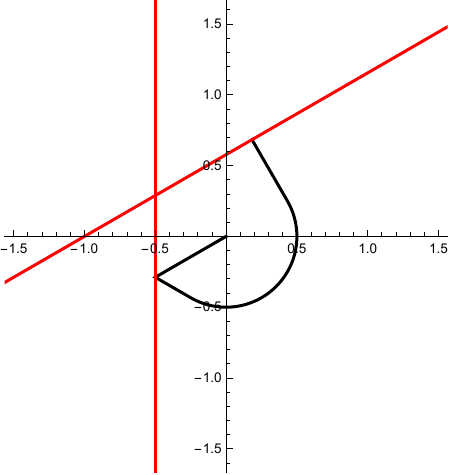}
  (d)\includegraphics[width=5cm]{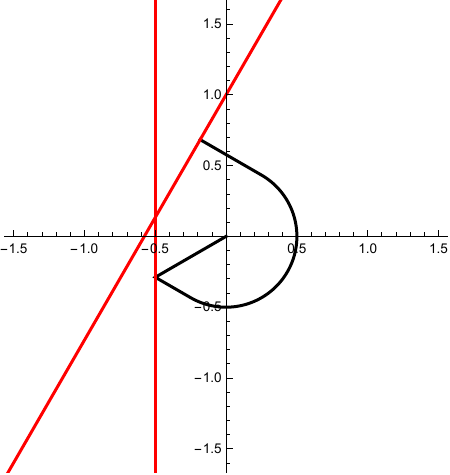}
  \caption{Search for two lines from the angle bisector in the middle at certain angle (a)$\pi/6, (b)\pi/3, (c)2\pi/3, (d)5\pi/6$. (Black curve is escape path, and red curve is forest boundary)}
\end{figure}

\subsubsection{Search for two parallel lines with unit distance (unit strip) from edge (Zalgaller Class 2)}

If assuming order of points ranges from $0$ to $N-1$.

As Figure 12 shows, assume that the path consists of two parts. The path first reaches the opposite straight line $line2$ at $(x_0,y_0)$, and rotates along the boundary straight line $line2$ until the other rotated boundary straight line $line1$ reaches the turning point $(x_0,y_0)$. The end of path is $(x_{N-1},y_{N-1})$. Figure 12 also shows the range of forest boundary rotation.

Discrete formulation of the minimization problem can be written as:
\begin{equation}
\begin{split}
\text{minimize } &\sqrt{x_0^2 + y_0^2}+\sum_{i=1}^{N-1} \sqrt{(x_i - x_{i-1})^2 + (y_i - y_{i-1})^2} ,\\
\text{subject to:} &x_i \cos \left( -\frac{\pi - \arctan \frac{y_0}{x_0}}{N} i -\frac{\pi}{2}\right) + y_i \sin \left(- \frac{\pi - \arctan \frac{y_0}{x_0}}{N} i-\frac{\pi}{2} \right) - 1 = 0, \\
&\forall i \in \{0, 1, 2, \ldots, N-1\}.
\end{split}
\end{equation}

By solving above optimization, the results are shown in Figure 12:
\begin{figure}[H]
  \centering
  \includegraphics[width=14cm]{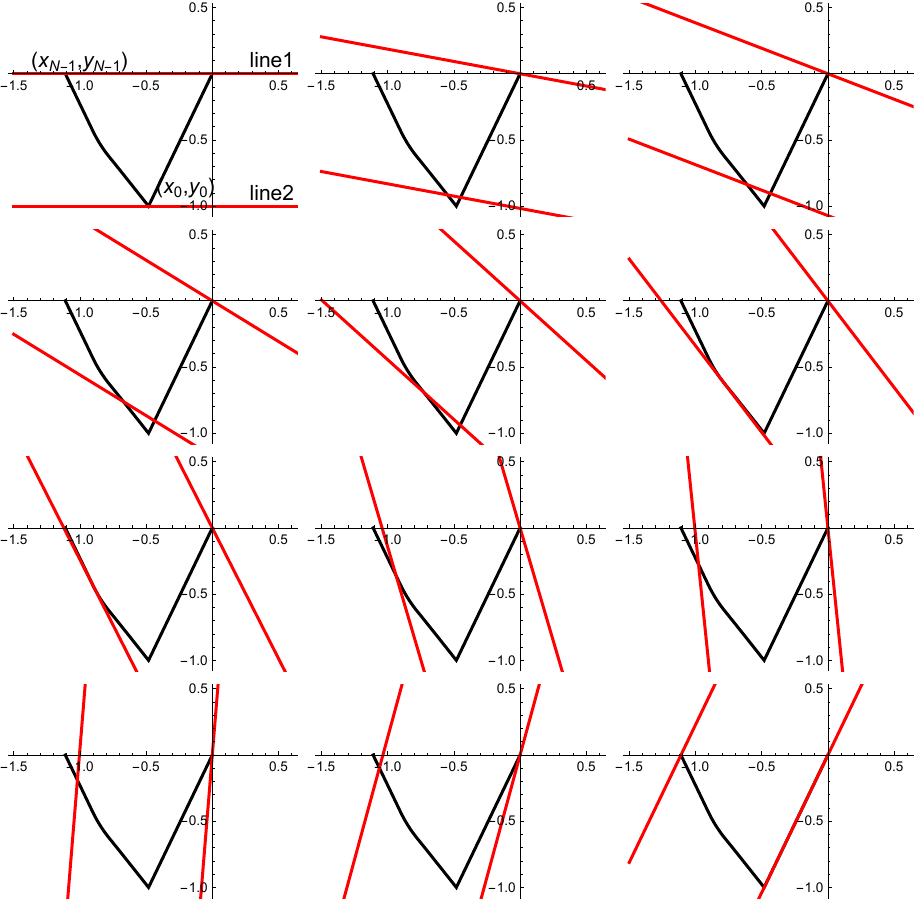}
  \caption{Results of search for two parallel lines with unit distance (unit strip) from edge (Zalgaller Class 2 \cite{Zalgaller2005}). (Black curve is escape path, and red curve is forest boundary)}
\end{figure}

\begin{boldremark}
The result in Figure 12 consists of three line segments and one circular arc. The length of obtained path is 2.297, same as Zalgaller’s paper section 5.5 \cite{Zalgaller2005}. It was obtained and discussed in Zalgaller's paper\cite{Zalgaller2005} for unit strip as Class 2 Fig 6(b). This is suboptimal. The optimal result of Zalgaller for unit strip is 2.278292 and it should be symmetrical with four line segments and two circular arcs. We will present this result in section 6.2.1. 
\end{boldremark}

Discrete formulation of 4.2.1-4.2.10 can be solved by global optimization solver. Mathematica code is in Appendix II.

\section{Weak Form II - Starting from several possible known points with unknown orientation}

Recall Definition 3.22, Weak Form II is starting from several known points with unknown orientation.

As in Section 4, the solution to Weak Form II requires finding the shortest path through all the escape points, as illustrated in Figure 2. To formulate the problem as a constrained optimization, we consider $N$ orientations and $M$ starting points. We then seek a minimal-length path that passes through escape point corresponding to each possible orientation and starting point.

\begin{definition}[One-to-one mapping]
We define a one-to-one mapping 
\begin{equation}
\left \{ i,k \right \} \to h,  \forall i \in \{0, 1, 2, \ldots, N-1\}, \, \forall k \in \{1, 2, \ldots, M\}, \, \exists h \in \{1, 2, \ldots, MN\}
\end{equation}
\end{definition}

\begin{definition}[With M starting points using all permutation to define]
With one-to-one mapping $\left \{ i,k \right \} \to h$, the Weak Form II can be written as follows:
\begin{equation}
\begin{split}
\text{minimize } &\sqrt{(x_{a_1} - 0)^2 + (y_{a_1} - 0)^2} + \sum_{i=2}^{MN} \sqrt{(x_{a_i} - x_{a_{i-1}})^2 + (y_{a_i} - y_{a_{i-1}})^2},\\
\text{subject to:} &F_{ki}(x_{a_h}, y_{a_h}) = 0, \forall i \in \{0, 1, 2, \ldots, N-1\}, \, \forall k \in \{1, 2, \ldots, M\}, \, \exists h \in \{1, 2, \ldots, MN\}.
\end{split}
\end{equation}
The optimized variable in the above equation is a permutation from 1 to $MN$, and $2MN$ continuous variables.
\end{definition}

\begin{definition}[Using subscript for Weak Form II]
Similar to Definition 4.2, Eq (57) can be reformulated as a binary optimization problem with subscripts:
\begin{equation}
\begin{split}
\text{minimize } &\sqrt{(x_{a_1} - 0)^2 + (y_{a_1} - 0)^2} + \sum_{i=2}^{MN} \sqrt{(x_{a_i} - x_{a_{i-1}})^2 + (y_{a_i} - y_{a_{i-1}})^2},\\
\text{subject to:} &F_{ki}(x_{a_h}, y_{a_h}) = 0, \forall i \in \{0, 1, 2, \ldots, N-1\}, \, \forall k \in \{1, 2, \ldots, M\}, \, \exists h \in \{1, 2, \ldots, MN\},\\
&a_h = \sum_{j=1}^{MN} j s_{hj}, \forall h \in \{1, 2, \ldots, MN\},\\
&\sum_{j=1}^{MN} s_{hj} = 1, \forall h \in \{1, 2, \ldots, MN\},\\
&\sum_{h=1}^{MN} s_{hj} = 1, \forall j \in \{1, 2, \ldots, MN\},\\
&s_{hj} = 
\begin{cases} 
&1, a_h = j, \forall j \in \{1, 2, \ldots, MN\}, \\
&0, \text{otherwise}.
\end{cases}
\end{split}
\end{equation}
The optimized variable in the above equation is $MN$ binary variables, and $2MN$ continuous variables.
\end{definition}

\begin{definition}[Using Miller–Tucker–Zemlin formulation for Weak Form II]
Using TSP and the Miller–Tucker–Zemlin formulation \cite{Miller1960}, Weak Form II can be expressed as:
\begin{equation}
\begin{split}
\text{minimize } &\sqrt{(x_1 - 0)^2 + (y_1 - 0)^2} + \sum_{i=1}^{MN} \sum_{j=1}^{MN} s_{i,j} \sqrt{(x_i - x_j)^2 + (y_i - y_j)^2},\\
\text{subject to:} &\sum_{j=1}^{MN} s_{i,j} = 1, \forall i \in \{1, 2, \ldots, MN\},\\
&\sum_{i=1}^{MN} s_{i,j} = 1, \forall j \in \{1, 2, \ldots, MN\},\\
&s_{i,i} = 0, \quad \forall i \in \{1, 2, \ldots, MN\},\\
&u_i - u_j + 1 \leq MN (1 - s_{i,j}), \forall i, j \in \{1, 2, \ldots, MN\},\\
&F_{ki}(x_h, y_h) = 0, \forall h \in \{1, 2, \ldots, MN\},\\
&0 \leq u_i \leq MN - 1,\\
&s_{i,j} \in \{0, 1\}.
\end{split}
\end{equation}
The optimized variable in the above equation is $MN$ binary variables, and $3MN$ continuous variables. It is a bigger TSP variation.

\begin{boldremark}
When $M=1$, Weak Form II degenerates to Weak Form I.
\end{boldremark}
\end{definition}

\section{Original Bellman’s lost-in-a-forest problem}
\subsection{General solution to Bellman’s lost-in-a-forest problem}

Recall Definition 3.22 and last Section 5, we provide the solution to Weak Form II. 

\begin{theorem}
When $M$ points are evenly distributed in the region like grid points with $x$ and $y$ direction increment $\Delta \to 0$, and when $N$ orientations are evenly distributed in $[0, 2\pi]$ with increment $2\pi/N \to 0$, and $M$ and $N$ become infinitely large, the solution to Weak Form II yields the solution to the original Bellman’s lost-in-a-forest problem.
\end{theorem}

\begin{proof}
Since the starting point region $R \in \mathbb{R}^2$ and interval $[0, 2\pi]$ are compact and continuous, the Weierstrass Extreme Value Theorem can guarantee the existence of a global minimum. Uniform continuity ensures that as the discretization gets finer, the discretized minimum values approach the global minimum value.

When $M$ and $N$ are infinitely large, the Riemann sum of curve length converges to the integral curve length. Then we can use $\Gamma$-convergence \cite{Braides2002} to prove and establish the convergence of the minimization.

Liminf Inequality: For any sequence $(x_{MN}, y_{MN}) \to (x, y)$:
\begin{equation}
\lim_{M, N \to \infty} \inf L_{MN}(x_{MN}, y_{MN}) \geq L(x, y).
\end{equation}
This follows from the lower semicontinuity of the length functional and the convergence of $L_{MN}$ to $L$.

Recovery Sequence: For any region $R \in \mathbb{R}^2$ and interval $[0, 2\pi]$, define $\tilde{x}_{MN}, \tilde{y}_{MN}$. Then:
\begin{equation}
\lim_{M, N \to \infty} \sup L_{MN}(\tilde{x}_{MN}, \tilde{y}_{MN}) \leq L(x, y).
\end{equation}
As $x$ and $y$ are continuous, this sequence approximates $(x, y)$ well.

Since both conditions are satisfied, $L_{MN} \, \Gamma$-converges to $L$. The polyline converges to the shortest curve. The proof is standard in shape optimization via piecewise approximation.
\end{proof}

\subsection{Examples of various forest boundary shapes}
\subsubsection{Unit strip}
Previously, Zalgaller established the solution for the unit strip using geometric methods \cite{Zalgaller2005}. Other authors also discussed this case \cite{Ward2008} \cite{Adhikari1989} \cite{Finch2004-1}. It was also called the broadworm. 

\begin{definition}[Rotated and translated unit strip]
Recall Definition 3.11, the rotated and translated unit strip boundary can be
\begin{equation}
\begin{aligned}
\left[
x_{a_h} \cos \frac{2\pi i}{N} + y_{a_h} \sin \frac{2\pi i}{N} + \frac{k-1}{2M}
\right]
\cdot
\left[
x_{a_h} \cos \frac{2\pi i}{N} + y_{a_h} \sin \frac{2\pi i}{N} + \frac{k-1}{2M} - 1
\right] = 0
\end{aligned}
\end{equation}
\end{definition}

The solution for escaping unit strip with discrete formulation can be expressed as: 
\begin{equation}
\begin{aligned}
\text{minimize} & \sqrt{(x_{a_1} - 0)^2 + (y_{a_1} - 0)^2} + \sum_{i=2}^{MN} \sqrt{(x_{a_i} - x_{a_{i-1}})^2 + (y_{a_i} - y_{a_{i-1}})^2}, \\
\text{subject to:} & 
\left[
x_{a_h} \cos \frac{2\pi i}{N} + y_{a_h} \sin \frac{2\pi i}{N} + \frac{k-1}{2M}
\right]
\cdot
\left[
x_{a_h} \cos \frac{2\pi i}{N} + y_{a_h} \sin \frac{2\pi i}{N} + \frac{k-1}{2M} - 1
\right] = 0, \\
& \forall i \in \{0, 1, 2, \ldots, N-1\}, \, \forall k \in \{1, 2, \ldots, M\}, \, \exists h \in \{1, 2, \ldots, MN\}.
\end{aligned}
\end{equation}

\begin{boldremark}
Considering symmetry, we can only consider starting point $s_y\in\left [ 0,0.5 \right ] $ as shown in Figure 14 (a). 
\end{boldremark}

For large $M$ and $N$ in above optimization, it is very difficult to solve for the NP-hard problem. We assume the order of the escape points and can calculate the solution of $x_{a_i}$ and $y_{a_i}$ for $N=12, M=26$ as shown in Figure 13, which is very close to Zalgaller's \cite{Zalgaller2005}.
\begin{figure}[H]
  \centering
  \includegraphics[width=10cm]{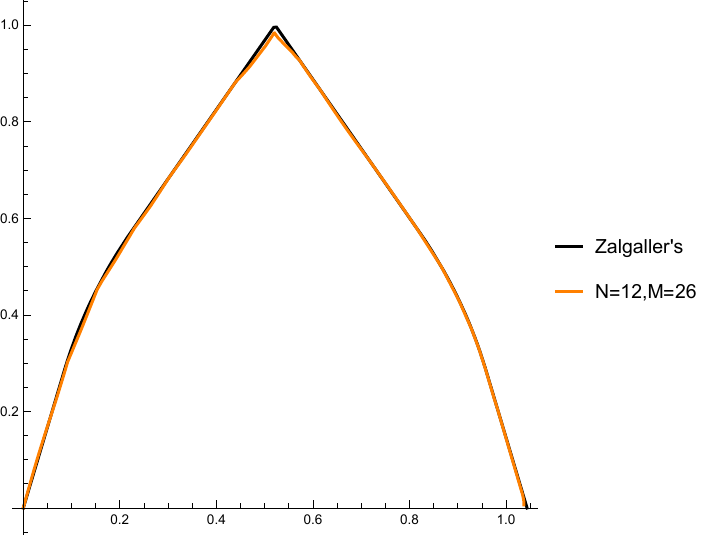}
  \caption{Results of search for two parallel lines with unit distance (unit strip) when $N=12, M=26$ comparing with Zalgaller's \cite{Zalgaller2005}). }
\end{figure}

Mathematica code for obtaining results in Figure 13 is in Appendix III.

\subsubsection{Unit circle}
The diameter has been proven to be the optimal path to escape from a circle (Theorem 4 in \cite{Finch2004}, Proposition 6 in \cite{Ward2008}). 

\begin{definition}[Rotated and translated circle]
Recall Definition 3.11, the rotated and translated circle can be
\begin{equation}
\left(  x_{a_h}-\frac{k}{M}\cos\frac{2\pi i}{N}  \right)^2
+ \left( y_{a_h}- \frac{k}{M}\sin\frac{2\pi i}{N} \right)^2  - 1 = 0
\end{equation}
\end{definition}

The solution for escaping a circle with discrete formulation can be
\begin{equation}
\begin{aligned}
\text{minimize} & \sqrt{(x_{a_1} - 0)^2 + (y_{a_1} - 0)^2} + \sum_{i=2}^{MN} \sqrt{(x_{a_i} - x_{a_{i-1}})^2 + (y_{a_i} - y_{a_{i-1}})^2}, \\
\text{subject to:}  & 
\left(  x_{a_h}-\frac{k}{M}\cos\frac{2\pi i}{N}  \right)^2
+ \left( y_{a_h}- \frac{k}{M}\sin\frac{2\pi i}{N} \right)^2  - 1 = 0, \\
& \forall i \in \{0, 1, 2, \ldots, N-1\}, \, \forall k \in \{1, 2, \ldots, M\}, \, \exists h \in \{1, 2, \ldots, MN\}.
\end{aligned}
\end{equation}

\begin{boldremark}
Considering symmetry, we can only consider starting points on a radius as shown in Figure 14 (b).
\end{boldremark}

\subsubsection{Equilateral triangle}
Besicovitch’s zigzag path is the optimal for equilateral triangle, and shorter than triangle side length \cite{Besicovitch1965}. This path has been discussed and extended in several papers \cite{Ward2008} \cite{Movshovich2012} \cite{Coulton2006} \cite{Movshovich2011}.

In this paper, the solution to escaping equilateral triangle with discrete formulation is 
\begin{equation}
\begin{aligned}
\text{minimize} & \sqrt{(x_{a_1} - 0)^2 + (y_{a_1} - 0)^2} + \sum_{i=2}^{MN} \sqrt{(x_{a_i} - x_{a_{i-1}})^2 + (y_{a_i} - y_{a_{i-1}})^2}, \\
\text{subject to:}  & 
\Bigg[
\left( x_{a_h} - s_{kx} \right) \cos\left( -\frac{\pi}{2} + \frac{2\pi i}{N} \right) 
+ \left( y_{a_h} - s_{ky} \right) \sin\left( -\frac{\pi}{2} + \frac{2\pi i}{N} \right) - \frac{\sqrt{3}}{6}
\Bigg] \\
& \cdot
\Bigg[
\left( x_{a_h} - s_{kx} \right) \cos\left( -\frac{\pi}{6} + \frac{2\pi i}{N} \right) 
+ \left( y_{a_h} - s_{ky} \right) \sin\left( -\frac{\pi}{6} + \frac{2\pi i}{N} \right) + \frac{\sqrt{3}}{6}
\Bigg] \\
& \cdot
\Bigg[
\left( x_{a_h} - s_{kx} \right) \cos\left( \frac{\pi}{6} + \frac{2\pi i}{N} \right) 
+ \left( y_{a_h} - s_{ky} \right) \sin\left( \frac{\pi}{6} + \frac{2\pi i}{N} \right) - \frac{\sqrt{3}}{6}
\Bigg] = 0, \\
& \forall i \in \{0, 1, 2, \ldots, N-1\}, \, \forall k \in \{1, 2, \ldots, M\}, \, \exists h \in \{1, 2, \ldots, MN\}.
\end{aligned}
\end{equation}

where $(s_{kx},s_{ky})$  are starting points evenly distributed as uniform grid within one-third of the equilateral triangle region, as shown in Figure 14 (c).

\begin{boldremark}
Considering symmetry, we can only consider starting from one-third of the equilateral triangle region.
\end{boldremark}

\subsubsection{Arbitrary triangle}
Although Gibbs applied numerical methods to determine escape paths for isosceles triangles \cite{Gibbs2016}, no general solution exists for arbitrary triangles. 

In this paper, the solution for escaping any triangle with discrete formulation is 
\begin{equation}
\begin{aligned}
\text{minimize} & \sqrt{(x_{a_1} - 0)^2 + (y_{a_1} - 0)^2} + \sum_{i=2}^{MN} \sqrt{(x_{a_i} - x_{a_{i-1}})^2 + (y_{a_i} - y_{a_{i-1}})^2}, \\
\text{subject to:}  & 
\Bigg[
\left( x_{a_h} - s_{kx} \right) \cos\left( \phi_1 + \frac{2\pi i}{N} \right) 
+ \left( y_{a_h} - s_{ky} \right) \sin\left( \phi_1 + \frac{2\pi i}{N} \right) - \delta_1
\Bigg] \\
& \cdot
\Bigg[
\left( x_{a_h} - s_{kx} \right) \cos\left( \phi_2 + \frac{2\pi i}{N} \right) 
+ \left( y_{a_h} - s_{ky} \right) \sin\left( \phi_2 + \frac{2\pi i}{N} \right) - \delta_2
\Bigg] \\
& \cdot
\Bigg[
\left( x_{a_h} - s_{kx} \right) \cos\left( \phi_3 + \frac{2\pi i}{N} \right) 
+ \left( y_{a_h} - s_{ky} \right) \sin\left( \phi_3 + \frac{2\pi i}{N} \right) - \delta_3
\Bigg] = 0, \\
& \forall i \in \{0, 1, 2, \ldots, N-1\}, \, \forall k \in \{1, 2, \ldots, M\}, \, \exists h \in \{1, 2, \ldots, MN\}.
\end{aligned}
\end{equation}

where $\phi_1, \phi_2, \phi_3, \delta_1, \delta_2, \delta_3$ are given parameters for the three lines of the triangle, $\left(s_{kx}, s_{ky}\right)$ are starting points evenly distributed as uniform grid within the triangle region, as shown in Figure 14 (d).

\begin{boldremark}
This approach can also be extended to polygons.
\end{boldremark}

\subsubsection{General circular sector}
Circular sectors have been extensively studied, with some conjectures relating them to the unit-length curve coverage problem \cite{Movshovich2017} \cite{Wetzel2019}. 

In this paper, the solution for escaping any circle sector with discrete formulation is 
\begin{equation}
\begin{aligned}
    \text{minimize} &\sqrt{(x_{a_1} - 0)^2 + (y_{a_1} - 0)^2} + \sum_{i=2}^{MN} \sqrt{(x_{a_i} - x_{a_{i-1}})^2 + (y_{a_i} - y_{a_{i-1}})^2}, \\
    \text{subject to:} &\left[ (x_{a_h} - s_{kx}) \cos\left( \phi_1 + \frac{2\pi i}{N} \right) + (y_{a_h} - s_{ky}) \sin\left( \phi_1 + \frac{2\pi i}{N} \right) \right] \\
        &\cdot\left[ (x_{a_h} - s_{kx}) \cos\left( \phi_2 + \frac{2\pi i}{N} \right) + (y_{a_h} - s_{ky}) \sin\left( \phi_2 + \frac{2\pi i}{N} \right) \right] \\
        &\cdot\left[ \left ( x-\sqrt{s_{kx}^2+s_{ky}^2}\cos \frac{2\pi i}{N}   \right )^2+ \left ( y-\sqrt{s_{kx}^2+s_{ky}^2}\sin \frac{2\pi i}{N}   \right )^2-r^2 \right] = 0,\\
        &\forall i \in \{0, 1, 2, \ldots, N-1\}, \forall k \in \{1, 2, \ldots, M\}, \exists h \in \{1, 2, \ldots, MN\}.
\end{aligned}
\end{equation}

where $\phi_1, \phi_2$ are given parameters for the two lines of circle sector, $r$ is radius, $(s_{kx},s_{ky})$ are starting points evenly distributed as uniform grid within circular sector, as shown in Figure 14 (e).

\begin{figure}[H]
  \centering
  (a)\includegraphics[width=5cm]{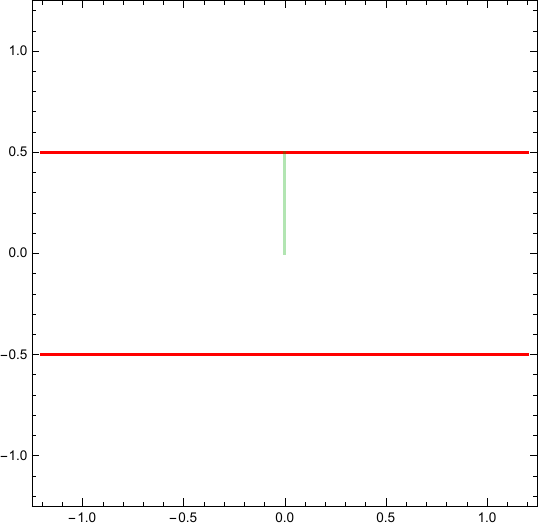} 
  (b)\includegraphics[width=5cm]{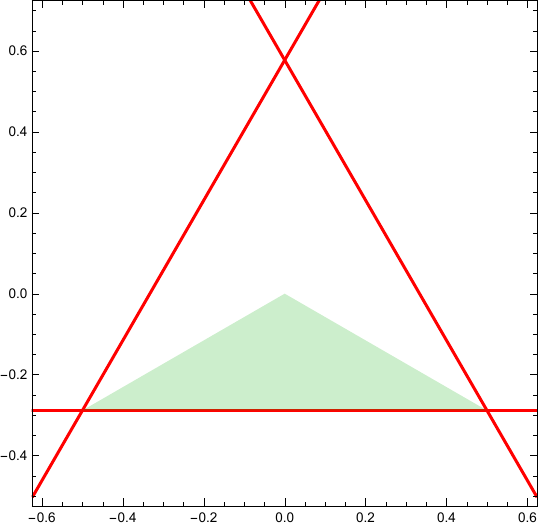} 
  (c)\includegraphics[width=5cm]{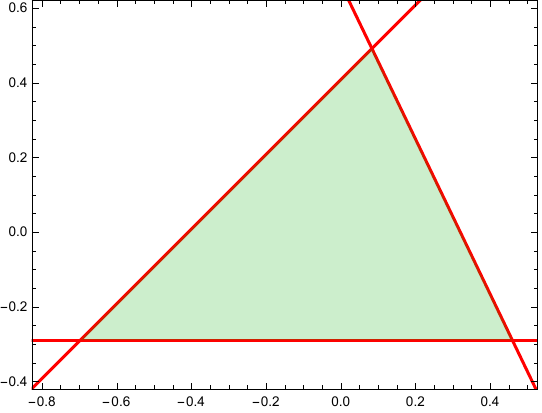}
  (d)\includegraphics[width=5cm]{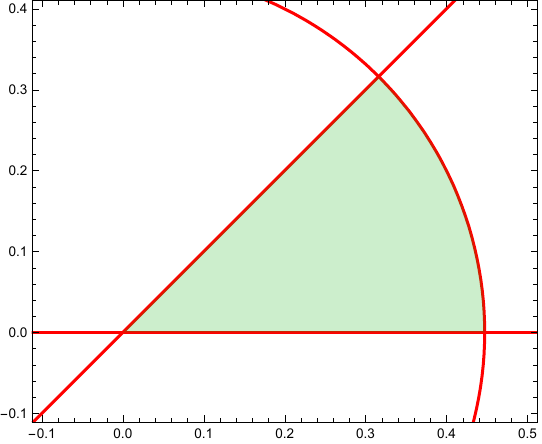} 
  (e)\includegraphics[width=5cm]{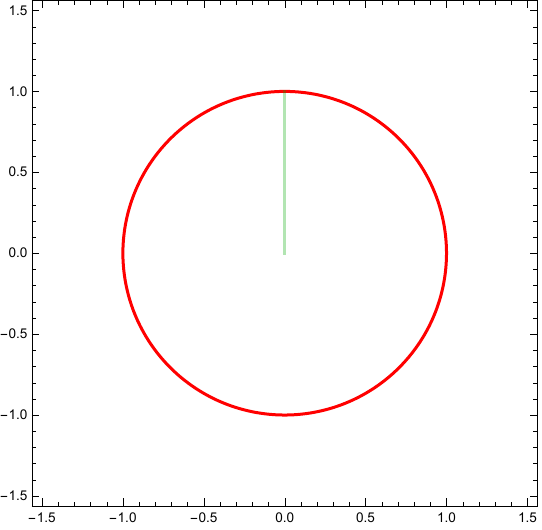} 
  \caption{Starting points evenly distributed with the region to escape: (a)unit stripe, (b)unit circle, (c) one third of equilateral triangle, (d) general triangle, (e) circular sector. (Green are starting points, and red curve is forest boundary). Note that xy increments approaching 0 for infinite $M$.}
\end{figure}

\section{Connection to Moser's worm problem}

\subsection{Introduction}
Moser's Worm Problem is another well-known open problem proposed by Leo Moser in 1966 \cite{worm wiki}. The problem seeks to find the smallest possible area of a convex set in the plane that can accommodate any planar curve of unit length. The unit-length curve is often referred to metaphorically as a ``worm", leading to the problem's informal name. Moser's Worm Problem is part of a lineage of ``universal covering problems", which ask for minimal sets that can contain or cover all instances of a geometric family.

While Moser’s Worm Problem to Bellman’s Lost-in-a-Forest Problem are distinct problems, they share a core conceptual structure involving invariance under planar isometries like translations and rotations, and also extremal constructions like smallest covers or shortest paths \cite{forest wiki} \cite{worm wiki}. 

While the problem remains open \cite{Croft2012}, research has narrowed the gap between known lower and upper bounds through a sequence of papers \cite{Norwood2003} \cite{Poole1973} \cite{Norwood1992} \cite{Johnson2004} \cite{Wang2006} \cite{Wetzel2013} \cite{Khandhawit2013} \cite{Som-Am2020}. Panraksa et al. \cite{Panraksa2007} demonstrated that no finite bound on the number of segments in a polychain suffices to solve the problem. Further studies explored sector covers for unit arcs \cite{Panraksa2021}, angleworm covers \cite{Sroysang2008}, and covers for closed curves \cite{Füredi2011}. However, it is still an unsolved and open problem.

\subsection{Connection of Bellman’s lost-in-a-forest problem to Moser's worm problem}
For a region $R\in \mathbb{R}^2$ with area $A$, the solution to Bellman’s lost-in-a-forest problem can find the optimal escape path with length $L$. Then $R$ is a universal cover for that escape path. And ratio $\frac{A}{L^2}$  is an upper bound of Moser's worm problem for a planar curve of unit length \cite{Wetzel2003} \cite{Panraksa2007-1}. Finch’s paper \cite{Finch2004} Theorem 3 states that “The escape length $\beta$ of a bounded forest $F$ is the largest $x$ for which $F$ is a cover of $C_x$.”

\begin{definition}[Upper bound of Moser's worm problem]
Formally, let simply connected closed region $R_{worm} \in \mathbb{R}^2$ with area $A$. Boundary $\partial R_{worm} = F(x, y) = 0$. Length of escape path solved by solution to Bellman’s lost-in-a-forest problem is $L$.

Then we obtain an upper bound of the worm problem. Let the minimum area of a shape that can cover every unit-length curve to be determined as $A_m$, then:
\begin{equation}
A_m \leq \frac{A}{L^2}
\end{equation}
 
Therefore, we can calculate solutions to Bellman’s Lost-in-a-Forest Problem for different shapes and iteratively refine the upper bound of Moser's Worm Problem. Philip Gibbs \cite{Gibbs2016} employed numerical methods to calculate the escape path for isosceles triangles, yielding an improved bound. By extending the calculations to polygons of arbitrary shapes, as the number of edges increases, the results converge toward the optimal value of Moser's Worm Problem. In practice, we need to test many different polygons and obtain new upper bounds.

\end{definition}

\section{Connection to shortest opaque set}

\subsection{Introduction}
An opaque set is a collection of curves or sets in the plane that block all lines of sight across a polygon, circle, or other shape \cite{Croft2012} \cite{Opaquewiki} \cite{Brakke1992}. The concept was introduced by Stefan Mazurkiewicz in 1916 \cite{Mazurkiewicz1916}, and Frederick Bagemihl posed the problem of minimizing their total length in 1959 \cite{Bagemihl1959}. Faber described the problem as finding “the shortest curve that meets all the lines that intersect a convex body” \cite{Faber1986}. Some previous paper discussed the weak form problem that “Curves intersecting certain sets of great-circles on the sphere” \cite{Croft1969}.

For a triangle or any convex polygon, the shortest connected opaque set is its minimum Steiner tree \cite{Faber1986}, associated with the Fermat point. However, without assuming connectivity, the Steiner tree's optimality remains unproven. Izumi demonstrated a slight improvement to the perimeter-halving lower bound for equilateral triangles \cite{Izumi2016}. For a unit square, a shorter, disconnected opaque forest is known, with length $\sqrt {2}+\tfrac {1}{2}\sqrt {6}$, comprising the Steiner tree of three vertices and a line segment connecting the fourth vertex to the center. This result, proposed by Jones in 1962 and 1964 \cite{Jones1962} \cite{Jones1964}, is optimal among forests with two components \cite{Kawohl1997} \cite{Kawohl2000} but remains conjectural for the general case \cite{Kawamura2019}. Other advancements include improvements to the perimeter-halving lower bound \cite{Kawamura2019} and studies on minimum opaque manifolds \cite{Asimov2008} and minimum opaque covers for polygonal regions \cite{Provan2012}. Dumitrescu et al. proposed linear-time approximation algorithms for convex polygons \cite{Dumitrescu2015}. However, the shortest opaque sets for unit squares and circles remain unsolved.

\subsection{Connection of Bellman’s lost-in-a-forest problem to shortest opaque set}
As ``opaque set is the shortest curve that meets all the lines that meet a convex body" \cite{Faber1986}. And connectivity is not necessary. People have found some unconnected opaque set with length $\sqrt {2}+\tfrac {1}{2}\sqrt {6}$ for unit square, and approximately 4.7998 for unit circle, although they were conjectured to be optimal.

\begin{definition}[2D region]
We define a two-dimensional region $R_O\in \mathbb{R}^2$ to find the opaque set.
\end{definition}

\begin{definition}[Inside points]
We define M points $(r_{kx},r_{ky}) $ inside the region $R_O$ as $(r_{kx},r_{ky}) \in R_O, \forall k \in \{1,2,\ldots,M\} $.
\end{definition}

\begin{definition}[Line orientation]
We define orientation $\alpha=\frac{2 \pi i}{N}, i \in \{0,1,2,\ldots,N-1\}$ to be orientation at the angle of x-axis.
\end{definition}

\begin{definition}[All lines within a region]
All lines within a region can be represented as:
\[
(x - r_{kx}) \cos \left( \frac{2 \pi i}{N} \right) + (y - r_{ky}) \sin \left( \frac{2 \pi i}{N} \right) = 0
\]
\[
\forall i \in \{0,1,2,\ldots,N-1\}, \forall k \in \{1,2,\ldots,M\}, M, N \to \infty
\]

When $M$ points are evenly distributed in the region like grid points, and when $N$ orientations are evenly distributed in $[0, 2\pi]$. 

A similar formulation of Eq (17) can be used to define the problem of finding the shortest opaque set.
\end{definition}

\begin{definition}[Shortest length of opaque set with one curve for M points and N orientations.]
For one curve, it is also similar to TSP or Hamiltonian path problem. It must pass through point $(x_{a_i},y_{a_i})$ on each line defined above to intersect or block all the lines. The starting point of this curve is not the origin $(0,0)$ but among these points. The problem can be written as
\begin{equation}
\begin{split}
\text{minimize } &\sum_{i=2}^{MN} \sqrt{(x_{a_i} - x_{a_{i-1}})^2 + (y_{a_i} - y_{a_{i-1}})^2}\\
\text{subject to:} &(x_{a_h} - r_{kx}) \cos \left( \frac{2 \pi i}{N} \right) + (y_{a_h} - r_{ky}) \sin \left( \frac{2 \pi i}{N} \right) = 0,\\
&\forall i \in \{0,1,2,\ldots,N-1\}, \forall k \in \{1,2,\ldots,M\}, \exists h \in \{1,2,\ldots,MN\}.
\end{split}
\end{equation}
\end{definition}

The optimized variable in the above equation is a permutation of $\{1, 2, ..., MN\}$, and $2MN$ continuous variables.

Therefore, when $M$ points are evenly distributed in the region like grid points with $x$ and $y$ direction increment $\Delta \to 0$ and $N$ orientations are evenly distributed in $[0, 2\pi]$ with increment $\theta \to 0$, when $M$ and $N$ are infinitely large, it yields the solution to opaque set with one curve for any region.

As opaque set can be unconnected with several curves, it is a variant of multiple TSP or multiple path cover problem. In the unconnected case, we allow multiple curves that collectively meet all lines.

\begin{definition}[$p$ subsets for a set $\{1, 2, ..., MN\}$ ]
We define $p \in \mathbb{N}$ and need to partition the set of $MN$ combinations into $p$ subsets then construct the opaque set. 
Let
\[
\bigcup_{i=1}^p P_i = \{1,2,\ldots,MN\}, P_i \neq \emptyset, P_i \cap P_j = \emptyset, \forall i,j \in \{1,2,\ldots,p\}.
\]
\end{definition}

Each subset $P_i$ is visited by one curve in opaque set.

\begin{definition}[Shortest length of opaque set with $p$ curves for $M$ points and $N$ orientations.]  

Then the problem can be written as:
\begin{equation}
\begin{split}
\text{minimize } &\sum_{j=1}^p \sum_{i \in P_j} \sqrt{(x_{a_i} - x_{a_{i-1}})^2 + (y_{a_i} - y_{a_{i-1}})^2}\\
\text{subject to:} &(x_{a_h} - r_{kx}) \cos \left( \frac{2 \pi i}{N} \right) + (y_{a_h} - r_{ky}) \sin \left( \frac{2 \pi i}{N} \right) = 0,\\
&\forall i \in \{0,1,2,\ldots,N-1\}, \forall k \in \{1,2,\ldots,M\}, \exists h \in \{1,2,\ldots,MN\}.
\end{split}
\end{equation}
\end{definition}
Therefore, when $M$ points are evenly distributed in the region like grid points with $x$ and $y$ direction increment $\Delta \to 0$ and $N$ orientations are evenly distributed in $[0, 2\pi]$ with increment $\theta \to 0$, when $M$ and $N$ are infinitely large, it yields the solution to opaque set with $p$ curves for any region. (Similar to Theorem 3)

The next two subsection will present the formulation to minimize opaque sets for the unit square and unit circle.

\subsection{Shortest opaque sets for the unit square}
Assume opaque set with $p$ curves for $M$ points and $N$ orientations. Let $M = m^2, m \in \mathbb{N}$. The solution formulation to find shortest opaque set for the unit square is:
\begin{equation}
\begin{split}
\text{minimize } &\sum_{j=1}^p \sum_{i \in P_j} \sqrt{(x_{a_i} - x_{a_{i-1}})^2 + (y_{a_i} - y_{a_{i-1}})^2}\\
\text{subject to:} &\left( x_{a_h} - \frac{u}{m} \right) \cos \left( \frac{2 \pi i}{N} \right) + \left( y_{a_h} - \frac{v}{m} \right) \sin \left( \frac{2 \pi i}{N} \right) = 0,\\
&\forall i \in \{0,1,2,\ldots,N-1\}, \forall u \in \{0,1,2,\ldots,m-1\}, \forall v \in \{0,1,2,\ldots,m-1\}, \\
&\exists h \in \{1,2,\ldots,MN\},
\end{split}
\end{equation}
where $\left( \frac{u}{m}, \frac{v}{m} \right) = (r_{ux}, r_{vy})$ are $M$ points within the unit square as uniform grid.

By solving this optimization, results such as the optimal opaque set for the unit square can be obtained. For $M, N \to \infty$, this approach yields highly accurate results.

\begin{boldremark}
Similar methods can extend to other polygons.
\end{boldremark}

\subsection{Shortest opaque sets for the unit circle (Beam Detector)}
The shortest length of the optimal opaque sets for the unit circle has been called the beam detection constant.

Assume opaque set with $p$ curves for $M$ points and $N$ orientations. Let $M = m^2, m \in \mathbb{N}$. The solution formulation to find shortest opaque set for the unit circle is:
\begin{equation}
\begin{split}
\text{minimize } &\sum_{j=1}^p \sum_{i \in P_j} \sqrt{(x_{a_i} - x_{a_{i-1}})^2 + (y_{a_i} - y_{a_{i-1}})^2}\\
\text{subject to:} &\left( x_{a_h} - \sqrt{\frac{u}{m}} \cos \frac{2 \pi v}{m} \right) \cos \left( \frac{2 \pi i}{N} \right)+ \left( y_{a_h} - \sqrt{\frac{u}{m}} \sin \frac{2 \pi v}{m} \right) \sin \left( \frac{2 \pi i}{N} \right) = 0,\\
&\forall i \in \{0,2,\ldots,N-1\}, \forall u \in \{0,1,2,\ldots,m-1\}, \forall v \in \{0,1,2,\ldots,m-1\},\\
& \exists h \in \{1,2,\ldots,MN\},
\end{split}
\end{equation}

where $\left(\sqrt{\frac{u}{m}} \cos \frac{2 \pi v}{m}, \sqrt{\frac{u}{m}} \sin \frac{2 \pi v}{m}\right) = (r_{uvx}, r_{uvy})$ are $M$ points within the unit circle as radial grid.

By solving this optimization, results such as the optimal opaque set for the unit circle can be obtained. For $M, N \to \infty$, this approach yields highly accurate results.

\subsection{Examples and results}
For the unit circle, some papers \cite{Joris1980} \cite{Croft1969} simplify it by only considering the intersection with all tangent lines around the circle. 

Assuming that the order of intersection points is $\left \{ 0,1,...,N-1 \right \}$

Discrete formulation of the minimization problem can be written as:
\begin{equation}
\begin{split}
\text{minimize } &\sum_{i=1}^{N-1} \sqrt{(x_i - x_{i-1})^2 + (y_i - y_{i-1})^2} \\
\text{subject to:} &x_i \cos \frac{2\pi i}{N} + y_i \sin \frac{2\pi i}{N} - 1 = 0, \forall i = 0, 1, 2, \ldots, N-1.
\end{split}
\end{equation}

By solving above optimization, the results are shown in Figure 15:
\begin{figure}[H]
  \centering
  \includegraphics[width=5cm]{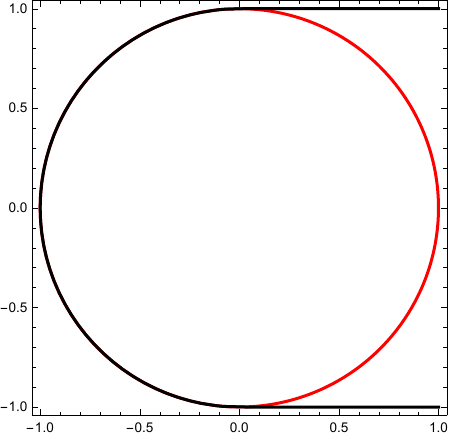}
  \caption{Opaque set with one curve for unit circle (Black curve is opaque set, and red curve is forest boundary).}
\end{figure}

\begin{boldremark}
The results in Figure 15 were obtained and mentioned in Fig 3 in \cite{Joris1980}, Fig 2 in \cite{Croft1969}, and Fig 16 in \cite{Faber1984}, but with geometric methods.
\end{boldremark}

Similarly, for an ellipse with the minor axis of 0.5 and the major axis of 1, if assuming only consider the tangent lines around, discrete formulation of the minimization problem can be written as:
\begin{equation}
\begin{split}
\text{minimize } &\sum_{i=1}^{N-1} \sqrt{(x_i - x_{i-1})^2 + (y_i - y_{i-1})^2} \\
\text{subject to:} &x_i \cos \left ( \frac{2\pi i}{N}+ \varphi\right )  +2 y_i \sin \left ( \frac{2\pi i}{N}+\varphi   \right )  - 1 = 0, \forall i = 0, 1, 2, \ldots, N-1.
\end{split}
\end{equation}
where $\varphi \in \left [ 0,2\pi \right ] $ is the starting angle of the tangent line to define the opening direction.

By solving above optimization, the results are shown in Figure 16:
\begin{figure}[H]
  \centering
  (a)\includegraphics[width=4cm]{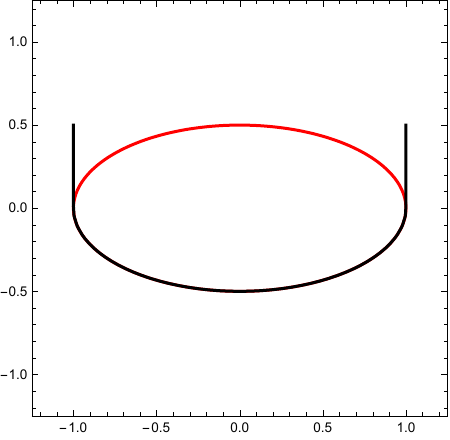}
  (b)\includegraphics[width=4cm]{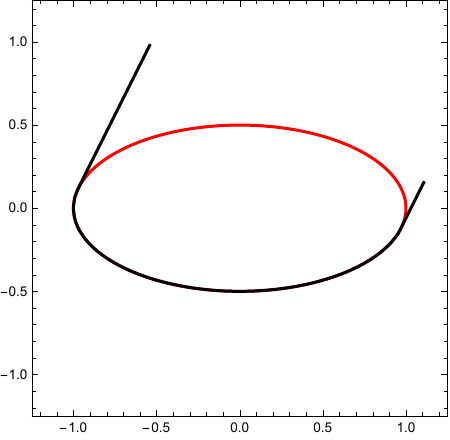}
  (c)\includegraphics[width=4cm]{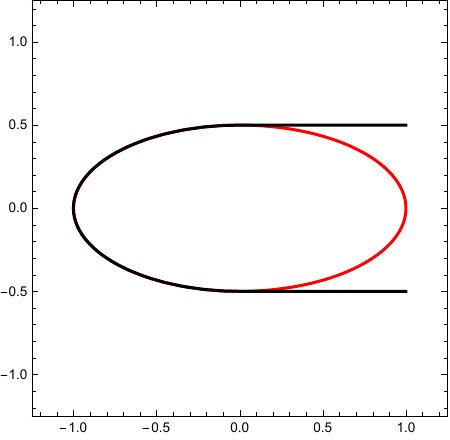}
  \caption{Opaque set with one curve for ellipse with different opening directions (a) $\pi/2$, (b) $\pi/4$, and (c) $0$ (Black curve is opaque set, and red curve is forest boundary).}
\end{figure}

Mathematica code for obtaining results in Figure 15 and 16 is in Appendix IV.

\begin{boldremark}
The results in Figure 16 were mentioned in page 40 of Ref \cite{Croft2012}, as the results for different directions of the ellipse are different.
\end{boldremark}

\section{More variations and extensions}
The solution in this paper is general, which can be used to optimally escape/search from one region defined by point, curve, surface, solid to another region, and optimal cover/fit of curves.

The contribution of this paper is providing an exact mathematical formulation and general solution to Bellman’s Lost-in-a-Forest Problem. The approach sets a framework for computational exploration and approximation for Moser’s Worm Problem and shortest opaque set. For the final solution, discretization become arbitrarily fine approaching a continuum. One more question remains: In order to obtain a reasonable approximation of the solution to the problems, how large do $M,N$ need to be?

There are more variations and extensions of the lost-in-a-forest problem, worm problem, and opaque set. For example, spiral search is conjectured to be the best search strategy or path for the entire plane without knowing the distance and direction of the target, but it has not been proven \cite{Zalgaller2005} \cite{Finch2005} \cite{Langetepe2010}. There are more cover and fit problems \cite{Wetzel2003} that can be solved by the method in this paper. The methods may provide ideas for universal covering problem. The forest, worm, and opaque set problems can also be extended to three dimensions \cite{Ghomi2021}. And the method can even be extended to non-Euclidean geometry and spherical and hyperbolic spaces.

\subsection{Lost-in-a-Forest Problem with closed path, and Worm problem for closed curve}
The original Bellman’s Lost-in-a-forest Problem does not require the escape path return to the starting point (closed). But in practical application of path exploration and planning, this may be required and useful. In addition, the closed path returning to the starting point may also be related to the variant of Moser’s worm problem for “closed worms”.

\begin{definition}[Length of closed escape path return to the starting point for Weak Form I]
Recall Definition 3.14, the length of closed escape path return to the starting point for Weak Form I can be written as
\begin{equation}
L_N=\sqrt{x_{a_0}^2 + y_{a_0}^2} + \sum_{i=1}^{N-1} \sqrt{(x_{a_i} - x_{a_{i-1}})^2 + (y_{a_i} - y_{a_{i-1}})^2} + \sqrt{x_{a_{N-1}}^2 + y_{a_{N-1}}^2}\\
\end{equation}
\end{definition}

\begin{definition}[Length of closed escape path return to the starting point for Weak Form II]
Recall Definition 3.21, the length of closed escape path return to the starting point for Weak Form II can be written as
\begin{equation}
L_{MN}=\sqrt{x_{a_1}^2 + y_{a_1}^2} + \sum_{i=2}^{MN} \sqrt{(x_{a_i} - x_{a_{i-1}})^2 + (y_{a_i} - y_{a_{i-1}})^2} +\sqrt{x_{a_{MN}}^2 + y_{a_{MN}}^2}
\end{equation}
\end{definition}

The other constraint definitions in the optimization for Weak Form I and II remain unchanged.

\begin{boldremark}
It is also possible to use TSP with “returning to the origin city” to define the optimization for closed path return to the starting point.
\end{boldremark}

There are some examples and results below. We use the same assumptions for the order of escape points as Section 4.2. And we only modified the objective function using Eq (76) above. By solving the optimization, the results are shown in Figures 17 and 18. Mathematica code is in Appendix V.

\begin{figure}[H]
  \centering
  (a)\includegraphics[width=4.5cm]{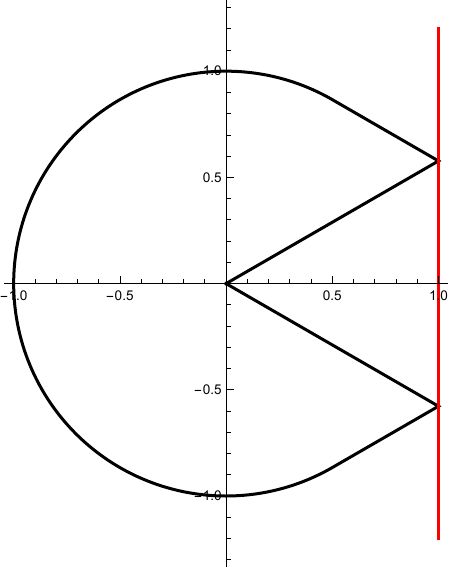}
  (b)\includegraphics[width=4.5cm]{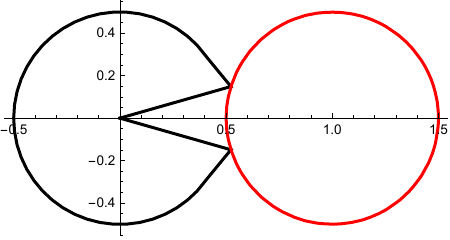}
  (c)\includegraphics[width=4.5cm]{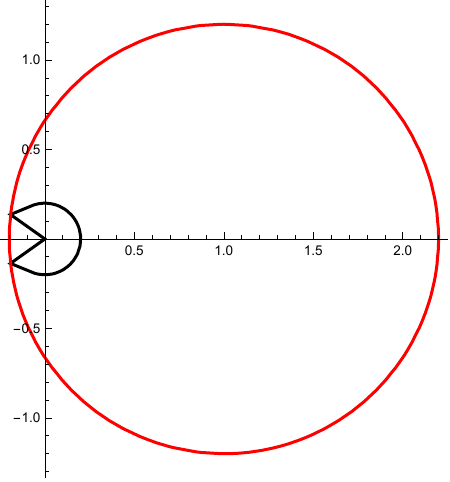}
  (d)\includegraphics[width=4.5cm]{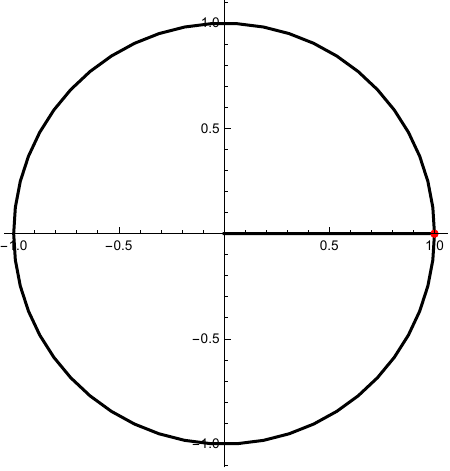}
  (e)\includegraphics[width=4.5cm]{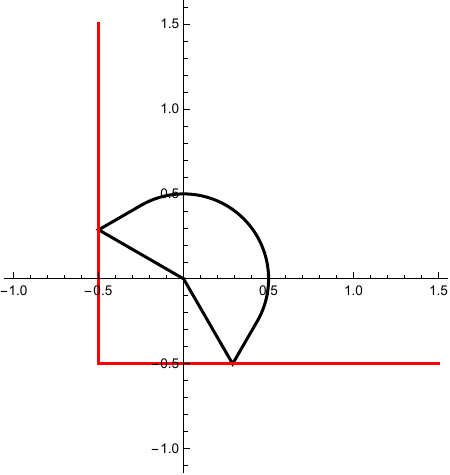}
  (f)\includegraphics[width=4.5cm]{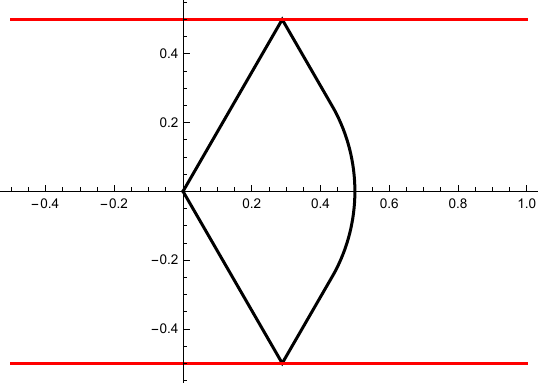}
  \caption{Results of search for (a) one line with unit distance; (b) a circle from exterior; (c) a circle from interior; (d) one point with given distance 1; (e) two perpendicular lines; (f) two parallel lines from the middle; and the path is return to the starting point (closed) (Black curve is escape path, and red curve is forest boundary).}
\end{figure}

\begin{figure}[H]
  \centering
  (a)\includegraphics[width=5cm]{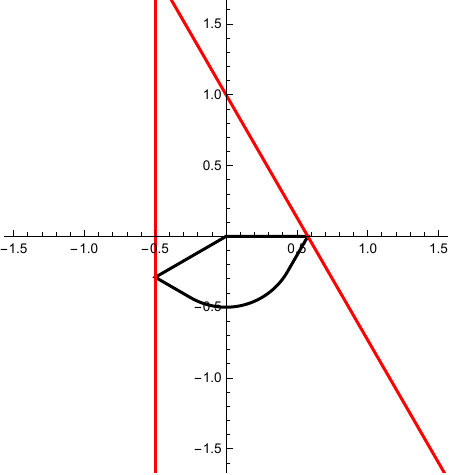}
  (b)\includegraphics[width=5cm]{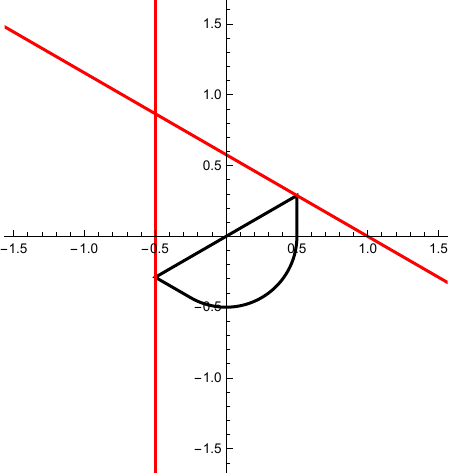}
  (c)\includegraphics[width=5cm]{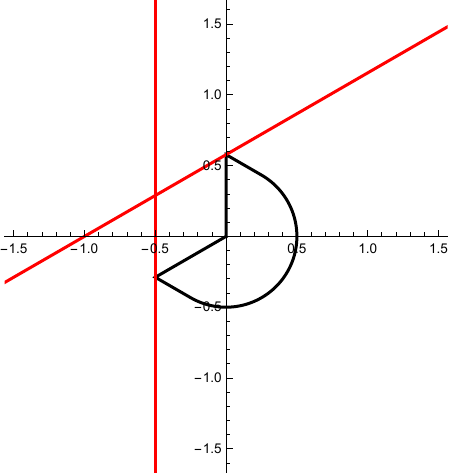}
  (d)\includegraphics[width=5cm]{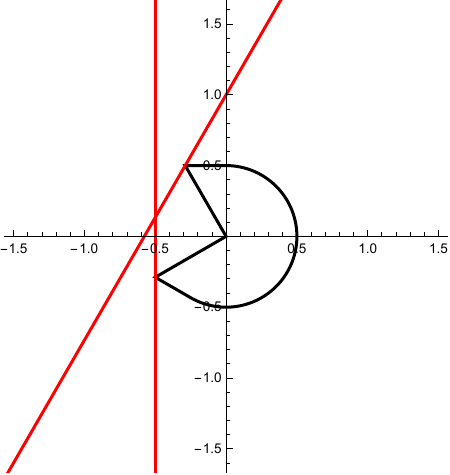}
  \caption{Results of search for two lines from the angle bisector in the middle at certain angle (a)$\pi/6$ (b)$\pi/3$ (c)$2\pi/3$ (d)$5\pi/6$ and the path is return to the starting point (closed) (Black curve is escape path, and red curve is forest boundary).}
\end{figure}
  
\begin{boldremark}
Result in Figure 17 (a) was obtained as Fig 2(a) in Ref \cite{Melzak2007}.
\end{boldremark}

And there are more descriptions and review in the book \cite{Brass2005}  for the variant of Moser’s worm problem for “closed worms” as follows:
\\

“What is the smallest area of a set in the plane that contains a congruent copy of every closed curve of length at most one? … What is the smallest area of a set in the plane that contains a congruent copy of every closed curve of length at most one? … What is the minimum area of a universal cover C for all plane convex sets of unit perimeter? … It would also be interesting to find the minimum perimeter of a universal cover for this class. It may be easier to solve these problems for “closed worms” than to answer Moser’s original question for “open” ones.”
\\

The above modification of the objective function for the optimal path to return to the starting point, and the method in Section 7 should be applicable to the above variant of the Moser’s worm problem for closed curve.

\subsection{Lost-in-a-Forest Problem in 3D}
The simplified lost-in-the-forest-problem for plane of unit distance in 3D was described in previous papers as:
\\

“What is the shortest closed orbit a satellite may take to inspect the entire surface of a round asteroid?” \cite{Ghomi2021}
\\

“What is the shortest curve in $\mathbb{R}^3$ that starts at the origin and intersects all planes at distance 1 from the origin. This precise question arises in the following contrived but amusing setting. Imagine that a space shop has already landed on the surface of a spherical asteroid, of unit radius. What is the shortest path that the spaceship can take in order to survey the entire surface of the asteroid?” \cite{Chan2003}
\\

“Perhaps more interestingly, what is the shortest continuous curve which intersects each plane cutting the sphere, that is the shortest curve with convex hull containing the sphere? Is it half a revolution of a helix with a straight line segment at each end?” \cite{Croft2012}
\\

“What is the form of the shortest curve $C$ going outside the unit sphere $S$ in $\mathbb R^3$ such that passing along $C$ we can see all points of $S$ from outside? How will the form of $C$ change if we require that $C$ have one of its (or both) endpoints on $S$?” \cite{Zalgaller2023}
\\

In this subsection we assume the three-dimensional xyz coordinate system.

\begin{definition}[Forest boundary in 3D]
Consider a three-dimensional body $\mathbb{R}^3$ and a surface defined implicitly by an equation:
\begin{equation}
F(x, y, z) = 0
\end{equation}

where $F: \mathbb{R}^3 \to \mathbb{R}$ is a given function. The surface $F(x, y ,z) = 0$ represents the ``forest boundary" aiming to find or escape in 3D.
\end{definition}

\begin{definition}[Escape path in 3D]
We aim to determine a continuous path as

\begin{equation}
f\left [ x(t),y(t), z(t)  \right ] 
\end{equation}
\end{definition}

In the original problem, we also do not know the starting point and orientation.

Similar to Weak Form I, we can assume start from origin $(0,0,0)$. And search the rotated spatial plane with unit distance.

\begin{definition}[Orientation in 3D]
We define arbitrary unknown initial heading orientation in 3D with $M, N \in \mathbb{N}$:
\begin{equation}
\frac{\pi i}{N} \in \left [ 0,\pi \right ] , i \in \left \{ 0, 1, ..., N-1 \right \} 
\end{equation}
\begin{equation}
\frac{2\pi k}{M} \in \left [ 0,2\pi \right ] , k \in \left \{ 1,2, ..., M \right \}
\end{equation}
\end{definition}

\begin{definition}[Escape points in 3D]
The escape points $\left ( x_i,y_i, z_i \right ) \in \mathbb{R}^3$ on each rotated forest boundary are denoted as 
\begin{equation*}
(x_0,y_0, z_0),(x_1,y_1,z_1),...,(x_{MN-1},y_{MN-1},z_{MN-1})
\end{equation*}
\end{definition}

\begin{definition}[Order of escape points on escape path]
We still only know escape points are on the escape path. However, the order is unknown. 

Let:
\begin{equation}
\mathbf{a} = \{a_0, a_1, \ldots, a_{MN-1}\}, \mathbf{a} \in S_{MN},
\end{equation}
where $S_{MN}$ is the symmetric group of all permutations of $\{0, 1, 2, \ldots, MN-1\}$.

$\mathbf{a}$ is the order in which the escape path visits the escape points. Then the order of escape points with different spatial orientations has $MN$ possible combinations.
\end{definition}

\begin{definition}[Length of escape path in 3D]
We define the length of escape path in 3D computed as the sum of the Euclidean distances from the origin to the first escape point and between the subsequent $MN$ escape points in order in 3D.
\end{definition}
$L_{MN}$ can be written as
\begin{equation}
L_{MN}=\sqrt{x_{a_1}^2 + y_{a_1}^2+ z_{a_1}^2} + \sum_{i=2}^{MN} \sqrt{(x_{a_i} - x_{a_{i-1}})^2 + (y_{a_i} - y_{a_{i-1}})^2+ (z_{a_i} - z_{a_{i-1}})^2},\\
\end{equation}

\begin{definition}[Solution to Bellman’s lost-in-a-forest problem to search the spatial plane with unit distance in 3D]
The objective is to find the minimal length $L_{MN}$ for $MN$ orientations to search the spatial plane with unit distance in 3D. Eqs (17) and (57) can be extended in 3D to minimize length of curve with escape point on every tangent plane of unit ball as:
\begin{equation}
\begin{aligned}
\text{minimize} & \sqrt{x_{a_1}^2 + y_{a_1}^2+ z_{a_1}^2} + \sum_{i=2}^{MN} \sqrt{(x_{a_i} - x_{a_{i-1}})^2 + (y_{a_i} - y_{a_{i-1}})^2+ (z_{a_i} - z_{a_{i-1}})^2}, \\
\text{subject to:} & 
x_{a_h} \sin \frac{\pi i}{N}\cos \frac{2\pi k}{M} + y_{a_h} \sin \frac{\pi i}{N}\sin \frac{2\pi k}{M}+ z_{a_h} \cos \frac{\pi i}{N} = 1, \\
& \forall i \in \{0, 1, 2, \ldots, N-1\}, \, \forall k \in \{1, 2, \ldots, M\}, \, \exists h \in \{1, 2, \ldots, MN\}.
\end{aligned}
\end{equation}

The optimized variables in the above equation are permutation from 1 to $MN$, and $3MN$ continuous variables.
\end{definition}

Therefore, more general cases in 3D can also be generalized and solved.

\subsection{Worm Problem and Opaque set in 3D}
Similarly, worm problem in 3D can be the minimum volume that can enclose a three-dimensional worm with unit length, described in previous papers as ``smallest sleeping bag for a baby snake" with more information \cite{Wetzel2003}, and “There are obvious 3-dimensional analogs of these problems, but very little seems to be known about the least volume of covering sets of certain types.” \cite{Croft2012}

Methods in Section 7.2 can be extended to three dimensions with above section ``Lost-in-a-Forest Problem in 3D".

Opaque set problem in 3D can be defined as ``the smallest surface area can block all space straight line in any space body". Eqs (70) and (71) can also be extended for shortest opaque set in 3D.

\section{Solve minimization problems}

The problems described in Sections 4–9 involve binary optimization. As $N$ and $M$ grows, the problem rapidly becomes NP-hard, and exact solution may fail to produce in reasonable time. The solutions for large $M,N$ are conceptually correct but computationally NP-hard to obtain global optimality. Global optimization is required but is extremely challenging. Solutions can be approached using:

If $M$ and $N$ are small, generate all possible permutations and enumerate all possible solutions. However, the complexity is factorial $O(n!)$ .

Use a branch-and-bound algorithm to systematically explore subsets of solutions. This method guarantees finding the global optimum by pruning suboptimal branches based on bounds.

Apply the Held-Karp algorithm, which computes the global optimal solution using dynamic programming. This approach has $O(n^22^n)$  complexity.

Use relaxations such as the Minimum Spanning Tree, or Linear Programming relaxation of TSP.

Other heuristics methods.

Although we obtained the formulations, it is still challenging to obtain the exact results and values for Moser's worm problem and Beam Detector due to high computational complexity.

\section{Appendices-Mathematica notebooks}
The appendices provide Mathematica notebooks that contain detailed equation derivations and numerical calculations presented in this paper.

We run notebooks in Mathematica 14.1. 

Please note that numerical global optimization function “NMinimize” in Mathematica might not always guarantee the global optimum for large number of variables.
\\
\\
\textbf{Appendix I} Mathematica notebook for continuous formulation solved by calculus of variations in Section 4.2.1.
\\
\\
\textbf{Appendix II} Mathematica notebook for discrete formulation solved in Section 4.2.1-4.2.10
\\
\\
\textbf{Appendix III} Mathematica notebook for discrete formulation solved in Section 6.2.1
\\
\\
\textbf{Appendix IV} Mathematica notebook for discrete formulation solved in Section 8.5
\\
\\
\textbf{Appendix V} Mathematica notebook for discrete formulation solved in Section 9.1
\\
\\

~\\
College of Engineering and Computer Science, University of Central Florida, Orlando, FL, USA

Email: \underline{zhipeng.deng@ucf.edu}



\section*{Supplementary material}

# Appendix 1

*In[ ]:=* `(*Weak form I, given forest boundary of any shape,`
`search/excape from one known point but orientation is unknown*)`

*In[ ]:=* `(*continuous formulation*)`

*In[ ]:=* `(*solve by calculus of variation: EL equation with Lagrange multiplier for constraints*)`

`(*4.2.1 search for one line/half plane with unit distance*)`

`(*Adding Lagrange multipliers to the objective functional*)`

*In[ ]:=* `H = √(x'[t]² + y'[t]²) + λ[t] * (x[t] * Cos[t] + y[t] * Sin[t] - 1);`

*In[ ]:=* `Needs["VariationalMethods`"]`

`(*Euler-Lagrange Equations*)`

*In[ ]:=* `EulerEquations[H, {x[t], y[t], λ[t]}, t]`

*Out[ ]=*
$$\left\{ \text{Cos}[t] \, \lambda[t] + \frac{y'[t] \, (-y'[t] \, x''[t] + x'[t] \, y''[t])}{\left(x'[t]^2 + y'[t]^2\right)^{3/2}} == 0, \right.$$
$$\left. \text{Sin}[t] \, \lambda[t] + \frac{x'[t] \, (y'[t] \, x''[t] - x'[t] \, y''[t])}{\left(x'[t]^2 + y'[t]^2\right)^{3/2}} == 0, \, -1 + \text{Cos}[t] \, x[t] + \text{Sin}[t] \, y[t] == 0 \right\}$$

`(*Substitution and simplification*)`

*In[ ]:=* `((D[H, x[t]] - D[D[H, x'[t]], t]) * Sin[t] - (D[H, y[t]] - D[D[H, y'[t]], t]) * Cos[t]) /.`
`{x'[t] → D[(1 - Sin[t] y[t])/Cos[t], t], x''[t] → D[(1 - Sin[t] y[t])/Cos[t], {t, 2}]} // Simplify`

*Out[ ]=*
$$\left( \text{Sec}[t] \, (\text{Sin}[t] - y[t]) \right.$$
$$\left. \left( \text{Sec}[t] \, (-3 + \text{Cos}[2\,t] + 4 \, \text{Sin}[t] \, y[t]) \, y'[t] + 4 \, y'[t]^2 + 2 \, (\text{Sin}[t] - y[t]) \, y''[t] \right) \right) \Big/$$
$$\left( 2 \left( \text{Tan}[t]^2 + \text{Sec}[t]^2 \, y[t]^2 - 2 \, \text{Tan}[t] \, y[t] \, (\text{Sec}[t] - y'[t]) - 2 \, \text{Sin}[t] \, \text{Tan}[t] \, y'[t] + y'[t]^2 \right) \right.$$
$$\left. \sqrt{y'[t]^2 + \left( \text{Sec}[t]^2 \, y[t] + \text{Tan}[t] \, (-\text{Sec}[t] + y'[t]) \right)^2} \right)$$

*In[ ]:=* `(*It is clear that the solution contain two parts`
`(Sin[t]-y[t])`
`and`
`(Sec[t] (-3+Cos[2 t]+4 Sin[t] y[t]) y'[t]+4 y'[t]²+2 (Sin[t]-y[t]) y''[t])*)`

*In[ ]:=* `(*numerical solution using Finite Difference method to solve above EL equations*)`

`(*Discretization points*)`

*In[ ]:=* `pts = 10 002;`

*In[ ]:=* `aa = N@Range[0, 2 π, 2 π / (pts - 1)];`
`grid = Flatten[Outer[List, aa], 1];`

`(*Finite difference method*)`

*In[ ]:=*  **opt** = **"DifferenceOrder"** → **4**;
**dfda** = **NDSolve`FiniteDifferenceDerivative**[{1}, {**aa**}, **opt**]["DifferentiationMatrix"];
**d2fda2** = **NDSolve`FiniteDifferenceDerivative**[{2}, {**aa**}, **opt**]["DifferentiationMatrix"];

*In[ ]:=*  **nL** = **Length**[**grid**];
**varsY** = **Table**[**Unique**[**y$**], {**aa**}];

(*Approximate segmentation for EL euqation*)

*In[ ]:=*  **eqns** = $\bigl($**Sec**[**aa**] (-3 + **Cos**[2 **aa**] + 4 **Sin**[**aa**] **varsY**) (**dfda.varsY**) + 4 (**dfda.varsY**)$^2$ +
2 (**Sin**[**aa**] - **varsY**) (**d2fda2.varsY**)$\bigr)$ * (**Boole**[0 ≤ # < 1.0472] & /@ **aa**) +
(**Sin**[**aa**] - **varsY**) * (**Boole**[1.0472 ≤ # < 3 π / 2] & /@ **aa**) +
$\bigl($**Sec**[**aa**] (-3 + **Cos**[2 **aa**] + 4 **Sin**[**aa**] **varsY**) (**dfda.varsY**) + 4 (**dfda.varsY**)$^2$ +
2 (**Sin**[**aa**] - **varsY**) (**d2fda2.varsY**)$\bigr)$ * (**Boole**[3 π / 2 ≤ # ≤ 2 π] & /@ **aa**);

*In[ ]:=*  **eqns**⟦**1**⟧ = **varsY**⟦**1**⟧ - **0.57735**; (*Approximate boundary value*)
**eqns**⟦**nL**⟧ = **varsY**⟦**nL**⟧ + **1**;

*In[ ]:=*  **govEq** = **eqns**;
**allVars** = **varsY**;
**root** = **FindRoot**[**govEq**, **Thread**[{**allVars**, **Flatten**[{**Table**[1, **nL**]}]}]];

··· **FindRoot**: The line search decreased the step size to within tolerance specified by AccuracyGoal and PrecisionGoal but was
unable to find a sufficient decrease in the merit function. You may need more than MachinePrecision digits of working
precision to meet these tolerances. ⓘ

(*Solve discretized equations*)

*In[ ]:=*  **sol** = **allVars** /. **root**;
{**ysoltemp**} = **Partition**[**sol**, **nL**];
**ysol** = **Interpolation**[**Transpose**[{**aa**, **ysoltemp**}]];

*In[ ]:=*  **Plot**$\Bigl[\Bigl\{\dfrac{1 - \text{Sin}[a]\ \text{ysol}[a]}{\text{Cos}[a]},$ **ysol**[**a**]$\Bigr\}$, {**a**, **0**, **2 π**},
**PerformanceGoal** → **"Speed"**, **PlotRange** → **All**, **PlotLegends** → {**"x[a]"**, **"y[a]"**}]

*Out[ ]=*

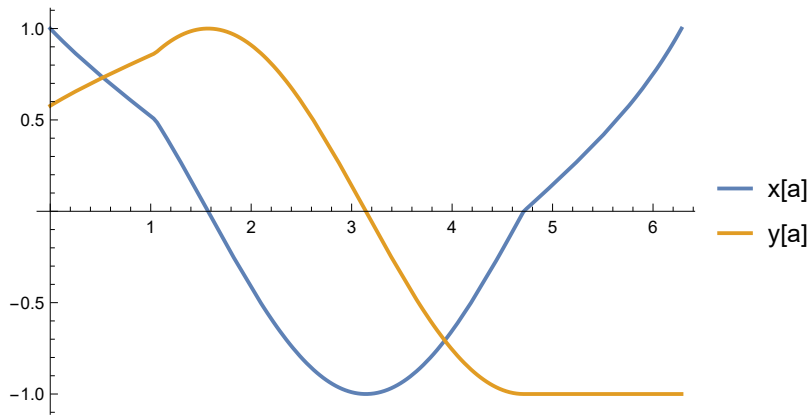

(*Boundary and path*)

$In[*]:=$ Show$\Big[$ListLinePlot$\Big[\Big\{\{0, 0\}, \Big\{\frac{1 - Sin[0] \ ysol[0]}{Cos[0]}, ysol[0]\Big\}\Big\}$, AspectRatio → Automatic,

PerformanceGoal → "Speed", PlotRange → All, PlotStyle → Black$\Big]$,

ParametricPlot$\Big[\Big\{\frac{1 - Sin[a] \ ysol[a]}{Cos[a]}, ysol[a]\Big\}, \{a, 0, 2 \pi\}$, PerformanceGoal → "Speed",

PlotRange → {{-1.5, 1.5}, {-1.5, 1.5}}, PlotStyle → Black$\Big]$,

ContourPlot$[x - 1 == 0, \{x, -1.2, 1.2\}, \{y, -1.2, 1.2\}$, ContourStyle → Red,

PerformanceGoal → "Speed"], PlotRange → All, ImageSize → 300$\Big]$

$Out[*]=$

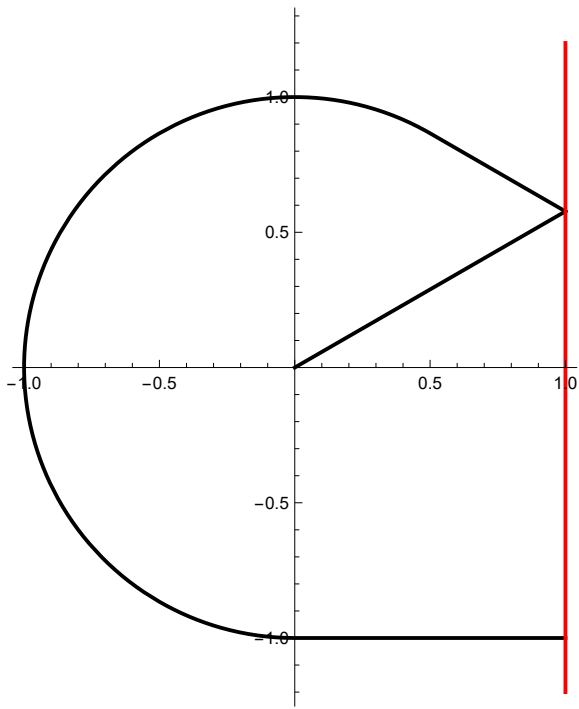

(*length of path*)

$In[*]:=$ $\sqrt{\Big(\frac{1 - Sin[0] \ ysol[0]}{Cos[0]}\Big)^2 + ysol[0]^2}$ + NIntegrate$\Big[\sqrt{D\Big[\frac{1 - Sin[a] \ ysol[a]}{Cos[a]}, a\Big]^2 + ysol'[a]^2}$,

{a, 0, 1.5 π}, WorkingPrecision → 50, MaxRecursion → 30, AccuracyGoal → 8$\Big]$ + 1

⋯ **NIntegrate:** Numerical integration converging too slowly; suspect one of the following: singularity, value of the integration is 0, highly oscillatory integrand, or WorkingPrecision too small. 🛈

$Out[*]=$

6.39724

# Appendix 2

```
(*WEAK FORM I, given forest boundary of any shape,
search/excape from one known point but orientation is unknown*)

(*discrete formulation*)

(*4.2.1 search a line with unit distance*)

n = 5000;(*Number of discretization points*)

(*Numerical minimization for global optimal*)

sol = NMinimize[{√(x[0]² + y[0]²) + Sum[√((x[i] - x[i-1])² + (y[i] - y[i-1])²), {i, 1, n}]
    (*objective*), Table[(x[i] * Cos[(2π/n) * i] + y[i] * Sin[(2π/n) * i] - 1) == 0, {i, 0, n}]
    (*constraints*)}, Flatten[Table[{x[i], y[i]}, {i, 0, n}]]];

(*length of path*)
```

*In[ ]:=* `sol[[1]]`

*Out[ ]=*
    6.39724

*In[ ]:=* `RealDigits[%]`

*Out[ ]=*
    {{6, 3, 9, 7, 2, 4, 1, 9, 9, 7, 3, 3, 0, 7, 9, 3}, 1}

```
In[*]:= Show[ListLinePlot[Prepend[Table[{x[i], y[i]}, {i, 0, n}] /. sol[[2]], {0, 0}],
        AspectRatio → Automatic, PerformanceGoal → "Speed", PlotRange → All, PlotStyle → Black],
      ContourPlot[x - 1 ⩵ 0, {x, -1.2, 1.2}, {y, -1.2, 1.2}, ContourStyle → Red,
        PerformanceGoal → "Speed"], PlotRange → All, ImageSize → 300]
```

Out[*]=

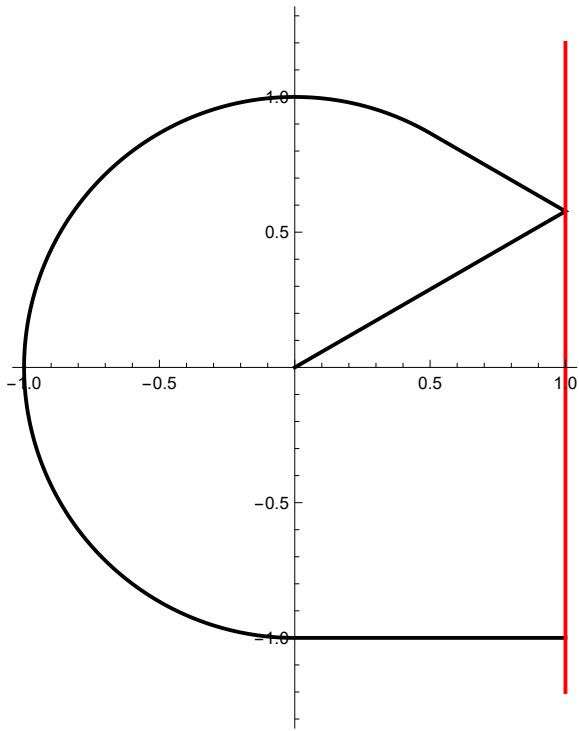

```
In[*]:= NotebookSave[];
```

(*4.2.2 search a circle from exterior*)

(*r/s=1*)

n = 100; (*Number of discretization points*)

(*Numerical minimization for global optimal*)

```
In[*]:= sol = NMinimize[{√(x[0]² + y[0]²) + Sum[√((x[i] - x[i-1])² + (y[i] - y[i-1])²), {i, 1, n}],
        Table[(x[i] - Cos[2π/n * i])² + (y[i] - Sin[2π/n * i])² ⩵ 0.5^2, {i, 0, n}]},
        Flatten[Table[{x[i], y[i]}, {i, 0, n}]]];
```

(*length*)

```
In[*]:= sol[[1]]
```

Out[*]=
      3.40002

*In[ ]:=* `RealDigits[%]`

*Out[ ]=*
{{3, 4, 0, 0, 0, 1, 6, 0, 1, 3, 6, 9, 7, 3, 8, 0}, 1}

*In[ ]:=* `Show[ListLinePlot[Prepend[Table[{x[i], y[i]}, {i, 0, n}] /. sol[[2]], {0, 0}],`
`    AspectRatio → Automatic, PerformanceGoal → "Speed", PlotRange → All, PlotStyle → Black],`
`   ContourPlot[(x - 1)^2 + (y - 0)^2 ⩵ 0.5^2, {x, -1, 2}, {y, -1.5, 1.5}, ContourStyle → Red],`
`   PlotRange → All, ImageSize → 300]`

*Out[ ]=*

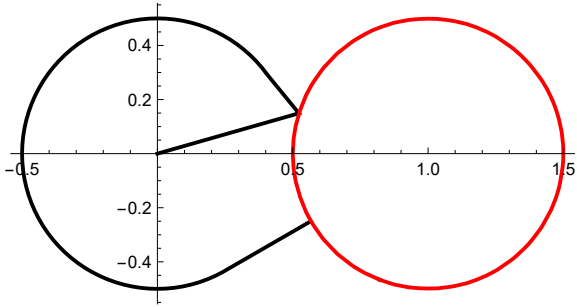

*In[ ]:=* `NotebookSave[];`

(*4.2.3 search a circle from interior*)

*In[ ]:=* (*r/s=0.2*)

`n = 100;` (*Number of discretization points*)

(*Numerical minimization for global optimal*)

*In[ ]:=* `sol = NMinimize[{√(x[0]^2 + y[0]^2) + Sum[√((x[i] - x[i-1])^2 + (y[i] - y[i-1])^2), {i, 1, n}]`,
`      Table[(x[i] - Cos[(2π)/n * i])^2 + (y[i] - Sin[(2π)/n * i])^2 ⩵ 1.2^2, {i, 0, n}]},`
`      Flatten[Table[{x[i], y[i]}, {i, 0, n}]]];`

(*length*)

*In[ ]:=* `sol[[1]]`

*Out[ ]=*
1.24738

*In[ ]:=* `RealDigits[%]`

*Out[ ]=*
{{1, 2, 4, 7, 3, 8, 4, 7, 6, 6, 2, 2, 4, 6, 9, 0}, 1}

*In[•]:=* `Show[ListLinePlot[Prepend[Table[{x[i], y[i]}, {i, 0, n}] /. sol[[2]], {0, 0}],`
    `AspectRatio → Automatic, PerformanceGoal → "Speed", PlotRange → All, PlotStyle → Black],`
    `ContourPlot[(x - 1)^2 + (y - 0)^2 ⩵ 1.2^2, {x, -1, 3}, {y, -1.5, 1.5}, ContourStyle → Red],`
    `PlotRange → All, ImageSize → 300]`

*Out[•]=*

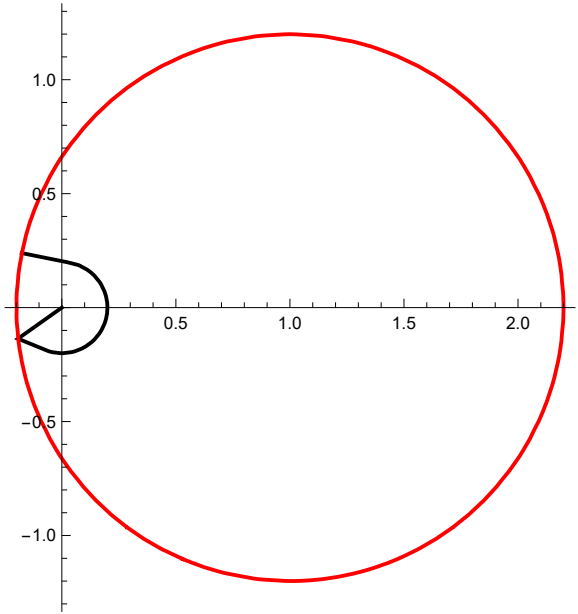

*In[•]:=* `NotebookSave[];`

---

    `(*4.2.4 search a circle from interior when shortest path is not unique*)`

*In[•]:=* `(*r/s=0.333454*)`

    `n = 100;(*Number of discretization points*)`

    `(*Numerical minimization for global optimal*)`

*In[•]:=* `sol = NMinimize[{{√(x[0]^2 + y[0]^2) + Sum[√((x[i] - x[i-1])^2 + (y[i] - y[i-1])^2), {i, 1, n}],`
    `Table[(x[i] - Cos[2π/n * i])^2 + (y[i] - Sin[2π/n * i])^2 ⩵ 1.500272^2, {i, 0, n}]},`
    `Flatten[Table[{x[i], y[i]}, {i, 0, n}]]];`

    `(*length*)`

*In[•]:=* `sol[[1]]`

*Out[•]=*
    `2.99957`

*In[•]:=* `RealDigits[%]`

*Out[•]=*
    `{{2, 9, 9, 9, 5, 6, 6, 6, 2, 3, 2, 8, 0, 1, 2, 4}, 1}`

*In[•]:=* `(*diameter 2r=3*)`

*In[∘]:=* `Show[ListLinePlot[Prepend[Table[{x[i], y[i]}, {i, 0, n}] /. sol⟦2⟧, {0, 0}],`
 `AspectRatio → Automatic, PerformanceGoal → "Speed", PlotRange → All, PlotStyle → Black],`
 `ContourPlot[(x - 1)^2 + (y - 0)^2 ⩵ 1.500272^2, {x, -1, 3}, {y, -2, 2}, ContourStyle → Red],`
 `PlotRange → All, ImageSize → 300]`

*Out[∘]=*

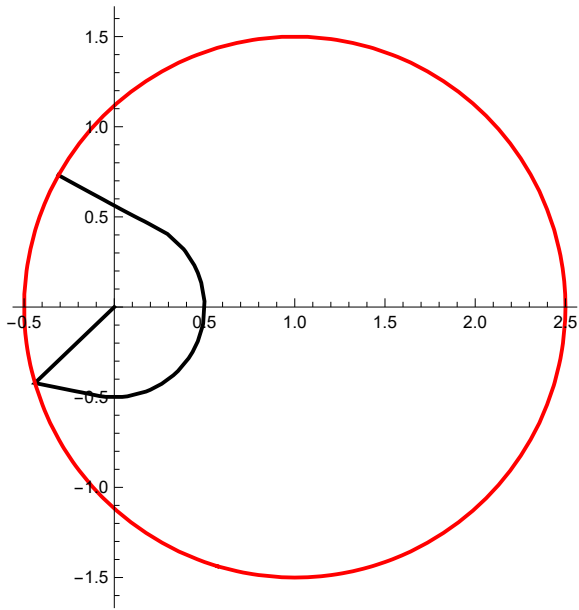

*In[∘]:=* `NotebookSave[];`

---

`(*4.2.5 search for one point with given distance 1*)`

`nn = 100;` `(*Number of discretization points*)`

`(*Numerical minimization for global optimal*)`

*In[∘]:=* `sol = NMinimize[{{√(x[0]² + y[0]²) + Sum[√((x[i] - x[i-1])² + (y[i] - y[i-1])²), {i, 1, nn}],`
 `Table[(x[i] - Cos[(2 π)/nn * i])² + (y[i] - Sin[(2 π)/nn * i])² ⩵ 0.0², {i, 0, nn}]},`
 `Flatten[Table[{x[i], y[i]}, {i, 0, nn}]]];`

`(*length*)`

*In[∘]:=* `sol⟦1⟧`

*Out[∘]=*
 `7.28215`

*In[∘]:=* `RealDigits[%]`

*Out[∘]=*
 `{{7, 2, 8, 2, 1, 5, 1, 8, 0, 2, 3, 9, 4, 6, 0, 6}, 1}`

*In[∘]:=* `(*the length should be 2π+1*)`

*In[ ]:=* `Show[ListLinePlot[Prepend[Table[{x[i], y[i]}, {i, 0, nn}] /. sol[[2]], {0, 0}],`
`    AspectRatio → Automatic, PerformanceGoal → "Speed", PlotRange → All, PlotStyle → Black],`
`  Graphics[{Red, Disk[{1, 0}, 0.02]}], PlotRange → All, ImageSize → 300]`

*Out[ ]=*

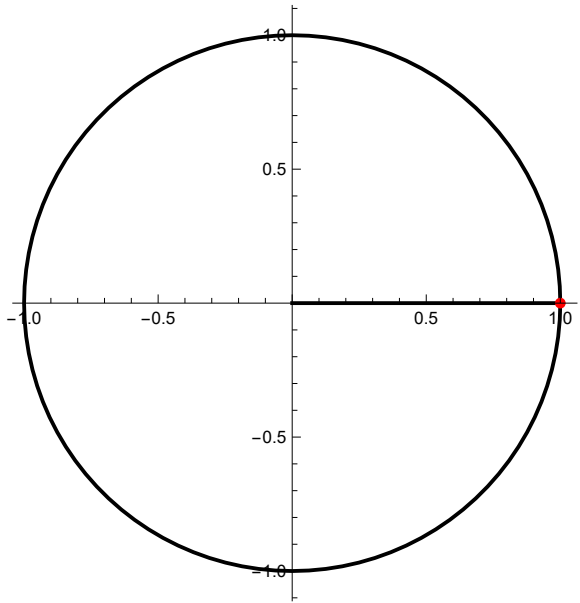

*In[ ]:=* `NotebookSave[];`

`(*4.2.6 Search for a shape from interior with distance`
` significantly short at certain angle, when shortest path is not unique*)`

`(*only search for circle*)`

`n = 30; (*Number of discretization points*)`

`(*Numerical minimization for global optimal*)`

*In[ ]:=* `sol = NMinimize[{{√(x[0]² + y[0]²) + Sum[√((x[i] - x[i-1])² + (y[i] - y[i-1])²), {i, 1, n}],`

`    Table[((x[i] * Cos[ (2π/n) * i] - y[i] * Sin[ (2π/n) * i] - 0)² +`

`      (x[i] * Sin[ (2π/n) * i] + y[i] * Cos[ (2π/n) * i] - 0)² - 1²) == 0,`

`    {i, 0, n}]}, Flatten[Table[{x[i], y[i]}, {i, 0, n}]]];`

*In[ ]:=* `sol[[1]]`

*Out[ ]=*

`1.`

```
In[ ]:=  Show[ListLinePlot[Prepend[Table[{x[i], y[i]}, {i, 0, n}] /. sol[[2]], {0, 0}],
          AspectRatio → Automatic, PerformanceGoal → "Speed", PlotRange → All, PlotStyle → Black],
         ContourPlot[(x - 0) ^2 + (y - 0) ^2 ⩵ 1^2, {x, -3, 3}, {y, -2, 2}, ContourStyle → Red],
         ContourPlot[y ⩵ 0, {x, -1 - (1. - 1/(1 + 2 π)), -1 + (1. - 1/(1 + 2 π))}, {y, -2, 2},
          ContourStyle → Red], PlotRange → All, ImageSize → 300, AspectRatio → Automatic]
```

Out[ ]=

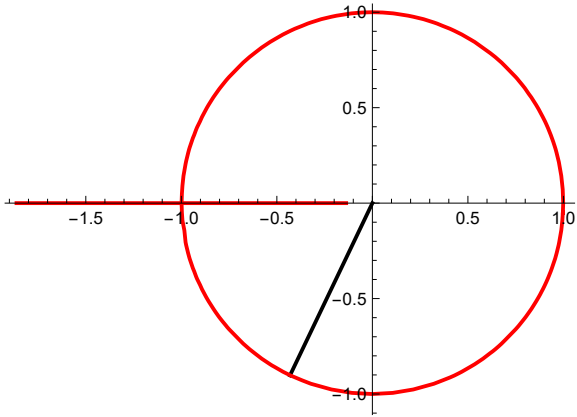

```
         (*only search for line segment*)

In[ ]:=  sol = NMinimize[{√(x[0]^2 + y[0]^2) + Sum[√((x[i] - x[i-1])^2 + (y[i] - y[i-1])^2), {i, 1, n}],
          Table[(x[i] * Sin[2 π/n * i] + y[i] * Cos[2 π/n * i] - 0)^2 +
            ((x[i] * Cos[2 π/n * i] - y[i] * Sin[2 π/n * i] - 1)^2 - (1. - 1/(1 + 2 π))^2)^2 ⩵ 0,
           {i, 0, n}]}, Flatten[Table[{x[i], y[i]}, {i, 0, n}]]];

In[ ]:=  sol[[1]]

Out[ ]=
         0.99674
```

```
In[*]:= Show[ListLinePlot[Prepend[Table[{x[i], y[i]}, {i, 0, n}] /. sol[[2]], {0, 0}],
        AspectRatio → Automatic, PerformanceGoal → "Speed", PlotRange → All, PlotStyle → Black],
       ContourPlot[(x - 0)^2 + (y - 0)^2 == 1^2, {x, -3, 3}, {y, -2, 2}, ContourStyle → Red],
       ContourPlot[y == 0, {x, -1 - (1. - 1/(1 + 2 π)), -1 + (1. - 1/(1 + 2 π))}, {y, -2, 2},
         ContourStyle → Red], PlotRange → All, ImageSize → 300, AspectRatio → Automatic]
```

Out[*]=

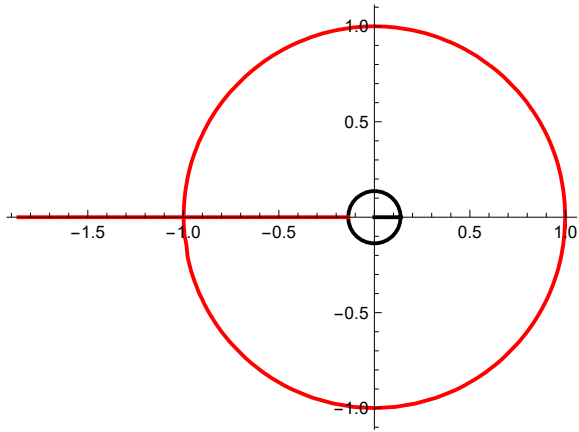

```
(*search for circle and line segment, *)
```

```
In[*]:= sol = NMinimize[{{√(x[0]^2 + y[0]^2) + Sum[√((x[i] - x[i-1])^2 + (y[i] - y[i-1])^2), {i, 1, n}],
        Table[((x[i] * Cos[2π/n * i] - y[i] * Sin[2π/n * i] - 0)^2 + (x[i] * Sin[2π/n * i] +
            y[i] * Cos[2π/n * i] - 0)^2 - 1^2) * ((x[i] * Sin[2π/n * i] + y[i] * Cos[2π/n * i] - 0)^2 +
            ((x[i] * Cos[2π/n * i] - y[i] * Sin[2π/n * i] - 1)^2 - (1. - 1/(1 + 2π))^2)^2) == 0,
         {i, 0, n}]}, Flatten[Table[{x[i], y[i]}, {i, 0, n}]]];
```

```
In[*]:= sol[[1]]
```

Out[*]=

```
0.998375
```

*In[ ]:=* **Show[ListLinePlot[Prepend[Table[{x[i], y[i]}, {i, 0, n}] /. sol⟦2⟧, {0, 0}],**
    **AspectRatio → Automatic, PerformanceGoal → "Speed", PlotRange → All, PlotStyle → Black],**
    **ContourPlot[(x - 0)^2 + (y - 0)^2 == 1^2, {x, -3, 3}, {y, -2, 2}, ContourStyle → Red],**
    **ContourPlot$\left[y == 0, \left\{x, -1 - \left(1. - \dfrac{1}{1 + 2\pi}\right), -1 + \left(1. - \dfrac{1}{1 + 2\pi}\right)\right\}, \{y, -2, 2\},\right.$**
    **ContourStyle → Red], PlotRange → All, ImageSize → 300, AspectRatio → Automatic]**

*Out[ ]=*

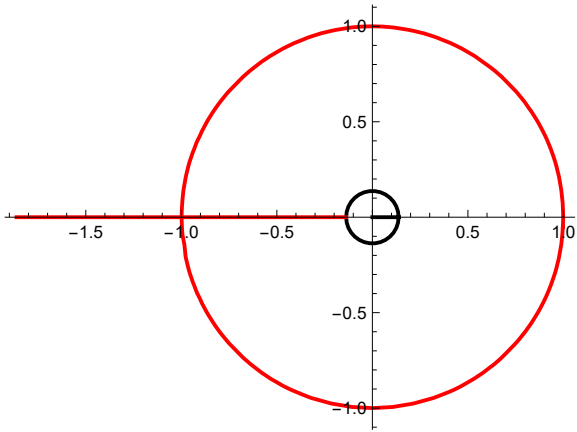

(∗shortest path is not unique, as length of both path is 1∗)

(∗4.2.7 two perpendicular lines with 0.5 distance∗)

**nn = 500;** (∗Number of discretization points∗)

(∗Numerical minimization for global optimal∗)

*In[ ]:=* **sol = ConvexOptimization[**
    $\sqrt{\text{x[0]}^2 + \text{y[0]}^2} + \text{Sum}\left[\sqrt{(\text{x[i]} - \text{x[i-1]})^2 + (\text{y[i]} - \text{y[i-1]})^2}, \{i, 1, \text{nn}\}\right],$
    **Table$\left[\left(-\text{x[i]} * \text{Sin}\left[\dfrac{2\pi}{\text{nn}} * \dfrac{3}{4} * i\right] + \text{y[i]} * \text{Cos}\left[\dfrac{2\pi}{\text{nn}} * \dfrac{3}{4} * i\right] + 0.5\right) == 0, \{i, 0, \text{nn}\}\right],$**
    **Flatten[Table[{x[i], y[i]}, {i, 0, nn}]], MaxIterations → 20 000];**

(∗length∗)

*In[ ]:=* $\sqrt{\text{x[0]}^2 + \text{y[0]}^2} + \text{Sum}\left[\sqrt{(\text{x[i]} - \text{x[i-1]})^2 + (\text{y[i]} - \text{y[i-1]})^2}, \{i, 1, \text{nn}\}\right]$ **/. sol**

*Out[ ]=*
    2.41322

*In[ ]:=* `Show[ListLinePlot[Prepend[Table[{x[i], y[i]}, {i, 0, nn}] /. sol, {0, 0}],`
`   AspectRatio → Automatic, PerformanceGoal → "Speed", PlotRange → All, PlotStyle → Black],`
`  ContourPlot[x + 0.5 ⩵ 0, {x, -1.5, 1.5}, {y, -0.5, 1.5}, ContourStyle → Red],`
`  ContourPlot[y + 0.5 ⩵ 0, {x, -0.5, 1.5}, {y, -1, 1.5}, ContourStyle → Red],`
`  PlotRange → {{-1, 1.5}, {-1, 1.5}}, ImageSize → 300]`

*Out[ ]=*

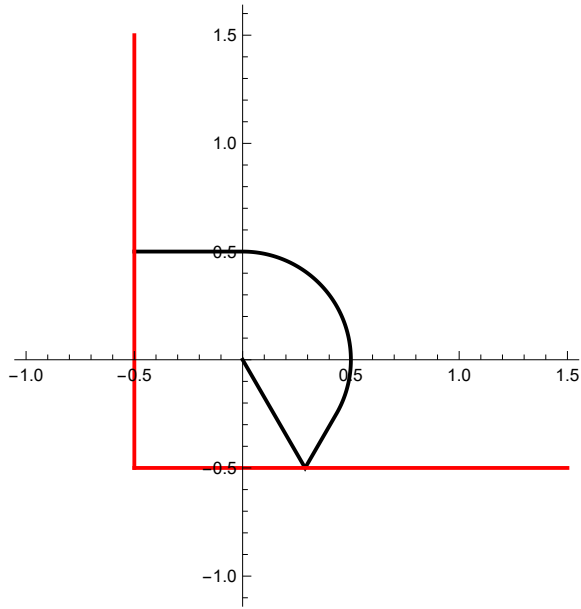

*In[ ]:=* `NotebookSave[];`

---

`(*4.2.8 two parallel lines from the middle*)`

`nn = 500; (*Number of discretization points*)`

`(*Numerical minimization for global optimal*)`

*In[ ]:=* `sol = ConvexOptimization[`
$$\sqrt{x[0]^2 + y[0]^2} + \text{Sum}\Big[\sqrt{(x[i] - x[i-1])^2 + (y[i] - y[i-1])^2}, \{i, 1, nn/2\}\Big],$$
`   Table[x[i] * Cos[`$\frac{2\pi}{nn}$` * i - π / 2] + y[i] * Sin[`$\frac{2\pi}{nn}$` * i - π / 2] ⩵ 0.5, {i, 0, nn / 2}],`
`   Flatten[Table[{x[i], y[i]}, {i, 0, nn / 2}]], MaxIterations → 20000];`

`(*length*)`

*In[ ]:=* $\sqrt{x[0]^2 + y[0]^2} + \text{Sum}\Big[\sqrt{(x[i] - x[i-1])^2 + (y[i] - y[i-1])^2}, \{i, 1, nn/2\}\Big]$ `/. sol`

*Out[ ]=*
`1.62782`

*In[•]:=* `Show[ListLinePlot[Prepend[Table[{x[i], y[i]}, {i, 0, nn / 2}] /. sol, {0, 0}],`
`    AspectRatio → Automatic, PerformanceGoal → "Speed", PlotRange → All, PlotStyle → Black],`
`  ContourPlot[y - 0.5 == 0, {x, -0.5, 1}, {y, -1.5, 1.5}, ContourStyle → Red],`
`  ContourPlot[y + 0.5 == 0, {x, -0.5, 1}, {y, -1.5, 1.5}, ContourStyle → Red], PlotRange → All]`

*Out[•]=*

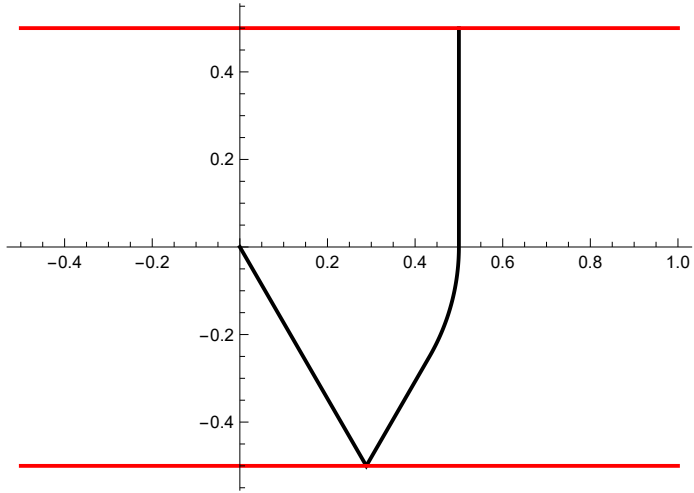

*In[•]:=* `NotebookSave[];`

`(*4.2.9 search for two lines from the angle bisector in the middle*)`

*In[•]:=* `(*π/6 - 30 degree*)`

*In[•]:=* `BB = 1 / 2 + π / 6 / (2 π);`

`nn = 300;(*Number of discretization points*)`

`(*Numerical minimization for global optimal*)`

*In[•]:=* `sol =`
`  NMinimize[{√(x[0]² + y[0]²) + Sum[√((x[i] - x[i-1])² + (y[i] - y[i-1])²), {i, 1, BB * nn}],`
`    Table[(x[i] * Cos[(2 π)/nn * i] + y[i] * Sin[(2 π)/nn * i] + 0.5) == 0, {i, 0, BB * nn}]},`
`    Flatten[Table[{x[i], y[i]}, {i, 0, BB * nn}]]];`

*In[•]:=* `(*length*)`

*In[•]:=* `sol[[1]]`

*Out[•]=*

`1.88961`

*In[ ]:=* Show[ListLinePlot[Prepend[Table[{x[i], y[i]}, {i, 0, BB * nn}] /. sol〚2〛, {0, 0}],
    AspectRatio → Automatic, PerformanceGoal → "Speed", PlotRange → All, PlotStyle → Black],
    ContourPlot[{x * Cos[0] + y * Sin[0] + 0.5 ⩵ 0, x * Cos[0 + π / 6] + y * Sin[0 + π / 6] - 0.5 ⩵ 0},
    {x, -2, 2}, {y, -2, 2}, ContourStyle → Red],
    PlotRange → {{-1.5, 1.5}, {-1.5, 1.5}}, ImageSize → 300]

*Out[ ]=*

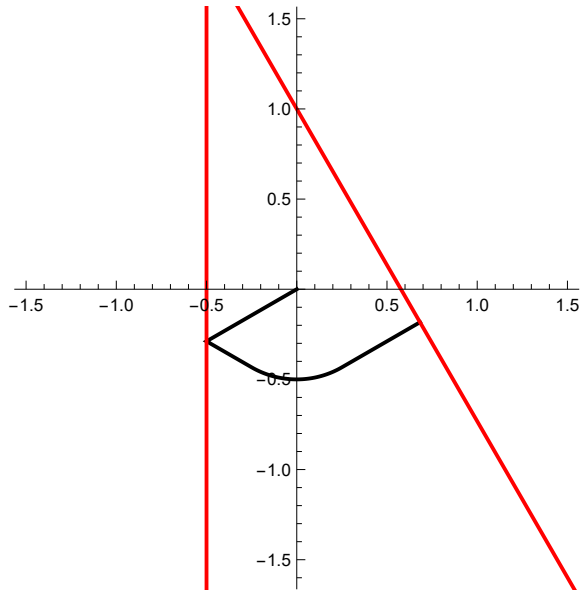

*In[ ]:=* (*π/3 - 60 degree*)

*In[ ]:=* BB = 1 / 2 + π / 3 / (2 π);

    nn = 300; (*Number of discretization points*)

*In[ ]:=* sol =
    NMinimize[{ $\sqrt{x[0]^2 + y[0]^2}$ + Sum[ $\sqrt{(x[i] - x[i-1])^2 + (y[i] - y[i-1])^2}$, {i, 1, BB * nn}],
        Table[ $\left( x[i] * Cos\left[\frac{2\pi}{nn} * i\right] + y[i] * Sin\left[\frac{2\pi}{nn} * i\right] + 0.5 \right)$ ⩵ 0, {i, 0, BB * nn}]},
        Flatten[Table[{x[i], y[i]}, {i, 0, BB * nn}]]];

*In[ ]:=* (*length*)

*In[ ]:=* sol〚1〛

*Out[ ]=*
    2.15141

*In[∘]:=* `Show[ListLinePlot[Prepend[Table[{x[i], y[i]}, {i, 0, BB * nn}] /. sol〚2〛, {0, 0}],`
`    AspectRatio → Automatic, PerformanceGoal → "Speed", PlotRange → All, PlotStyle → Black],`
`  ContourPlot[{x * Cos[0] + y * Sin[0] + 0.5 ⩵ 0, x * Cos[0 + π / 3] + y * Sin[0 + π / 3] - 0.5 ⩵ 0},`
`   {x, -2, 2}, {y, -2, 2}, ContourStyle → Red],`
`  PlotRange → {{-1.5, 1.5}, {-1.5, 1.5}}, ImageSize → 300]`

*Out[∘]=*

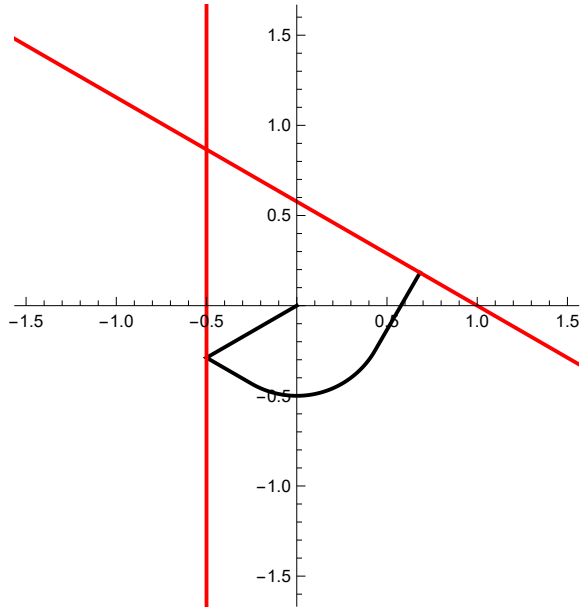

*In[∘]:=* `(*2π/3 - 120 degree*)`

*In[∘]:=* `BB = 1 / 2 + 2 π / 3 / (2 π);`

`nn = 300;` `(*Number of discretization points*)`

*In[∘]:=* `sol =`
`   NMinimize[{√{x[0]² + y[0]²} + Sum[√{(x[i] - x[i - 1])² + (y[i] - y[i - 1])²}, {i, 1, BB * nn}],`
`     Table[(x[i] * Cos[ (2 π)/nn * i] + y[i] * Sin[ (2 π)/nn * i] + 0.5) ⩵ 0, {i, 0, BB * nn}]},`
`     Flatten[Table[{x[i], y[i]}, {i, 0, BB * nn}]]];`

*In[∘]:=* `(*length*)`

*In[∘]:=* `sol〚1〛`

*Out[∘]=*

`2.675`

```
In[ ]:= Show[ListLinePlot[Prepend[Table[{x[i], y[i]}, {i, 0, BB * nn}] /. sol[[2]], {0, 0}],
        AspectRatio → Automatic, PerformanceGoal → "Speed", PlotRange → All, PlotStyle → Black],
      ContourPlot[{x * Cos[0] + y * Sin[0] + 0.5 ⩵ 0, x * Cos[0 + 2 π / 3] + y * Sin[0 + 2 π / 3] - 0.5 ⩵ 0},
        {x, -2, 2}, {y, -2, 2}, ContourStyle → Red],
      PlotRange → {{-1.5, 1.5}, {-1.5, 1.5}}, ImageSize → 300]
```

Out[ ]=

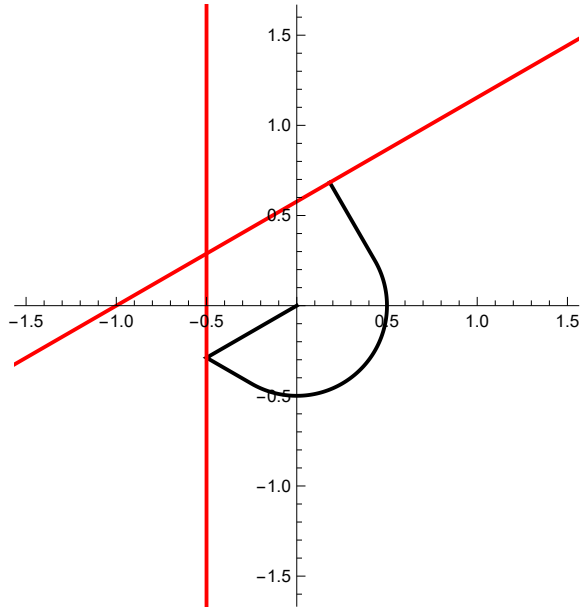

```
In[ ]:= (*5π/6 - 150 degree*)
```

```
In[ ]:= BB = 1 / 2 + 5 π / 6 / (2 π);
```

```
     nn = 300; (*Number of discretization points*)
```

```
In[ ]:= sol =
        NMinimize[{ {√(x[0]² + y[0]²) + Sum[√((x[i] - x[i-1])² + (y[i] - y[i-1])²), {i, 1, BB * nn}]},
          Table[(x[i] * Cos[(2 π / nn) * i] + y[i] * Sin[(2 π / nn) * i] + 0.5) ⩵ 0, {i, 0, BB * nn}]},
          Flatten[Table[{x[i], y[i]}, {i, 0, BB * nn}]]];
```

```
In[ ]:= (*length*)
```

```
In[ ]:= sol[[1]]
```

Out[ ]=

```
     2.93679
```

```
In[•]:= Show[ListLinePlot[Prepend[Table[{x[i], y[i]}, {i, 0, BB * nn}] /. sol[[2]], {0, 0}],
         AspectRatio → Automatic, PerformanceGoal → "Speed", PlotRange → All, PlotStyle → Black],
       ContourPlot[{x * Cos[0] + y * Sin[0] + 0.5 == 0, x * Cos[0 + 5 π / 6] + y * Sin[0 + 5 π / 6] - 0.5 == 0},
         {x, -2, 2}, {y, -2, 2}, ContourStyle → Red],
       PlotRange → {{-1.5, 1.5}, {-1.5, 1.5}}, ImageSize → 300]
```

Out[•]=

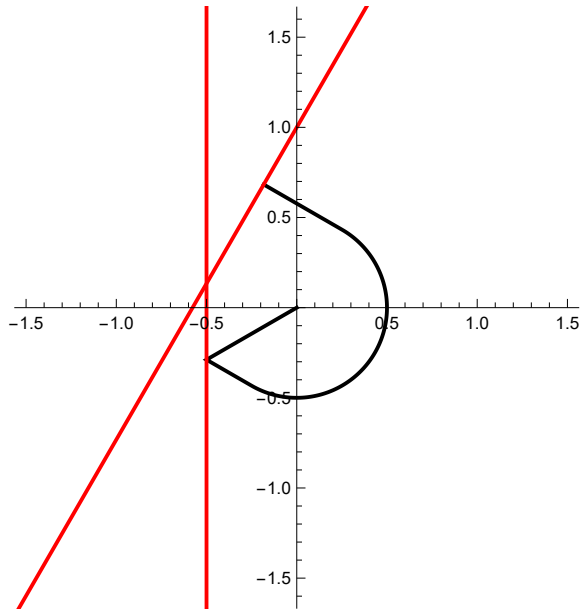

```
In[•]:= NotebookSave[];
```

```
(* 4.2.10 search for two parallel lines with unit distance (unit strip) from one line
   (Zalgaller class 2) *)
```

```
In[•]:= nn = 56; (*Number of discretization points*)
```

```
In[•]:= sol = NMinimize[
         {√(x[0]² + y[0]²) + Sum[√((x[i] - x[i-1])² + (y[i] - y[i-1])²), {i, 1, nn}], Prepend[
           Table[x[i] * Cos[-(π - ArcTan[y[0]/x[0]])/nn * i - π / 2] + y[i] * Sin[-(π - ArcTan[y[0]/x[0]])/nn * i - π / 2] ==
             1, {i, 0, nn}], y[nn] == 0]}, Flatten[Table[{x[i], y[i]}, {i, 0, nn}]]];
```

```
In[•]:= (*length*)
```

```
In[•]:= sol[[1]]
```

Out[•]=

```
2.29747
```

```mathematica
In[*]:= Grid[Partition[Prepend[
    Table[Show[ListLinePlot[Prepend[Table[{x[i], y[i]}, {i, 0, nn}] /. sol〚2〛, {0, 0}],
        PlotRange → All, PlotStyle → Black, PerformanceGoal → "Speed"],
      ContourPlot[{x * Cos[a] + y * Sin[a] - 1 ⩵ 0, x * Cos[a] + y * Sin[a] ⩵ 0}, {x, -1.5, 1.5},
        {y, -1.5, 1.5}, ContourStyle → Red], PlotRange → {{-1.5, 0.6}, {-1, 0.5}},
      AspectRatio → Automatic, PerformanceGoal → "Speed", ImageSize → 200],
     {a, 3 π / 2 - 0.18358, ArcTan[y[0]/x[0]] + π / 2 /. sol〚2〛, -0.18358}],

    Show[ListLinePlot[Prepend[Table[{x[i], y[i]}, {i, 0, nn}] /. sol〚2〛, {0, 0}],
      PlotRange → All, PlotStyle → Black, PerformanceGoal → "Speed"],
     ContourPlot[{x * Cos[3 π / 2] + y * Sin[3 π / 2] - 1 ⩵ 0, x * Cos[3 π / 2] + y * Sin[3 π / 2] ⩵ 0},
       {x, -1.5, 1.5}, {y, -1.5, 1.5}, ContourStyle → Red],
     PlotRange → {{-1.5, 0.6}, {-1, 0.5}}, AspectRatio → Automatic,
     PerformanceGoal → "Speed", ImageSize → 200,
     Epilog → {Text[Style["(x_{N-1},y_{N-1})", 14], {-1, 0.1}],
       Text[Style["(x_0,y_0)", 14], {-0.2, -0.85}], Text[Style["line1", 14], {0.45, 0.1}],
       Text[Style["line2", 14], {0.45, -0.9}]}], 3]]
```

Out[ ]=

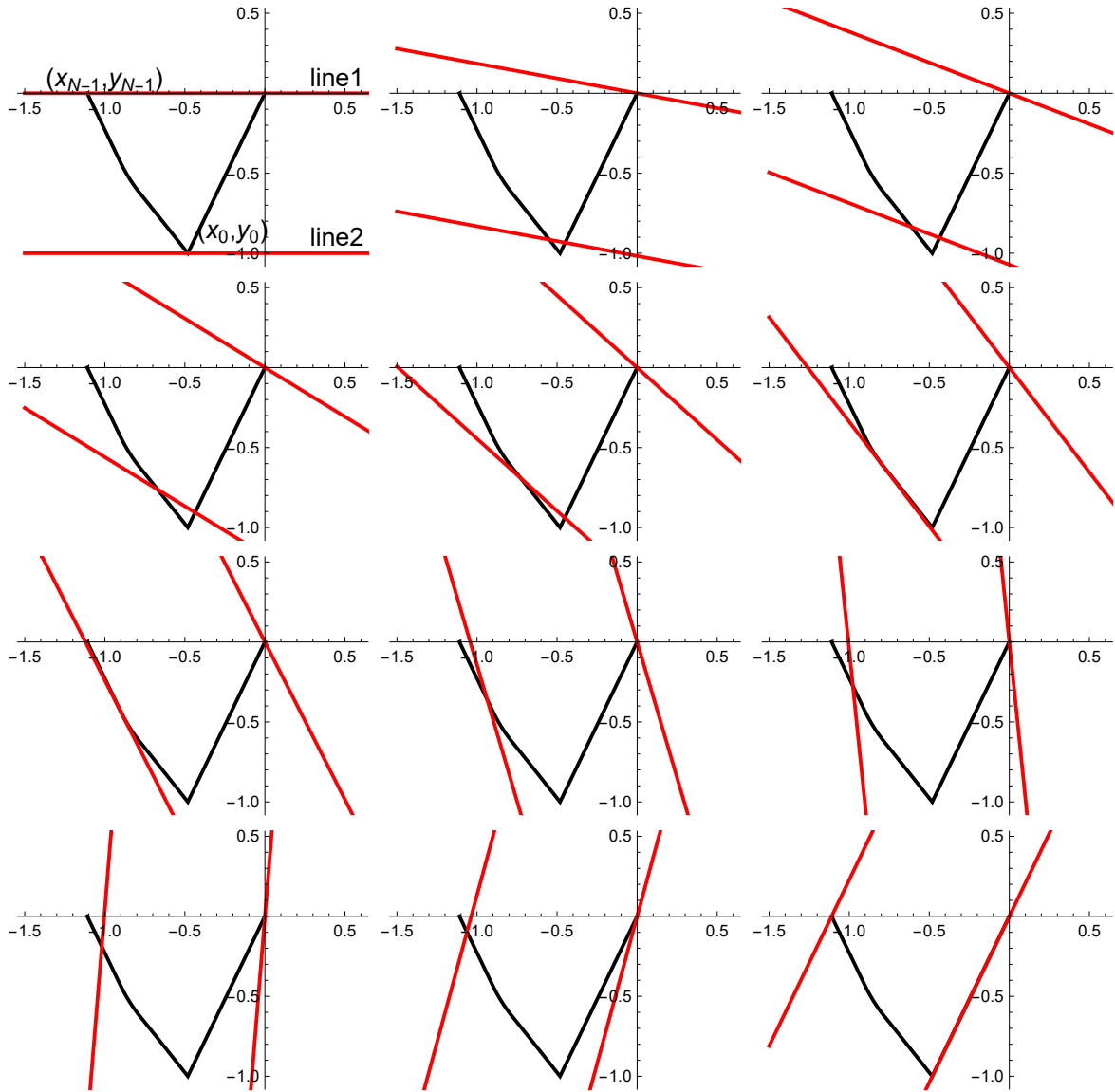

In[ ]:= **NotebookSave[];**

# Appendix 3

(*Section 6.2.1
 solve for unit strip forest and compare with Zalgaller's*)

In[*]:= $\phi$ = ArcSin$\left[\frac{1}{6} + \frac{4}{3} \text{ Sin}\left[\frac{1}{3} * \text{ArcSin}\left[\frac{17}{64}\right]\right]\right]$ (*Parameter in Zalgaller's paper*)

Out[*]=

ArcSin$\left[\frac{1}{6} + \frac{4}{3} \text{ Sin}\left[\frac{1}{3} \text{ ArcSin}\left[\frac{17}{64}\right]\right]\right]$

In[*]:= $\psi$ = ArcTan$\left[\frac{1}{2} \text{ Sec}[\phi]\right]$ (*Parameter in Zalgaller's paper*)

Out[*]=

ArcTan$\left[\dfrac{1}{2 \sqrt{1 - \left(\frac{1}{6} + \frac{4}{3} \text{ Sin}\left[\frac{1}{3} \text{ ArcSin}\left[\frac{17}{64}\right]\right]\right)^2}}\right]$

In[*]:= xx = Sec$[\phi]$ (*Parameter in Zalgaller's paper*)

Out[*]=

$\dfrac{1}{\sqrt{1 - \left(\frac{1}{6} + \frac{4}{3} \text{ Sin}\left[\frac{1}{3} \text{ ArcSin}\left[\frac{17}{64}\right]\right]\right)^2}}$

In[*]:= nn = 12;

In[*]:= mm = 26;

In[*]:= i =.; k =.;

(*Assumed order of points*)

BB = {237, 254, 269, 283, 294, 308, 108, 105, 103, 101, 98, 96, 94, 91, 89, 87, 84, 82, 78, 74,
71, 69, 65, 62, 59, 55, 153, 161, 169, 175, 182, 190, 196, 201, 206, 212, 218, 223,
229, 233, 239, 247, 251, 257, 262, 267, 275, 279, 282, 287, 291, 299, 289, 293, 297,
300, 303, 306, 311, 123, 127, 129, 132, 135, 136, 140, 143, 147, 151, 155, 160, 165,
167, 172, 176, 178, 183, 186, 181, 185, 189, 192, 195, 198, 202, 205, 208, 211, 214,
217, 221, 226, 228, 231, 235, 241, 244, 248, 252, 259, 263, 266, 271, 276, 45, 52, 58,
67, 77, 86, 93, 100, 107, 110, 112, 115, 118, 121, 124, 128, 131, 133, 137, 141, 146,
149, 154, 158, 164, 168, 3, 8, 13, 17, 20, 23, 25, 28, 31, 33, 36, 38, 41, 43, 47, 49,
51, 53, 54, 57, 61, 64, 68, 73, 76, 81, 2, 5, 7, 10, 12, 15, 18, 19, 21, 22, 24, 26,
27, 29, 30, 32, 34, 35, 37, 39, 40, 42, 44, 46, 48, 50, 1, 4, 6, 9, 11, 14, 16, 312,
307, 304, 301, 298, 292, 288, 285, 280, 277, 274, 268, 264, 260, 253, 249, 245, 242,
236, 224, 220, 215, 213, 210, 207, 203, 200, 197, 194, 191, 188, 184, 179, 177, 173,
171, 166, 163, 159, 156, 152, 148, 145, 142, 138, 125, 122, 119, 116, 114, 80, 310,
305, 302, 296, 290, 286, 281, 278, 273, 265, 261, 256, 250, 246, 240, 234, 230, 227,
222, 219, 272, 258, 243, 232, 225, 216, 209, 204, 199, 193, 187, 180, 174, 170, 162,
157, 150, 144, 139, 134, 130, 126, 120, 117, 113, 111, 238, 255, 270, 284, 295, 309,
109, 106, 104, 102, 99, 97, 95, 92, 90, 88, 85, 83, 79, 75, 72, 70, 66, 63, 60, 56};

In[*]:= X = Table[x[h], {h, 1, nn * mm}];

In[*]:= Y = Table[y[h], {h, 1, nn * mm}];

```
(*Numerical minimization may run for very long time to finish*)

res =
  NMinimize[{√(x[1]² + y[1]²) + Sum[√((x[h] - x[h-1])² + (y[h] - y[h-1])²), {h, 2, nn*mm}],

    Flatten[{Table[(x[BB[[(i-1)*mm+k]]] * Cos[2π/(nn-1) * i] + y[BB[[(i-1)*mm+k]]] *

        Sin[2π/(nn-1) * i] + k/(2 mm)) * (x[BB[[(i-1)*mm+k]]] * Cos[2π/(nn-1) * i] +

        y[BB[[(i-1)*mm+k]]] * Sin[2π/(nn-1) * i] + k/(2 mm) - 1) == 0, {i, 1, nn}, {k, 1, mm}],

      Table[0 ≤ x[BB[[(i-1)*mm+k]]] ≤ 1.1, {i, 1, nn}, {k, 1, mm}],
      Table[0 ≤ y[BB[[(i-1)*mm+k]]] ≤ 1, {i, 1, nn}, {k, 1, mm}],
      Table[x[s] ≤ x[s+1], {s, 1, nn*mm-1}]}]},

    {X ∈ Vectors[nn*mm, Reals], Y ∈ Vectors[nn*mm, Reals]},
    AccuracyGoal → 4,
    PrecisionGoal → 3,
    Method → "DifferentialEvolution",
    MaxIterations → 3000
    (*Use Differential Evolution algorithm
      with a maximum of 3000 iterations*)]
```

Out[*]=
```
{2.24853, {x[1] → 0.00514583, x[2] → 0.00565734, x[3] → 0.00970798, x[4] → 0.0103182,
  x[5] → 0.0113529, x[6] → 0.0155856, x[7] → 0.0171309, x[8] → 0.0195447, x[9] → 0.0209309,
  x[10] → 0.0230229, x[11] → 0.0264231, x[12] → 0.0290203, x[13] → 0.0295924,
  x[14] → 0.0320375, x[15] → 0.0351909, x[16] → 0.0377956, x[17] → 0.0398768,
  x[18] → 0.0414588, x[19] → 0.0478358, x[20] → 0.0503304, x[21] → 0.0542553,
  x[22] → 0.0606942, x[23] → 0.0608449, x[24] → 0.067076, x[25] → 0.0712933,
  x[26] → 0.0733886, x[27] → 0.0795627, x[28] → 0.0815477, x[29] → 0.0856743,
  x[30] → 0.0916156, x[31] → 0.0916175, x[32] → 0.0995044, x[33] → 0.103411,
  x[34] → 0.107276, x[35] → 0.115122, x[36] → 0.115209, x[37] → 0.122817, x[38] → 0.126861,
  x[39] → 0.130413, x[40] → 0.137984, x[41] → 0.13843, x[42] → 0.145487, x[43] → 0.149946,
  x[44] → 0.152968, x[45] → 0.152976, x[46] → 0.162722, x[47] → 0.162722, x[48] → 0.173204,
  x[49] → 0.17646, x[50] → 0.183628, x[51] → 0.190167, x[52] → 0.191053, x[53] → 0.203627,
  x[54] → 0.217061, x[55] → 0.226166, x[56] → 0.226166, x[57] → 0.23049, x[58] → 0.230491,
  x[59] → 0.237889, x[60] → 0.237889, x[61] → 0.244798, x[62] → 0.250031, x[63] → 0.250031,
  x[64] → 0.259069, x[65] → 0.262135, x[66] → 0.262135, x[67] → 0.264101, x[68] → 0.273069,
  x[69] → 0.273795, x[70] → 0.273795, x[71] → 0.285292, x[72] → 0.285292, x[73] → 0.286875,
  x[74] → 0.296731, x[75] → 0.296731, x[76] → 0.300627, x[77] → 0.300673, x[78] → 0.308522,
  x[79] → 0.308522, x[80] → 0.310257, x[81] → 0.314833, x[82] → 0.320616, x[83] → 0.320616,
  x[84] → 0.3327, x[85] → 0.3327, x[86] → 0.334563, x[87] → 0.344625, x[88] → 0.344625,
  x[89] → 0.356493, x[90] → 0.356493, x[91] → 0.368352, x[92] → 0.368352, x[93] → 0.369282,
  x[94] → 0.3802, x[95] → 0.3802, x[96] → 0.392058, x[97] → 0.392058, x[98] → 0.403908,
  x[99] → 0.403908, x[100] → 0.404099, x[101] → 0.415673, x[102] → 0.415673,
  x[103] → 0.42746, x[104] → 0.42746, x[105] → 0.439248, x[106] → 0.439248,
  x[107] → 0.439248, x[108] → 0.452739, x[109] → 0.452739, x[110] → 0.468885, x[111] → 0.5,
```

x[112] → 0.501388, x[113] → 0.519231, x[114] → 0.519231, x[115] → 0.526854,
x[116] → 0.531805, x[117] → 0.538462, x[118] → 0.542185, x[119] → 0.545566,
x[120] → 0.557692, x[121] → 0.557692, x[122] → 0.559267, x[123] → 0.563937,
x[124] → 0.572784, x[125] → 0.572784, x[126] → 0.576923, x[127] → 0.578049,
x[128] → 0.586662, x[129] → 0.590093, x[130] → 0.596154, x[131] → 0.600775,
x[132] → 0.602552, x[133] → 0.615069, x[134] → 0.615385, x[135] → 0.615402,
x[136] → 0.627067, x[137] → 0.628935, x[138] → 0.633075, x[139] → 0.634615,
x[140] → 0.638849, x[141] → 0.642851, x[142] → 0.643911, x[143] → 0.650796,
x[144] → 0.653846, x[145] → 0.654866, x[146] → 0.65685, x[147] → 0.662991,
x[148] → 0.66602, x[149] → 0.670955, x[150] → 0.673077, x[151] → 0.675389,
x[152] → 0.677302, x[153] → 0.683181, x[154] → 0.685156, x[155] → 0.688041,
x[156] → 0.688756, x[157] → 0.692308, x[158] → 0.699442, x[159] → 0.700329,
x[160] → 0.700873, x[161] → 0.706302, x[162] → 0.711538, x[163] → 0.711947,
x[164] → 0.713751, x[165] → 0.713755, x[166] → 0.723345, x[167] → 0.726313,
x[168] → 0.727943, x[169] → 0.728924, x[170] → 0.730769, x[171] → 0.734667,
x[172] → 0.738792, x[173] → 0.745955, x[174] → 0.75, x[175] → 0.750958, x[176] → 0.751232,
x[177] → 0.757194, x[178] → 0.763583, x[179] → 0.768398, x[180] → 0.769231,
x[181] → 0.77123, x[182] → 0.77266, x[183] → 0.775884, x[184] → 0.779543,
x[185] → 0.782249, x[186] → 0.788112, x[187] → 0.788462, x[188] → 0.79065,
x[189] → 0.793242, x[190] → 0.794021, x[191] → 0.801709, x[192] → 0.804196,
x[193] → 0.807692, x[194] → 0.812757, x[195] → 0.815156, x[196] → 0.815176,
x[197] → 0.823394, x[198] → 0.825671, x[199] → 0.826923, x[200] → 0.833916,
x[201] → 0.834455, x[202] → 0.836186, x[203] → 0.844423, x[204] → 0.846154,
x[205] → 0.846693, x[206] → 0.853712, x[207] → 0.854937, x[208] → 0.857195,
x[209] → 0.865385, x[210] → 0.865397, x[211] → 0.867403, x[212] → 0.871163,
x[213] → 0.874456, x[214] → 0.876672, x[215] → 0.883434, x[216] → 0.884615,
x[217] → 0.885855, x[218] → 0.885863, x[219] → 0.888334, x[220] → 0.892278,
x[221] → 0.894902, x[222] → 0.898461, x[223] → 0.899974, x[224] → 0.901004,
x[225] → 0.903846, x[226] → 0.903852, x[227] → 0.908143, x[228] → 0.91197,
x[229] → 0.912256, x[230] → 0.917455, x[231] → 0.92006, x[232] → 0.923077,
x[233] → 0.923861, x[234] → 0.926753, x[235] → 0.928183, x[236] → 0.928185,
x[237] → 0.930035, x[238] → 0.930035, x[239] → 0.934272, x[240] → 0.935329,
x[241] → 0.935408, x[242] → 0.935916, x[243] → 0.942308, x[244] → 0.942571,
x[245] → 0.943508, x[246] → 0.943676, x[247] → 0.943683, x[248] → 0.949119,
x[249] → 0.950176, x[250] → 0.951372, x[251] → 0.951782, x[252] → 0.955534,
x[253] → 0.956739, x[254] → 0.956782, x[255] → 0.956782, x[256] → 0.958907,
x[257] → 0.959542, x[258] → 0.961538, x[259] → 0.961676, x[260] → 0.962832,
x[261] → 0.966127, x[262] → 0.966781, x[263] → 0.967645, x[264] → 0.968809,
x[265] → 0.973324, x[266] → 0.973631, x[267] → 0.974004, x[268] → 0.974819,
x[269] → 0.976442, x[270] → 0.976442, x[271] → 0.979891, x[272] → 0.980769,
x[273] → 0.980813, x[274] → 0.981247, x[275] → 0.981731, x[276] → 0.985937,
x[277] → 0.987214, x[278] → 0.987989, x[279] → 0.988872, x[280] → 0.993144,
x[281] → 0.995148, x[282] → 0.996043, x[283] → 0.996425, x[284] → 0.996426,
x[285] → 0.999093, x[286] → 1.0023, x[287] → 1.00318, x[288] → 1.00502, x[289] → 1.0057,
x[290] → 1.00945, x[291] → 1.0103, x[292] → 1.01095, x[293] → 1.01115, x[294] → 1.0155,
x[295] → 1.0155, x[296] → 1.01661, x[297] → 1.01663, x[298] → 1.01691, x[299] → 1.01748,
x[300] → 1.02211, x[301] → 1.02287, x[302] → 1.02379, x[303] → 1.02758, x[304] → 1.02882,

x[305] → 1.03094, x[306] → 1.03305, x[307] → 1.03479, x[308] → 1.0348, x[309] → 1.0348,

x[310] → 1.03711, x[311] → 1.03711, x[312] → 1.03711, y[1] → 0.0186887, y[2] → 0.0205438,

y[3] → 0.0351966, y[4] → 0.0373735, y[5] → 0.0410546, y[6] → 0.0560447, y[7] → 0.0614939,

y[8] → 0.0699546, y[9] → 0.0747047, y[10] → 0.0818343, y[11] → 0.0933435,

y[12] → 0.102084, y[13] → 0.103994, y[14] → 0.111965, y[15] → 0.122183, y[16] → 0.130565,

y[17] → 0.137228, y[18] → 0.142197, y[19] → 0.162118, y[20] → 0.169885,

y[21] → 0.182001, y[22] → 0.201868, y[23] → 0.202335, y[24] → 0.221784, y[25] → 0.23501,

y[26] → 0.24176, y[27] → 0.261856, y[28] → 0.268346, y[29] → 0.282006, y[30] → 0.302304,

y[31] → 0.30231, y[32] → 0.320914, y[33] → 0.330405, y[34] → 0.339626, y[35] → 0.358274,

y[36] → 0.358483, y[37] → 0.377052, y[38] → 0.387058, y[39] → 0.395916, y[40] → 0.414801,

y[41] → 0.415918, y[42] → 0.433746, y[43] → 0.444957, y[44] → 0.452709, y[45] → 0.452729,

y[46] → 0.469703, y[47] → 0.469703, y[48] → 0.486067, y[49] → 0.491174, y[50] → 0.50248,

y[51] → 0.512752, y[52] → 0.514148, y[53] → 0.53517, y[54] → 0.557677, y[55] → 0.572908,

y[56] → 0.572908, y[57] → 0.580202, y[58] → 0.580203, y[59] → 0.590236, y[60] → 0.590236,

y[61] → 0.599733, y[62] → 0.606913, y[63] → 0.606913, y[64] → 0.619388, y[65] → 0.62365,

y[66] → 0.62365, y[67] → 0.62641, y[68] → 0.639969, y[69] → 0.641076, y[70] → 0.641076,

y[71] → 0.658757, y[72] → 0.658757, y[73] → 0.661207, y[74] → 0.676529, y[75] → 0.676529,

y[76] → 0.682633, y[77] → 0.682703, y[78] → 0.693751, y[79] → 0.693751, y[80] → 0.696192,

y[81] → 0.702509, y[82] → 0.710502, y[83] → 0.710502, y[84] → 0.72727, y[85] → 0.72727,

y[86] → 0.729863, y[87] → 0.744284, y[88] → 0.744284, y[89] → 0.761388, y[90] → 0.761388,

y[91] → 0.778505, y[92] → 0.778505, y[93] → 0.779846, y[94] → 0.79564, y[95] → 0.79564,

y[96] → 0.812759, y[97] → 0.812759, y[98] → 0.82989, y[99] → 0.82989, y[100] → 0.830164,

y[101] → 0.847153, y[102] → 0.847153, y[103] → 0.864383, y[104] → 0.864383,

y[105] → 0.881611, y[106] → 0.881611, y[107] → 0.881611, y[108] → 0.896189,

y[109] → 0.896189, y[110] → 0.914287, y[111] → 0.954891, y[112] → 0.956724,

y[113] → 0.983306, y[114] → 0.98579, y[115] → 0.975192, y[116] → 0.969786,

y[117] → 0.96303, y[118] → 0.959146, y[119] → 0.955629, y[120] → 0.9437,

y[121] → 0.9437, y[122] → 0.941378, y[123] → 0.935937, y[124] → 0.926839,

y[125] → 0.926839, y[126] → 0.920206, y[127] → 0.918537, y[128] → 0.905843,

y[129] → 0.900841, y[130] → 0.892116, y[131] → 0.88565, y[132] → 0.883203,

y[133] → 0.866072, y[134] → 0.865649, y[135] → 0.865622, y[136] → 0.847871,

y[137] → 0.845035, y[138] → 0.838788, y[139] → 0.836471, y[140] → 0.830137,

y[141] → 0.824169, y[142] → 0.822596, y[143] → 0.812426, y[144] → 0.807941,

y[145] → 0.806457, y[146] → 0.803587, y[147] → 0.794751, y[148] → 0.79041,

y[149] → 0.783366, y[150] → 0.780355, y[151] → 0.777105, y[152] → 0.774421,

y[153] → 0.766206, y[154] → 0.763472, y[155] → 0.759495, y[156] → 0.758511,

y[157] → 0.75364, y[158] → 0.743867, y[159] → 0.742655, y[160] → 0.741912,

y[161] → 0.734506, y[162] → 0.727382, y[163] → 0.726819, y[164] → 0.72434,

y[165] → 0.724335, y[166] → 0.710884, y[167] → 0.706712, y[168] → 0.704414,

y[169] → 0.703033, y[170] → 0.700426, y[171] → 0.694913, y[172] → 0.689078,

y[173] → 0.678926, y[174] → 0.673191, y[175] → 0.67183, y[176] → 0.671438,

y[177] → 0.662918, y[178] → 0.653785, y[179] → 0.646893, y[180] → 0.645702,

y[181] → 0.642829, y[182] → 0.640777, y[183] → 0.636125, y[184] → 0.630842,

y[185] → 0.626931, y[186] → 0.618455, y[187] → 0.61795, y[188] → 0.614773,

y[189] → 0.61101, y[190] → 0.609881, y[191] → 0.598682, y[192] → 0.595056,

y[193] → 0.58996, y[194] → 0.582586, y[195] → 0.579106, y[196] → 0.579079,

y[197] → 0.566303, y[198] → 0.562772, y[199] → 0.560829, y[200] → 0.549967,

y[201] → 0.549133, y[202] → 0.546437, y[203] → 0.533624, y[204] → 0.530935,
y[205] → 0.530096, y[206] → 0.519197, y[207] → 0.517284, y[208] → 0.51375,
y[209] → 0.500944, y[210] → 0.50092, y[211] → 0.49715, y[212] → 0.490086,
y[213] → 0.483916, y[214] → 0.479735, y[215] → 0.466875, y[216] → 0.464627,
y[217] → 0.462246, y[218] → 0.462232, y[219] → 0.457446, y[220] → 0.449772,
y[221] → 0.44464, y[222] → 0.437634, y[223] → 0.434646, y[224] → 0.432616,
y[225] → 0.426962, y[226] → 0.426949, y[227] → 0.417129, y[228] → 0.408537,
y[229] → 0.407896, y[230] → 0.396048, y[231] → 0.390101, y[232] → 0.383233,
y[233] → 0.381455, y[234] → 0.374946, y[235] → 0.371694, y[236] → 0.371689,
y[237] → 0.366924, y[238] → 0.366923, y[239] → 0.355559, y[240] → 0.35272,
y[241] → 0.352508, y[242] → 0.351149, y[243] → 0.334001, y[244] → 0.333269,
y[245] → 0.330629, y[246] → 0.330138, y[247] → 0.33012, y[248] → 0.313497,
y[249] → 0.310241, y[250] → 0.306543, y[251] → 0.30528, y[252] → 0.29361,
y[253] → 0.289869, y[254] → 0.289735, y[255] → 0.289735, y[256] → 0.282697,
y[257] → 0.280595, y[258] → 0.273955, y[259] → 0.273486, y[260] → 0.269565,
y[261] → 0.258362, y[262] → 0.256148, y[263] → 0.253212, y[264] → 0.249277,
y[265] → 0.23399, y[266] → 0.232953, y[267] → 0.231708, y[268] → 0.228984,
y[269] → 0.223573, y[270] → 0.223573, y[271] → 0.212931, y[272] → 0.210218,
y[273] → 0.210072, y[274] → 0.208631, y[275] → 0.207037, y[276] → 0.192725,
y[277] → 0.188345, y[278] → 0.185668, y[279] → 0.182635, y[280] → 0.168064,
y[281] → 0.161238, y[282] → 0.158219, y[283] → 0.156907, y[284] → 0.156907,
y[285] → 0.14778, y[286] → 0.1368, y[287] → 0.133821, y[288] → 0.127499,
y[289] → 0.12517, y[290] → 0.112347, y[291] → 0.109425, y[292] → 0.107219,
y[293] → 0.106524, y[294] → 0.0916548, y[295] → 0.0916545, y[296] → 0.0879295,
y[297] → 0.0878836, y[298] → 0.0869322, y[299] → 0.0850037, y[300] → 0.0692433,
y[301] → 0.0666466, y[302] → 0.0635256, y[303] → 0.0506013, y[304] → 0.0463632,
y[305] → 0.0390785, y[306] → 0.0319585, y[307] → 0.0260759, y[308] → 0.0260554,
y[309] → 0.026055, y[310] → 0.0131144, y[311] → 0.0131145, y[312] → 0.00631395}}

```
In[*]:= Show[Plot[Piecewise[
        {{Tan[π / 2 - ϕ] * x, 0 ≤ x ≤ xx - Cos[ϕ]},
         {√(1 - (x - xx) ^ 2), Cos[π - ϕ] + xx ≤ x ≤ Cos[π - (ArcTan[1 / (xx / 2)] - ψ)] + xx},
         {Tan[π / 2 - (ArcTan[1 / (xx / 2)] - ψ)] * (x - (xx - Cos[ArcTan[1 / (xx / 2)] - ψ])) +
           Sin[ArcTan[1 / (xx / 2)] - ψ], xx - (Cos[ArcTan[1 / (xx / 2)] - ψ]) ≤ x ≤ xx / 2},
         {-1/Tan[ArcTan[1 / (xx / 2)] - ψ] * (x - Cos[ArcTan[1 / (xx / 2)] - ψ]) +
           Sin[ArcTan[1 / (xx / 2)] - ψ], xx / 2 ≤ x ≤ Cos[ArcTan[1 / (xx / 2)] - ψ]},
         {√(1 - x ^ 2), Cos[ArcTan[1 / (xx / 2)] - ψ] ≤ x ≤ Cos[ϕ]},
         {-1/Tan[ϕ] * (x - Cos[ϕ]) + Sin[ϕ], Cos[ϕ] ≤ x ≤ xx}}], {x, 0, xx},
       PerformanceGoal → "Speed", PlotStyle → Black, PlotLegends → {"Zalgaller's"},
       PlotRange → All, AspectRatio → Automatic, PerformanceGoal → "Speed", AxesOrigin → {0, 0}],
      ListLinePlot[Prepend[Table[{x[i], y[i]} /. res[[2]], {i, 1, nn * mm}], {0, 0}],
       PlotRange → {{-0.2, 1.2}, {-0.2, 1.2}}, AspectRatio → Automatic,
       PerformanceGoal → "Speed", PlotStyle → Orange, PlotLegends → {"N=12,M=26"}]]
```

Out[*]=

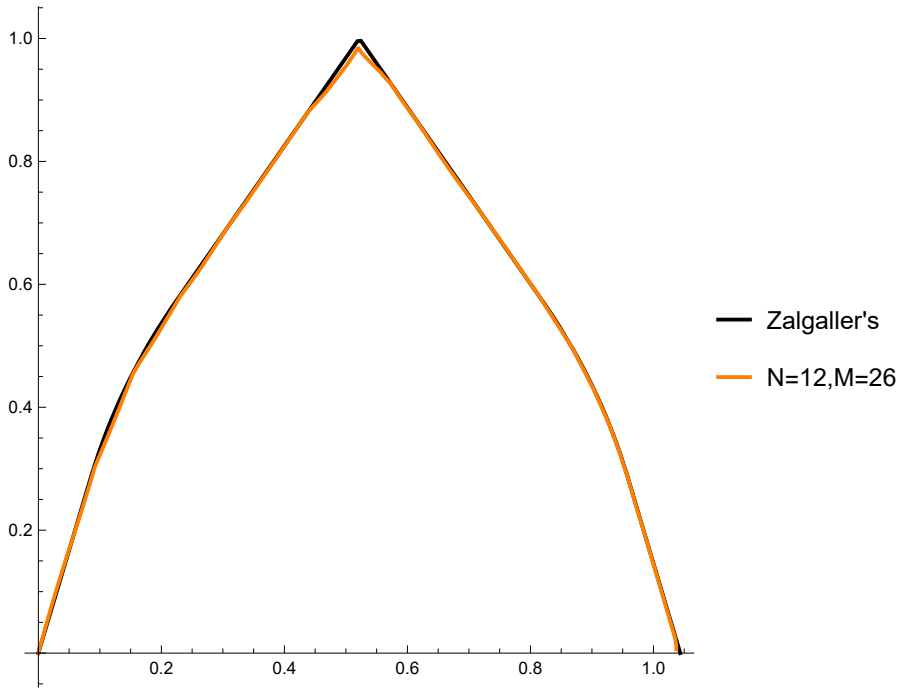

# Appendix 4

(*Section 8.5*)

*In[ ]:=* (*opaque set with one curve for circle and ellipse*)

*In[ ]:=* **n = 100;** (*Number of discretization points*)

*In[ ]:=* (*Numerical minimization for global optimal*)

*In[ ]:=* **sol = NMinimize** $\Big[\Big\{$ **Sum** $\Big[\sqrt{(x[i] - x[i-1])^2 + (y[i] - y[i-1])^2}, \{i, 1, n\}\Big]$ (*objective*),

**Table** $\Big[\Big(x[i] * \text{Cos}\Big[\frac{2\pi}{n} * i\Big] + y[i] * \text{Sin}\Big[\frac{2\pi}{n} * i\Big] - 1\Big) == 0, \{i, 0, n\}\Big]$ (*constraints*) $\Big\}$,

**Flatten[Table[{x[i], y[i]}, {i, 0, n}]]** $\Big];$

*In[ ]:=* (*length of path*)

*In[ ]:=* **sol[[1]]**

*Out[ ]=*
 5.14105

*In[ ]:=* **RealDigits[%]**

*Out[ ]=*
 {{5, 1, 4, 1, 0, 4, 8, 1, 2, 3, 1, 9, 1, 1, 4, 6}, 1}

*In[ ]:=* **Show[ContourPlot[x^2 + y^2 - 1 == 0, {x, -1.2, 1.2}, {y, -1.2, 1.2},**
  **ContourStyle → Red, PerformanceGoal → "Speed", PlotPoints → 100],**
  **ListLinePlot[Table[{x[i], y[i]}, {i, 0, n}] /. sol[[2]], AspectRatio → Automatic,**
  **PerformanceGoal → "Speed", PlotRange → All, PlotStyle → Black],**
  **PlotRange → All, ImageSize → 300, AspectRatio → Automatic]**

*Out[ ]=*

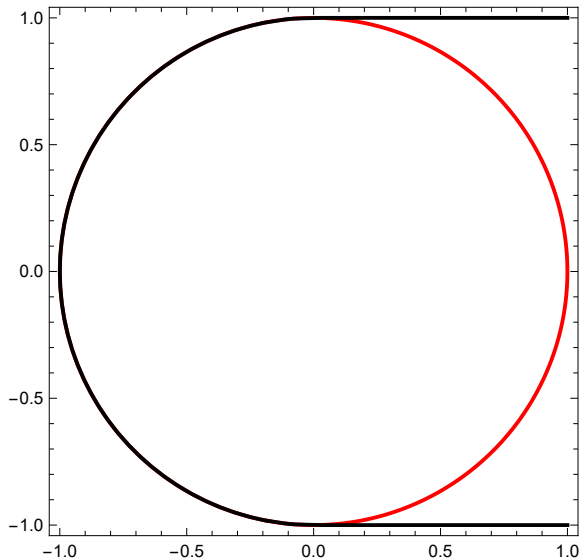

*In[ ]:=* **n = 100;** (*Number of discretization points*)

(*opening direction π/2*)

*In[ ]:=*  (*Numerical minimization for global optimal*)

*In[ ]:=*  **sol = NMinimize[{Sum[ $\sqrt{(x[i] - x[i-1])^2 + (y[i] - y[i-1])^2}$ , {i, 1, n}]** (*objective*),

$\quad$ **Table[$\left(x[i] * Cos\left[\dfrac{2\pi}{n} * i + \pi/2\right] + 2 * y[i] * Sin\left[\dfrac{2\pi}{n} * i + \pi/2\right] - 1\right) == 0$, {i, 0, n}]**

$\quad$ (*constraints*)**}, Flatten[Table[{x[i], y[i]}, {i, 0, n}]]];**

*In[ ]:=*  (*length of path*)

*In[ ]:=*  **sol[[1]]**

*Out[ ]=*

$\quad$ 3.42165

*In[ ]:=*  **RealDigits[%]**

*Out[ ]=*

$\quad$ {{3, 4, 2, 1, 6, 5, 4, 2, 0, 9, 9, 8, 5, 5, 0, 4}, 1}

*In[ ]:=*  **Show[ContourPlot[$x^2 + 4 y^2 - 1 == 0$, {x, -1.2, 1.2}, {y, -1.2, 1.2},**

$\quad$ **ContourStyle → Red, PerformanceGoal → "Speed", PlotPoints → 100],**

$\quad$ **ListLinePlot[Table[{x[i], y[i]}, {i, 0, n}] /. sol[[2]], AspectRatio → Automatic,**

$\quad$ **PerformanceGoal → "Speed", PlotRange → All, PlotStyle → Black],**

$\quad$ **PlotRange → {{-1.2, 1.2}, {-1.2, 1.2}}, ImageSize → 300, AspectRatio → Automatic]**

*Out[ ]=*

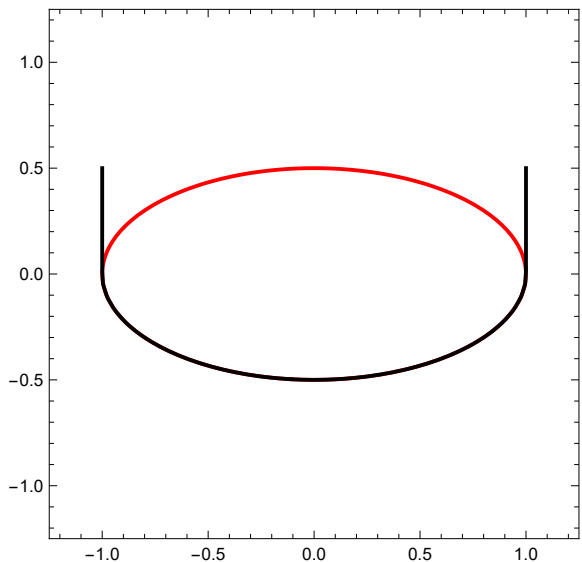

(*different opening directions π/4*)

*In[ ]:=*  **sol = NMinimize[{Sum[ $\sqrt{(x[i] - x[i-1])^2 + (y[i] - y[i-1])^2}$ , {i, 1, n}]** (*objective*),

$\quad$ **Table[$\left(x[i] * Cos\left[\dfrac{2\pi}{n} * i + \pi/4\right] + 2 * y[i] * Sin\left[\dfrac{2\pi}{n} * i + \pi/4\right] - 1\right) == 0$, {i, 0, n}]**

$\quad$ (*constraints*)**}, Flatten[Table[{x[i], y[i]}, {i, 0, n}]]];**

*In[*]:=* **sol[[1]]**

*Out[*]=*
3.68661

*In[*]:=* **Show[ContourPlot[x^2 + 4 y^2 - 1 == 0, {x, -1.2, 1.2}, {y, -1.2, 1.2},**
**ContourStyle → Red, PerformanceGoal → "Speed", PlotPoints → 100],**
**ListLinePlot[Table[{x[i], y[i]}, {i, 0, n}] /. sol[[2]], AspectRatio → Automatic,**
**PerformanceGoal → "Speed", PlotRange → All, PlotStyle → Black],**
**PlotRange → {{-1.2, 1.2}, {-1.2, 1.2}}, ImageSize → 300, AspectRatio → Automatic]**

*Out[*]=*

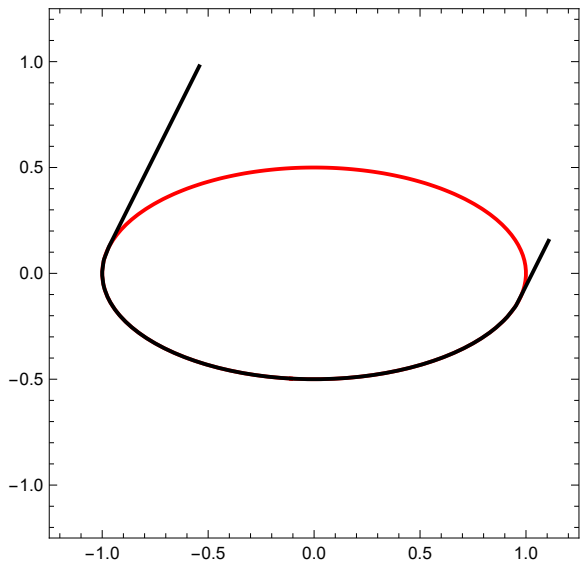

(*different opening directions θ*)

*In[*]:=* **sol = NMinimize[{Sum[$\sqrt{(x[i] - x[i-1])^2 + (y[i] - y[i-1])^2}$, {i, 1, n}] (*objective*),**

**Table[$\left(x[i] * Cos\left[\frac{2\pi}{n} * i + \theta\right] + 2 * y[i] * Sin\left[\frac{2\pi}{n} * i + \theta\right] - 1\right)$ == 0, {i, 0, n}]**

**(*constraints*)}, Flatten[Table[{x[i], y[i]}, {i, 0, n}]]];**

*In[*]:=* **sol[[1]]**

*Out[*]=*
4.42171

*In[ ]:=* `Show[ContourPlot[x^2 + 4 y^2 - 1 == 0, {x, -1.2, 1.2}, {y, -1.2, 1.2},`
`ContourStyle → Red, PerformanceGoal → "Speed", PlotPoints → 100],`
`ListLinePlot[Table[{x[i], y[i]}, {i, 0, n}] /. sol[[2]], AspectRatio → Automatic,`
`PerformanceGoal → "Speed", PlotRange → All, PlotStyle → Black],`
`PlotRange → {{-1.2, 1.2}, {-1.2, 1.2}}, ImageSize → 300, AspectRatio → Automatic]`

*Out[ ]=*

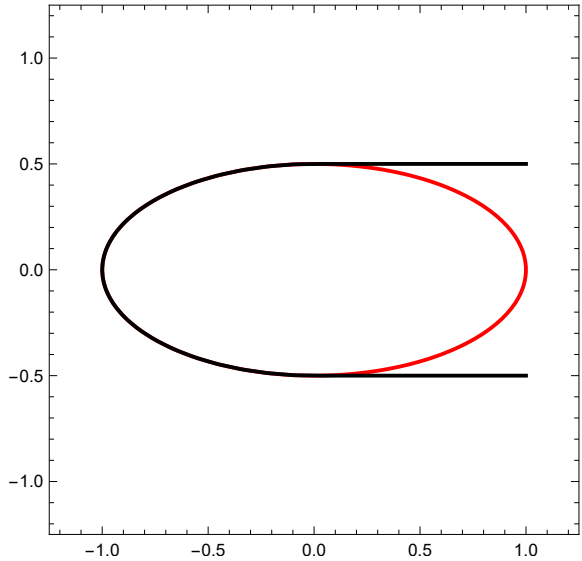

# Appendix 5

```
(*Section 9.1 Lost-in-a-Forest Problem with closed path return to starting point*)

(*discrete formulation*)

(*search a line with unit distance*)
```

*In[ ]:=* `n = 500;` *(*Number of discretization points*)*

*In[ ]:=* *(*Numerical minimization for global optimal*)*

*In[ ]:=* `sol = NMinimize[`

$$\left\{ \sqrt{x[0]^2 + y[0]^2} + \text{Sum}\left[ \sqrt{(x[i] - x[i-1])^2 + (y[i] - y[i-1])^2}, \{i, 1, n\} \right] + \sqrt{x[n]^2 + y[n]^2} \right.$$

*(*objective*)*, `Table[` $\left( x[i] * \text{Cos}\left[\frac{2\pi}{n} * i\right] + y[i] * \text{Sin}\left[\frac{2\pi}{n} * i\right] - 1 \right)$ == 0, {i, 0, n}]

*(*constraints*)* $\left.\right\}$, `Flatten[Table[{x[i], y[i]}, {i, 0, n}]]];`

*In[ ]:=* *(*length of path*)*

*In[ ]:=* `sol[[1]]`

*Out[ ]=*

7.65286

*In[ ]:=* `RealDigits[%]`

*Out[ ]=*

{{7, 6, 5, 2, 8, 6, 4, 2, 7, 6, 7, 3, 0, 7, 2, 8}, 1}

*In[ ]:=* **Show[**
    **ListLinePlot[Append[Prepend[Table[{x[i], y[i]}, {i, 0, n}] /. sol〚2〛, {0, 0}], {0, 0}],**
      **AspectRatio → Automatic, PerformanceGoal → "Speed", PlotRange → All, PlotStyle → Black],**
    **ContourPlot[x - 1 ⩵ 0, {x, -1.2, 1.2}, {y, -1.2, 1.2}, ContourStyle → Red,**
      **PerformanceGoal → "Speed"], PlotRange → All, ImageSize → 300]**

*Out[ ]=*

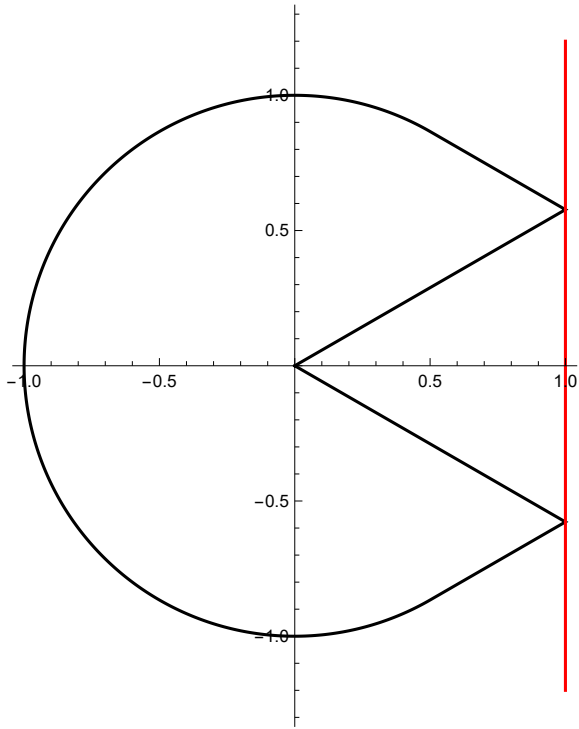

*In[ ]:=* **NotebookSave[];**

(*search a circle from exterior*)

(*r/s=1*)

*In[ ]:=* **n = 50;** (*Number of discretization points*)

*In[ ]:=* (*Numerical minimization for global optimal*)

*In[ ]:=* **sol = NMinimize[{ {√(x[0]² + y[0]²) +**
      **Sum[√((x[i] - x[i-1])² + (y[i] - y[i-1])²), {i, 1, n}] + √(x[n]² + y[n]²),**
    **Table[(x[i] - Cos[ (2π/n) * i])^2 + (y[i] - Sin[ (2π/n) * i])^2 ⩵ 0.5^2, {i, 0, n}]},**
    **Flatten[Table[{x[i], y[i]}, {i, 0, n}]]}];**

*In[ ]:=* (*length*)

*In[ ]:=* **sol〚1〛**

*Out[ ]=*
    3.97269

*In[ ]:=* **RealDigits[%]**

*Out[ ]=*
{{3, 9, 7, 2, 6, 8, 8, 6, 0, 1, 8, 1, 1, 8, 9}, 1}

*In[ ]:=* **Show[**
**ListLinePlot[Append[Prepend[Table[{x[i], y[i]}, {i, 0, n}] /. sol〚2〛, {0, 0}], {0, 0}],**
**AspectRatio → Automatic, PerformanceGoal → "Speed", PlotRange → All, PlotStyle → Black],**
**ContourPlot[(x - 1) ^2 + (y - 0) ^2 ⩵ 0.5^2, {x, -1, 2}, {y, -1.5, 1.5}, ContourStyle → Red],**
**PlotRange → All, ImageSize → 300]**

*Out[ ]=*

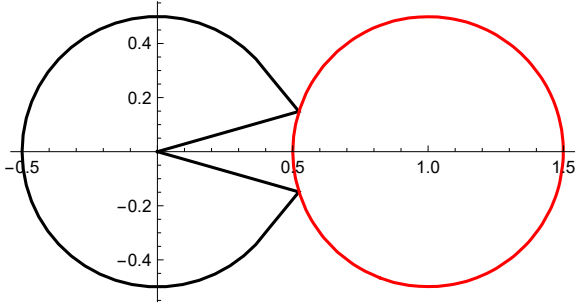

*In[ ]:=* **NotebookSave[];**

(*search a circle from interior*)

*In[ ]:=* (*r/s=0.2*)

*In[ ]:=* **n = 50;** (*Number of discretization points*)

*In[ ]:=* (*Numerical minimization for global optimal*)

*In[ ]:=* **sol = NMinimize[{{$\sqrt{x[0]^2 + y[0]^2}$ +**
**Sum[$\sqrt{(x[i] - x[i-1])^2 + (y[i] - y[i-1])^2}$, {i, 1, n}] + $\sqrt{x[n]^2 + y[n]^2}$,**
**Table[$\left(x[i] - Cos\left[\frac{2\pi}{n} * i\right]\right)$^2 + $\left(y[i] - Sin\left[\frac{2\pi}{n} * i\right]\right)$^2 ⩵ 1.2^2, {i, 0, n}]},**
**Flatten[Table[{x[i], y[i]}, {i, 0, n}]]];**

*In[ ]:=* (*length*)

*In[ ]:=* **sol〚1〛**

*Out[ ]=*
1.50702

*In[ ]:=* **RealDigits[%]**

*Out[ ]=*
{{1, 5, 0, 7, 0, 1, 7, 2, 4, 6, 6, 2, 8, 7, 6, 7}, 1}

```
In[*]:= Show[
    ListLinePlot[Append[Prepend[Table[{x[i], y[i]}, {i, 0, n}] /. sol[[2]], {0, 0}], {0, 0}],
     AspectRatio → Automatic, PerformanceGoal → "Speed", PlotRange → All, PlotStyle → Black],
    ContourPlot[(x - 1)^2 + (y - 0)^2 == 1.2^2, {x, -1, 3}, {y, -1.5, 1.5}, ContourStyle → Red],
    PlotRange → All, ImageSize → 300]
```

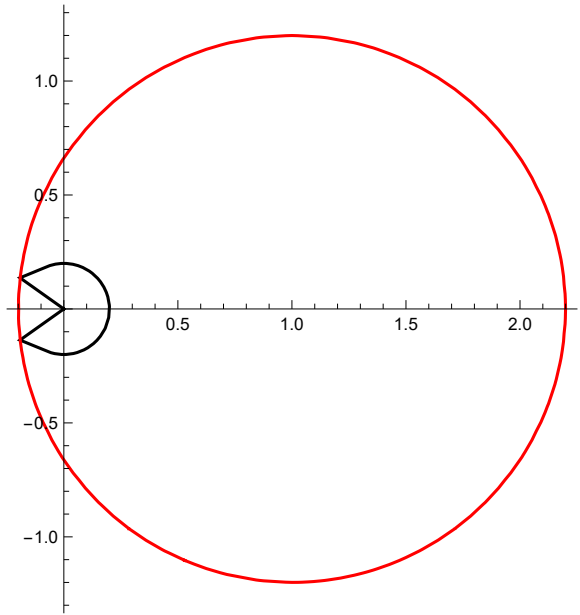

```
In[*]:= NotebookSave[];
```

```
(*search for one point with given distance 1*)
```

```
In[*]:= nn = 50; (*Number of discretization points*)
```

```
In[*]:= (*Numerical minimization for global optimal*)
```

```
In[*]:= sol = NMinimize[{{√(x[0]^2 + y[0]^2) +
       Sum[√((x[i] - x[i-1])^2 + (y[i] - y[i-1])^2), {i, 1, nn}] + √(x[nn]^2 + y[nn]^2),
      Table[(x[i] - Cos[2π/nn * i])^2 + (y[i] - Sin[2π/nn * i])^2 == 0.0^2, {i, 0, nn}]},
     Flatten[Table[{x[i], y[i]}, {i, 0, nn}]]];
```

```
In[*]:= (*length*)
```

```
In[*]:= sol[[1]]
```

```
Out[*]=
    8.27905
```

```
In[*]:= RealDigits[%]
```

```
Out[*]=
    {{8, 2, 7, 9, 0, 5, 1, 9, 3, 2, 8, 6, 2, 9, 9, 8}, 1}
```

*In[ ]:=* (*the length should be 2π+1*)

*In[ ]:=* **Show[**
  **ListLinePlot[Append[Prepend[Table[{x[i], y[i]}, {i, 0, nn}] /. sol〚2〛, {0, 0}], {0, 0}],**
   **AspectRatio → Automatic, PerformanceGoal → "Speed", PlotRange → All, PlotStyle → Black],**
  **Graphics[{Red, Disk[{1, 0}, 0.02]}], PlotRange → All, ImageSize → 300]**

*Out[ ]=*

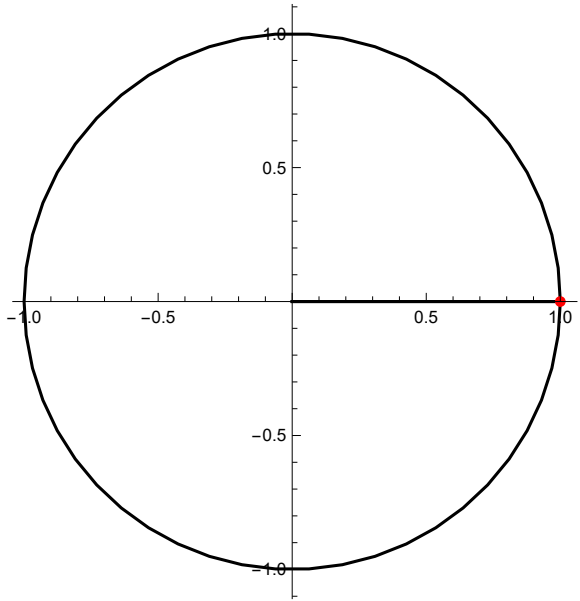

*In[ ]:=* **NotebookSave[];**

---

(*two perpendicular lines with 0.5 distance*)

*In[ ]:=* **nn = 200;** (*Number of discretization points*)

*In[ ]:=* (*Numerical minimization for global optimal*)

*In[ ]:=* **sol = NMinimize[{{ $\sqrt{x[0]^2 + y[0]^2}$ +**
    **Sum[ $\sqrt{(x[i] - x[i-1])^2 + (y[i] - y[i-1])^2}$ , {i, 1, nn}] + $\sqrt{x[nn]^2 + y[nn]^2}$ ,**
    **Table[ $\left[-x[i] * Sin\left[\frac{2\pi}{nn} * \frac{3}{4} * i\right] + y[i] * Cos\left[\frac{2\pi}{nn} * \frac{3}{4} * i\right] + 0.5\right] == 0$ , {i, 0, nn}]},**
    **Flatten[Table[{x[i], y[i]}, {i, 0, nn}]]];**

*In[ ]:=* (*length*)

*In[ ]:=* **sol〚1〛**

*Out[ ]=*
  3.04102

*In[ ]:=* **Show[**
    **ListLinePlot[Append[Prepend[Table[{x[i], y[i]}, {i, 0, nn}] /. sol[[2]], {0, 0}], {0, 0}],**
     **AspectRatio → Automatic, PerformanceGoal → "Speed", PlotRange → All, PlotStyle → Black],**
    **ContourPlot[x + 0.5 == 0, {x, -1.5, 1.5}, {y, -0.5, 1.5}, ContourStyle → Red],**
    **ContourPlot[y + 0.5 == 0, {x, -0.5, 1.5}, {y, -1, 1.5}, ContourStyle → Red],**
    **PlotRange → {{-1, 1.5}, {-1, 1.5}}, ImageSize → 300]**

*Out[ ]=*

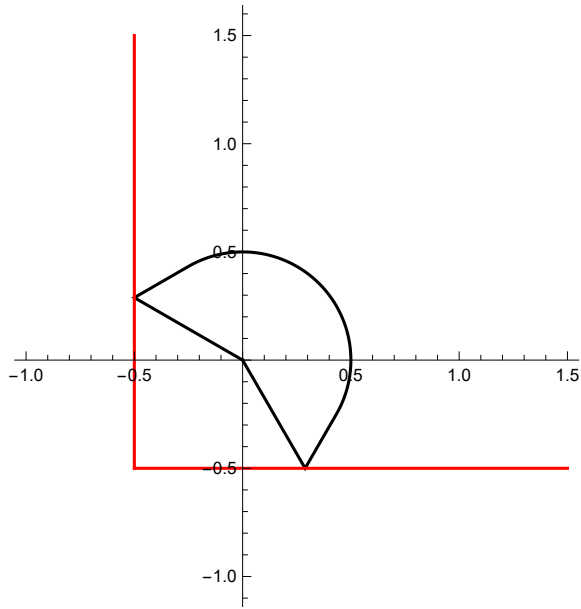

*In[ ]:=* **NotebookSave[];**

    (*two parallel lines from the middle*)

*In[ ]:=* **nn = 200;** (*Number of discretization points*)

*In[ ]:=* (*Numerical minimization for global optimal*)

*In[ ]:=* **sol = NMinimize[{ { $\sqrt{x[0]^2 + y[0]^2}$ +**
    **Sum[ $\sqrt{(x[i] - x[i-1])^2 + (y[i] - y[i-1])^2}$ , {i, 1, nn / 2}] + $\sqrt{x[nn/2]^2 + y[nn/2]^2}$ ,**
    **Table[ $x[i] * Cos\left[\dfrac{2\pi}{nn} * i - \pi/2\right] + y[i] * Sin\left[\dfrac{2\pi}{nn} * i - \pi/2\right]$ == 0.5, {i, 0, nn / 2}]},**
    **Flatten[Table[{x[i], y[i]}, {i, 0, nn / 2}]]}];**

*In[ ]:=* (*length*)

*In[ ]:=* **sol[[1]]**

*Out[ ]=*
    2.25563

*In[ ]:=* **Show[ListLinePlot[**
    **Append[Prepend[Table[{x[i], y[i]}, {i, 0, nn / 2}] /. sol〚2〛, {0, 0}], {0, 0}],**
    **AspectRatio → Automatic, PerformanceGoal → "Speed", PlotRange → All, PlotStyle → Black],**
  **ContourPlot[y - 0.5 ⩵ 0, {x, -0.5, 1}, {y, -1.5, 1.5}, ContourStyle → Red],**
  **ContourPlot[y + 0.5 ⩵ 0, {x, -0.5, 1}, {y, -1.5, 1.5}, ContourStyle → Red, PlotRange → All]**

*Out[ ]=*

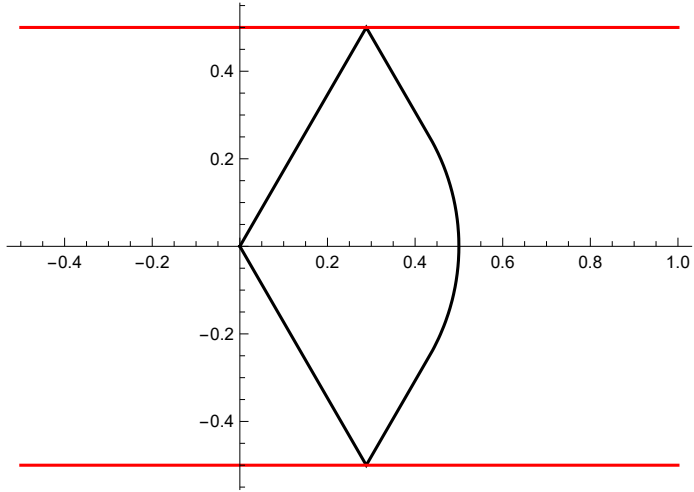

*In[ ]:=* **NotebookSave[];**

(*search for two lines from the angle bisector in the middle*)

*In[ ]:=* (*π/6 - 30 degree*)

*In[ ]:=* **BB = 1 / 2 + π / 6 / (2 π);**

*In[ ]:=* **nn = 300;** (*Number of discretization points*)

*In[ ]:=* (*Numerical minimization for global optimal*)

*In[ ]:=* **sol = NMinimize[{$\sqrt{x[0]^2 + y[0]^2}$ +**
    **Sum[$\sqrt{(x[i] - x[i-1])^2 + (y[i] - y[i-1])^2}$, {i, 1, BB * nn}] + $\sqrt{x[BB * nn]^2 + y[BB * nn]^2}$,**
    **Table[$\left(x[i] * Cos\left[\frac{2\pi}{nn} * i\right] + y[i] * Sin\left[\frac{2\pi}{nn} * i\right] + 0.5\right)$ ⩵ 0, {i, 0, BB * nn}]},**
    **Flatten[Table[{x[i], y[i]}, {i, 0, BB * nn}]]];**

*In[ ]:=* (*length*)

*In[ ]:=* **sol〚1〛**

*Out[ ]=*
  2.51743

```
In[ ]:= Show[ListLinePlot[
         Append[Prepend[Table[{x[i], y[i]}, {i, 0, BB * nn}] /. sol[[2]], {0, 0}], {0, 0}],
         AspectRatio → Automatic, PerformanceGoal → "Speed", PlotRange → All, PlotStyle → Black],
       ContourPlot[{x * Cos[0] + y * Sin[0] + 0.5 == 0, x * Cos[0 + π / 6] + y * Sin[0 + π / 6] - 0.5 == 0},
        {x, -2, 2}, {y, -2, 2}, ContourStyle → Red],
       PlotRange → {{-1.5, 1.5}, {-1.5, 1.5}}, ImageSize → 300]
```

Out[ ]=

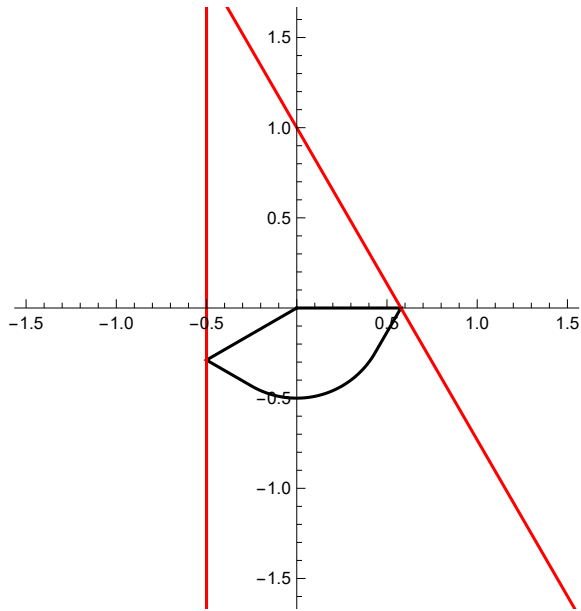

```
In[ ]:= (*π/3 - 60 degree*)

In[ ]:= BB = 1 / 2 + π / 3 / (2 π);

In[ ]:= nn = 300; (*Number of discretization points*)

In[ ]:= sol = NMinimize[{{√(x[0]^2 + y[0]^2) +
         Sum[√((x[i] - x[i - 1])^2 + (y[i] - y[i - 1])^2), {i, 1, BB * nn}] + √(x[BB * nn]^2 + y[BB * nn]^2),
        Table[(x[i] * Cos[(2 π)/nn * i] + y[i] * Sin[(2 π)/nn * i] + 0.5) == 0, {i, 0, BB * nn}]},
        Flatten[Table[{x[i], y[i]}, {i, 0, BB * nn}]]];

In[ ]:= (*length*)

In[ ]:= sol[[1]]
```

Out[ ]=
   2.77923

*In[ ]:=* `Show[ListLinePlot[`
        `Append[Prepend[Table[{x[i], y[i]}, {i, 0, BB * nn}] /. sol[[2]], {0, 0}], {0, 0}],`
        `AspectRatio → Automatic, PerformanceGoal → "Speed", PlotRange → All, PlotStyle → Black],`
      `ContourPlot[{x * Cos[0] + y * Sin[0] + 0.5 ⩵ 0, x * Cos[0 + π / 3] + y * Sin[0 + π / 3] - 0.5 ⩵ 0},`
        `{x, -2, 2}, {y, -2, 2}, ContourStyle → Red],`
      `PlotRange → {{-1.5, 1.5}, {-1.5, 1.5}}, ImageSize → 300]`

*Out[ ]:=*

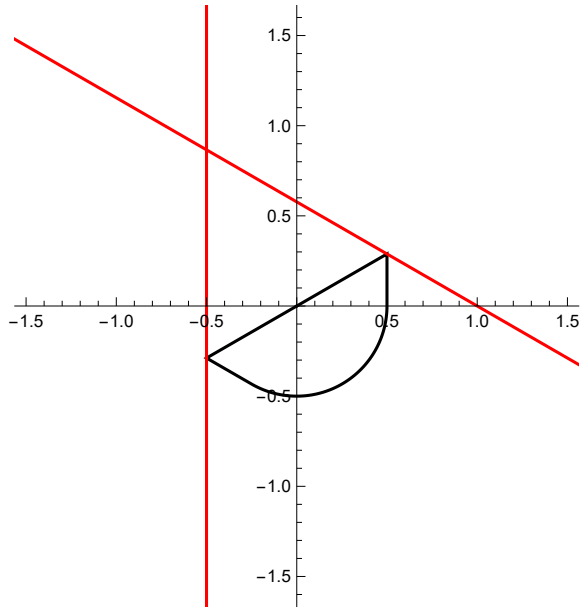

*In[ ]:=* `(*2π/3 - 120 degree*)`

*In[ ]:=* `BB = 1 / 2 + 2 π / 3 / (2 π);`

*In[ ]:=* `nn = 300;` `(*Number of discretization points*)`

*In[ ]:=* `sol = NMinimize[{{ √(x[0]² + y[0]²) +`
        `Sum[ √((x[i] - x[i - 1])² + (y[i] - y[i - 1])²) , {i, 1, BB * nn}] + √(x[BB * nn]² + y[BB * nn]²) ,`
        `Table[(x[i] * Cos[ (2 π / nn) * i] + y[i] * Sin[ (2 π / nn) * i] + 0.5) ⩵ 0, {i, 0, BB * nn}]},`
        `Flatten[Table[{x[i], y[i]}, {i, 0, BB * nn}]]];`

*In[ ]:=* `(*length*)`

*In[ ]:=* `sol[[1]]`

*Out[ ]:=*
      `3.30282`

*In[ ]:=* ```Show[ListLinePlot[
    Append[Prepend[Table[{x[i], y[i]}, {i, 0, BB * nn}]] /. sol[[2]], {0, 0}], {0, 0}],
    AspectRatio → Automatic, PerformanceGoal → "Speed", PlotRange → All, PlotStyle → Black],
  ContourPlot[{x * Cos[0] + y * Sin[0] + 0.5 ⩵ 0, x * Cos[0 + 2 π / 3] + y * Sin[0 + 2 π / 3] - 0.5 ⩵ 0},
    {x, -2, 2}, {y, -2, 2}, ContourStyle → Red],
  PlotRange → {{-1.5, 1.5}, {-1.5, 1.5}}, ImageSize → 300]```

*Out[ ]:=*

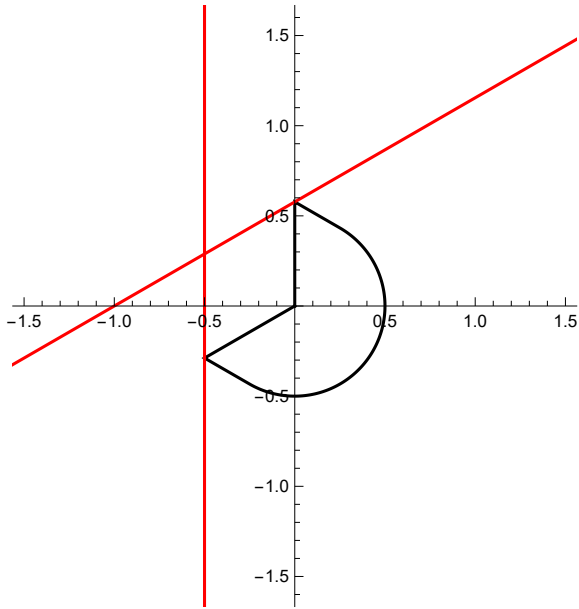

*In[ ]:=* `(*5π/6 - 150 degree*)`

*In[ ]:=* `BB = 1 / 2 + 5 π / 6 / (2 π);`

*In[ ]:=* `nn = 300;` `(*Number of discretization points*)`

*In[ ]:=* ```sol = NMinimize[{{ √(x[0]² + y[0]²) +
    Sum[ √((x[i] - x[i - 1])² + (y[i] - y[i - 1])²), {i, 1, BB * nn}] + √(x[BB * nn]² + y[BB * nn]²),
    Table[(x[i] * Cos[(2 π)/nn * i] + y[i] * Sin[(2 π)/nn * i] + 0.5) ⩵ 0, {i, 0, BB * nn}]},
    Flatten[Table[{x[i], y[i]}, {i, 0, BB * nn}]]];```

*In[ ]:=* `(*length*)`

*In[ ]:=* `sol[[1]]`

*Out[ ]:=*

`3.56461`

*In[ ]:=* `Show[ListLinePlot[`
`    Append[Prepend[Table[{x[i], y[i]}, {i, 0, BB * nn}] /. sol[[2]], {0, 0}], {0, 0}],`
`    AspectRatio → Automatic, PerformanceGoal → "Speed", PlotRange → All, PlotStyle → Black],`
`  ContourPlot[{x * Cos[0] + y * Sin[0] + 0.5 ⩵ 0, x * Cos[0 + 5 π / 6] + y * Sin[0 + 5 π / 6] - 0.5 ⩵ 0},`
`   {x, -2, 2}, {y, -2, 2}, ContourStyle → Red],`
`  PlotRange → {{-1.5, 1.5}, {-1.5, 1.5}}, ImageSize → 300]`

*Out[ ]=*

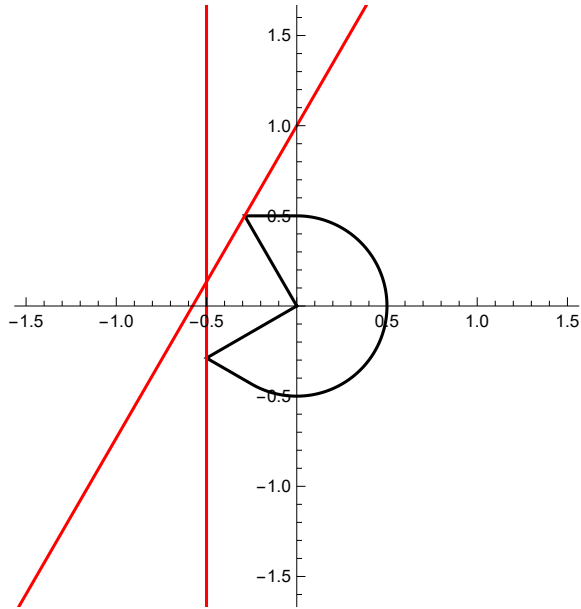

*In[ ]:=* `NotebookSave[];`



\end{document}